\input amstex
%

\def\next{AMS-SEKR}\ifx\styname\next \endinput\fi
\catcode`\@=11
\def\styname{AMS-SEKR}
\def\styversion{2.0}
{\W@{}\W@{\styname.STY - Version \styversion}\W@{}}
\hyphenation{acad-e-my acad-e-mies af-ter-thought anom-aly anom-alies
an-ti-deriv-a-tive an-tin-o-my an-tin-o-mies apoth-e-o-ses apoth-e-o-sis
ap-pen-dix ar-che-typ-al as-sign-a-ble as-sist-ant-ship as-ymp-tot-ic
asyn-chro-nous at-trib-uted at-trib-ut-able bank-rupt bank-rupt-cy
bi-dif-fer-en-tial blue-print busier busiest cat-a-stroph-ic
cat-a-stroph-i-cally con-gress cross-hatched data-base de-fin-i-tive
de-riv-a-tive dis-trib-ute dri-ver dri-vers eco-nom-ics econ-o-mist
elit-ist equi-vari-ant ex-quis-ite ex-tra-or-di-nary flow-chart
for-mi-da-ble forth-right friv-o-lous ge-o-des-ic ge-o-det-ic geo-met-ric
griev-ance griev-ous griev-ous-ly hexa-dec-i-mal ho-lo-no-my ho-mo-thetic
ideals idio-syn-crasy in-fin-ite-ly in-fin-i-tes-i-mal ir-rev-o-ca-ble
key-stroke lam-en-ta-ble light-weight mal-a-prop-ism man-u-script
mar-gin-al meta-bol-ic me-tab-o-lism meta-lan-guage me-trop-o-lis
met-ro-pol-i-tan mi-nut-est mol-e-cule mono-chrome mono-pole mo-nop-oly
mono-spline mo-not-o-nous mul-ti-fac-eted mul-ti-plic-able non-euclid-ean
non-iso-mor-phic non-smooth par-a-digm par-a-bol-ic pa-rab-o-loid
pa-ram-e-trize para-mount pen-ta-gon phe-nom-e-non post-script pre-am-ble
pro-ce-dur-al pro-hib-i-tive pro-hib-i-tive-ly pseu-do-dif-fer-en-tial
pseu-do-fi-nite pseu-do-nym qua-drat-ics quad-ra-ture qua-si-smooth
qua-si-sta-tion-ary qua-si-tri-an-gu-lar quin-tes-sence quin-tes-sen-tial
re-arrange-ment rec-tan-gle ret-ri-bu-tion retro-fit retro-fit-ted
right-eous right-eous-ness ro-bot ro-bot-ics sched-ul-ing se-mes-ter
semi-def-i-nite semi-ho-mo-thet-ic set-up se-vere-ly side-step sov-er-eign
spe-cious spher-oid spher-oid-al star-tling star-tling-ly
sta-tis-tics sto-chas-tic straight-est strange-ness strat-a-gem strong-hold
sum-ma-ble symp-to-matic syn-chro-nous topo-graph-i-cal tra-vers-a-ble
tra-ver-sal tra-ver-sals treach-ery turn-around un-at-tached un-err-ing-ly
white-space wide-spread wing-spread wretch-ed wretch-ed-ly Brown-ian
Eng-lish Euler-ian Feb-ru-ary Gauss-ian Grothen-dieck Hamil-ton-ian
Her-mit-ian Jan-u-ary Japan-ese Kor-te-weg Le-gendre Lip-schitz
Lip-schitz-ian Mar-kov-ian Noe-ther-ian No-vem-ber Rie-mann-ian
Schwarz-schild Sep-tem-ber
form per-iods Uni-ver-si-ty cri-ti-sism for-ma-lism}
\Invalid@\nofrills
\Invalid@\usualspace
\newif\ifnofrills@
\def\nofrills@#1#2{\relaxnext@
  \DN@{\ifx\next\nofrills
    \nofrills@true\let#2\relax\DN@\nofrills{\nextii@}%
  \else
    \nofrills@false\def#2{#1}\let\next@\nextii@\fi
\next@}}
\def\usualspace@#1{\ifnofrills@\def\usualspace{#1}\fi}
\def\addto#1#2{\csname \expandafter\eat@\string#1@\endcsname
  \expandafter{\the\csname \expandafter\eat@\string#1@\endcsname#2}}
\newdimen\bigsize@
\def\big@#1#2{{\hbox{$\left#2\vcenter to#1\bigsize@{}%
  \right.\nulldelimiterspace\z@\m@th$}}}
\def\big{\big@\@ne}
\def\Big{\big@{1.5}}
\def\bigg{\big@\tw@}
\def\Bigg{\big@{2.5}}
\def\raggedcenter@{\leftskip\z@ plus.4\hsize \rightskip\leftskip
 \parfillskip\z@ \parindent\z@ \spaceskip.3333em \xspaceskip.5em
 \pretolerance9999\tolerance9999 \exhyphenpenalty\@M
 \hyphenpenalty\@M \let\\\linebreak}
\def\upperspecialchars{\def\ss{SS}\let\i=I\let\j=J\let\ae\AE\let\oe\OE
  \let\o\O\let\aa\AA\let\l\L}
\def\uppercasetext@#1{%
  {\spaceskip1.2\fontdimen2\the\font plus1.2\fontdimen3\the\font
   \upperspecialchars\uctext@#1$\m@th\aftergroup\eat@$}}
\def\uctext@#1$#2${\endash@#1-\endash@$#2$\uctext@}
\def\endash@#1-#2\endash@{\uppercase{#1}\if\notempty{#2}--\endash@#2\endash@\fi}
\def\runaway@#1{\DN@{#1}\ifx\envir@\next@
  \Err@{You seem to have a missing or misspelled \string\end#1 ...}%
  \let\envir@\empty\fi}
\newif\iftemp@
\def\notempty#1{TT\fi\def\test@{#1}\ifx\test@\empty\temp@false
  \else\temp@true\fi \iftemp@}
\font@\tensmc=cmcsc10
\font@\sevenex=cmex7
\font@\sevenit=cmti7
\font@\eightrm=cmr8 
\font@\sixrm=cmr6 
\font@\eighti=cmmi8     \skewchar\eighti='177 
\font@\sixi=cmmi6       \skewchar\sixi='177   
\font@\eightsy=cmsy8    \skewchar\eightsy='60 
\font@\sixsy=cmsy6      \skewchar\sixsy='60   
\font@\eightex=cmex8
\font@\eightbf=cmbx8 
\font@\sixbf=cmbx6   
\font@\eightit=cmti8 
\font@\eightsl=cmsl8 
\font@\eightsmc=cmcsc8
\font@\eighttt=cmtt8 


\loadmsam
\loadmsbm
\loadeufm
\UseAMSsymbols
\newtoks\tenpoint@
\def\tenpoint{\normalbaselineskip12\p@
 \abovedisplayskip12\p@ plus3\p@ minus9\p@
 \belowdisplayskip\abovedisplayskip
 \abovedisplayshortskip\z@ plus3\p@
 \belowdisplayshortskip7\p@ plus3\p@ minus4\p@
 \textonlyfont@\rm\tenrm \textonlyfont@\it\tenit
 \textonlyfont@\sl\tensl \textonlyfont@\bf\tenbf
 \textonlyfont@\smc\tensmc \textonlyfont@\tt\tentt
 \textonlyfont@\bsmc\tenbsmc
 \ifsyntax@ \def\big##1{{\hbox{$\left##1\right.$}}}%
  \let\Big\big \let\bigg\big \let\Bigg\big
 \else
  \textfont\z@=\tenrm  \scriptfont\z@=\sevenrm  \scriptscriptfont\z@=\fiverm
  \textfont\@ne=\teni  \scriptfont\@ne=\seveni  \scriptscriptfont\@ne=\fivei
  \textfont\tw@=\tensy \scriptfont\tw@=\sevensy \scriptscriptfont\tw@=\fivesy
  \textfont\thr@@=\tenex \scriptfont\thr@@=\sevenex
        \scriptscriptfont\thr@@=\sevenex
  \textfont\itfam=\tenit \scriptfont\itfam=\sevenit
        \scriptscriptfont\itfam=\sevenit
  \textfont\bffam=\tenbf \scriptfont\bffam=\sevenbf
        \scriptscriptfont\bffam=\fivebf
  \setbox\strutbox\hbox{\vrule height8.5\p@ depth3.5\p@ width\z@}%
  \setbox\strutbox@\hbox{\lower.5\normallineskiplimit\vbox{%
        \kern-\normallineskiplimit\copy\strutbox}}%
 \setbox\z@\vbox{\hbox{$($}\kern\z@}\bigsize@=1.2\ht\z@
 \fi
 \normalbaselines\rm\ex@.2326ex\jot3\ex@\the\tenpoint@}
\newtoks\eightpoint@
\def\eightpoint{\normalbaselineskip10\p@
 \abovedisplayskip10\p@ plus2.4\p@ minus7.2\p@
 \belowdisplayskip\abovedisplayskip
 \abovedisplayshortskip\z@ plus2.4\p@
 \belowdisplayshortskip5.6\p@ plus2.4\p@ minus3.2\p@
 \textonlyfont@\rm\eightrm \textonlyfont@\it\eightit
 \textonlyfont@\sl\eightsl \textonlyfont@\bf\eightbf
 \textonlyfont@\smc\eightsmc \textonlyfont@\tt\eighttt
 \textonlyfont@\bsmc\eightbsmc
 \ifsyntax@\def\big##1{{\hbox{$\left##1\right.$}}}%
  \let\Big\big \let\bigg\big \let\Bigg\big
 \else
  \textfont\z@=\eightrm \scriptfont\z@=\sixrm \scriptscriptfont\z@=\fiverm
  \textfont\@ne=\eighti \scriptfont\@ne=\sixi \scriptscriptfont\@ne=\fivei
  \textfont\tw@=\eightsy \scriptfont\tw@=\sixsy \scriptscriptfont\tw@=\fivesy
  \textfont\thr@@=\eightex \scriptfont\thr@@=\sevenex
   \scriptscriptfont\thr@@=\sevenex
  \textfont\itfam=\eightit \scriptfont\itfam=\sevenit
   \scriptscriptfont\itfam=\sevenit
  \textfont\bffam=\eightbf \scriptfont\bffam=\sixbf
   \scriptscriptfont\bffam=\fivebf
 \setbox\strutbox\hbox{\vrule height7\p@ depth3\p@ width\z@}%
 \setbox\strutbox@\hbox{\raise.5\normallineskiplimit\vbox{%
   \kern-\normallineskiplimit\copy\strutbox}}%
 \setbox\z@\vbox{\hbox{$($}\kern\z@}\bigsize@=1.2\ht\z@
 \fi
 \normalbaselines\eightrm\ex@.2326ex\jot3\ex@\the\eightpoint@}

\font@\twelverm=cmr10 scaled\magstep1
\font@\twelveit=cmti10 scaled\magstep1
\font@\twelvesl=cmsl10 scaled\magstep1
\font@\twelvesmc=cmcsc10 scaled\magstep1
\font@\twelvett=cmtt10 scaled\magstep1
\font@\twelvebf=cmbx10 scaled\magstep1
\font@\twelvei=cmmi10 scaled\magstep1
\font@\twelvesy=cmsy10 scaled\magstep1
\font@\twelveex=cmex10 scaled\magstep1
\font@\twelvemsa=msam10 scaled\magstep1
\font@\twelveeufm=eufm10 scaled\magstep1
\font@\twelvemsb=msbm10 scaled\magstep1
\newtoks\twelvepoint@
\def\twelvepoint{\normalbaselineskip15\p@
 \abovedisplayskip15\p@ plus3.6\p@ minus10.8\p@
 \belowdisplayskip\abovedisplayskip
 \abovedisplayshortskip\z@ plus3.6\p@
 \belowdisplayshortskip8.4\p@ plus3.6\p@ minus4.8\p@
 \textonlyfont@\rm\twelverm \textonlyfont@\it\twelveit
 \textonlyfont@\sl\twelvesl \textonlyfont@\bf\twelvebf
 \textonlyfont@\smc\twelvesmc \textonlyfont@\tt\twelvett
 \textonlyfont@\bsmc\twelvebsmc
 \ifsyntax@ \def\big##1{{\hbox{$\left##1\right.$}}}%
  \let\Big\big \let\bigg\big \let\Bigg\big
 \else
  \textfont\z@=\twelverm  \scriptfont\z@=\tenrm  \scriptscriptfont\z@=\sevenrm
  \textfont\@ne=\twelvei  \scriptfont\@ne=\teni  \scriptscriptfont\@ne=\seveni
  \textfont\tw@=\twelvesy \scriptfont\tw@=\tensy \scriptscriptfont\tw@=\sevensy
  \textfont\thr@@=\twelveex \scriptfont\thr@@=\tenex
        \scriptscriptfont\thr@@=\tenex
  \textfont\itfam=\twelveit \scriptfont\itfam=\tenit
        \scriptscriptfont\itfam=\tenit
  \textfont\bffam=\twelvebf \scriptfont\bffam=\tenbf
        \scriptscriptfont\bffam=\sevenbf
  \setbox\strutbox\hbox{\vrule height10.2\p@ depth4.2\p@ width\z@}%
  \setbox\strutbox@\hbox{\lower.6\normallineskiplimit\vbox{%
        \kern-\normallineskiplimit\copy\strutbox}}%
 \setbox\z@\vbox{\hbox{$($}\kern\z@}\bigsize@=1.4\ht\z@
 \fi
 \normalbaselines\rm\ex@.2326ex\jot3.6\ex@\the\twelvepoint@}

\def\headfonts{\twelvepoint\bf}

\font@\fourteenrm=cmr10 scaled\magstep2
\font@\fourteenit=cmti10 scaled\magstep2
\font@\fourteensl=cmsl10 scaled\magstep2
\font@\fourteensmc=cmcsc10 scaled\magstep2
\font@\fourteentt=cmtt10 scaled\magstep2
\font@\fourteenbf=cmbx10 scaled\magstep2
\font@\fourteeni=cmmi10 scaled\magstep2
\font@\fourteensy=cmsy10 scaled\magstep2
\font@\fourteenex=cmex10 scaled\magstep2
\font@\fourteenmsa=msam10 scaled\magstep2
\font@\fourteeneufm=eufm10 scaled\magstep2
\font@\fourteenmsb=msbm10 scaled\magstep2
\newtoks\fourteenpoint@
\def\fourteenpoint{\normalbaselineskip15\p@
 \abovedisplayskip18\p@ plus4.3\p@ minus12.9\p@
 \belowdisplayskip\abovedisplayskip
 \abovedisplayshortskip\z@ plus4.3\p@
 \belowdisplayshortskip10.1\p@ plus4.3\p@ minus5.8\p@
 \textonlyfont@\rm\fourteenrm \textonlyfont@\it\fourteenit
 \textonlyfont@\sl\fourteensl \textonlyfont@\bf\fourteenbf
 \textonlyfont@\smc\fourteensmc \textonlyfont@\tt\fourteentt
 \textonlyfont@\bsmc\fourteenbsmc
 \ifsyntax@ \def\big##1{{\hbox{$\left##1\right.$}}}%
  \let\Big\big \let\bigg\big \let\Bigg\big
 \else
  \textfont\z@=\fourteenrm  \scriptfont\z@=\twelverm  \scriptscriptfont\z@=\tenrm
  \textfont\@ne=\fourteeni  \scriptfont\@ne=\twelvei  \scriptscriptfont\@ne=\teni
  \textfont\tw@=\fourteensy \scriptfont\tw@=\twelvesy \scriptscriptfont\tw@=\tensy
  \textfont\thr@@=\fourteenex \scriptfont\thr@@=\twelveex
        \scriptscriptfont\thr@@=\twelveex
  \textfont\itfam=\fourteenit \scriptfont\itfam=\twelveit
        \scriptscriptfont\itfam=\twelveit
  \textfont\bffam=\fourteenbf \scriptfont\bffam=\twelvebf
        \scriptscriptfont\bffam=\tenbf
  \setbox\strutbox\hbox{\vrule height12.2\p@ depth5\p@ width\z@}%
  \setbox\strutbox@\hbox{\lower.72\normallineskiplimit\vbox{%
        \kern-\normallineskiplimit\copy\strutbox}}%
 \setbox\z@\vbox{\hbox{$($}\kern\z@}\bigsize@=1.7\ht\z@
 \fi
 \normalbaselines\rm\ex@.2326ex\jot4.3\ex@\the\fourteenpoint@}

\def\chapheadfonts{\fourteenpoint\bf}

\font@\seventeenrm=cmr10 scaled\magstep3
\font@\seventeenit=cmti10 scaled\magstep3
\font@\seventeensl=cmsl10 scaled\magstep3
\font@\seventeensmc=cmcsc10 scaled\magstep3
\font@\seventeentt=cmtt10 scaled\magstep3
\font@\seventeenbf=cmbx10 scaled\magstep3
\font@\seventeeni=cmmi10 scaled\magstep3
\font@\seventeensy=cmsy10 scaled\magstep3
\font@\seventeenex=cmex10 scaled\magstep3
\font@\seventeenmsa=msam10 scaled\magstep3
\font@\seventeeneufm=eufm10 scaled\magstep3
\font@\seventeenmsb=msbm10 scaled\magstep3
\newtoks\seventeenpoint@
\def\seventeenpoint{\normalbaselineskip18\p@
 \abovedisplayskip21.6\p@ plus5.2\p@ minus15.4\p@
 \belowdisplayskip\abovedisplayskip
 \abovedisplayshortskip\z@ plus5.2\p@
 \belowdisplayshortskip12.1\p@ plus5.2\p@ minus7\p@
 \textonlyfont@\rm\seventeenrm \textonlyfont@\it\seventeenit
 \textonlyfont@\sl\seventeensl \textonlyfont@\bf\seventeenbf
 \textonlyfont@\smc\seventeensmc \textonlyfont@\tt\seventeentt
 \textonlyfont@\bsmc\seventeenbsmc
 \ifsyntax@ \def\big##1{{\hbox{$\left##1\right.$}}}%
  \let\Big\big \let\bigg\big \let\Bigg\big
 \else
  \textfont\z@=\seventeenrm  \scriptfont\z@=\fourteenrm  \scriptscriptfont\z@=\twelverm
  \textfont\@ne=\seventeeni  \scriptfont\@ne=\fourteeni  \scriptscriptfont\@ne=\twelvei
  \textfont\tw@=\seventeensy \scriptfont\tw@=\fourteensy \scriptscriptfont\tw@=\twelvesy
  \textfont\thr@@=\seventeenex \scriptfont\thr@@=\fourteenex
        \scriptscriptfont\thr@@=\fourteenex
  \textfont\itfam=\seventeenit \scriptfont\itfam=\fourteenit
        \scriptscriptfont\itfam=\fourteenit
  \textfont\bffam=\seventeenbf \scriptfont\bffam=\fourteenbf
        \scriptscriptfont\bffam=\twelvebf
  \setbox\strutbox\hbox{\vrule height14.6\p@ depth6\p@ width\z@}%
  \setbox\strutbox@\hbox{\lower.86\normallineskiplimit\vbox{%
        \kern-\normallineskiplimit\copy\strutbox}}%
 \setbox\z@\vbox{\hbox{$($}\kern\z@}\bigsize@=2\ht\z@
 \fi
 \normalbaselines\rm\ex@.2326ex\jot5.2\ex@\the\seventeenpoint@}

\font@\rrrrrm=cmr10 scaled\magstep4
\font@\bigtitlefont=cmbx10 scaled\magstep4

\parindent1pc
\normallineskiplimit\p@
\newdimen\indenti \indenti=2pc
\def\pageheight#1{\vsize#1}
\def\pagewidth#1{\hsize#1%
   \captionwidth@\hsize \advance\captionwidth@-2\indenti}
\pagewidth{30pc} \pageheight{47pc}
\def\topmatter{%
 \ifx\undefined\msafam
 \else\font@\eightmsa=msam8 \font@\sixmsa=msam6
   \ifsyntax@\else \addto\tenpoint{\textfont\msafam=\tenmsa
              \scriptfont\msafam=\sevenmsa \scriptscriptfont\msafam=\fivemsa}%
     \addto\eightpoint{\textfont\msafam=\eightmsa \scriptfont\msafam=\sixmsa
              \scriptscriptfont\msafam=\fivemsa}%
   \fi
 \fi
 \ifx\undefined\msbfam
 \else\font@\eightmsb=msbm8 \font@\sixmsb=msbm6
   \ifsyntax@\else \addto\tenpoint{\textfont\msbfam=\tenmsb
         \scriptfont\msbfam=\sevenmsb \scriptscriptfont\msbfam=\fivemsb}%
     \addto\eightpoint{\textfont\msbfam=\eightmsb \scriptfont\msbfam=\sixmsb
         \scriptscriptfont\msbfam=\fivemsb}%
   \fi
 \fi
 \ifx\undefined\eufmfam
 \else \font@\eighteufm=eufm8 \font@\sixeufm=eufm6
   \ifsyntax@\else \addto\tenpoint{\textfont\eufmfam=\teneufm
       \scriptfont\eufmfam=\seveneufm \scriptscriptfont\eufmfam=\fiveeufm}%
     \addto\eightpoint{\textfont\eufmfam=\eighteufm
       \scriptfont\eufmfam=\sixeufm \scriptscriptfont\eufmfam=\fiveeufm}%
   \fi
 \fi
 \ifx\undefined\eufbfam
 \else \font@\eighteufb=eufb8 \font@\sixeufb=eufb6
   \ifsyntax@\else \addto\tenpoint{\textfont\eufbfam=\teneufb
      \scriptfont\eufbfam=\seveneufb \scriptscriptfont\eufbfam=\fiveeufb}%
    \addto\eightpoint{\textfont\eufbfam=\eighteufb
      \scriptfont\eufbfam=\sixeufb \scriptscriptfont\eufbfam=\fiveeufb}%
   \fi
 \fi
 \ifx\undefined\eusmfam
 \else \font@\eighteusm=eusm8 \font@\sixeusm=eusm6
   \ifsyntax@\else \addto\tenpoint{\textfont\eusmfam=\teneusm
       \scriptfont\eusmfam=\seveneusm \scriptscriptfont\eusmfam=\fiveeusm}%
     \addto\eightpoint{\textfont\eusmfam=\eighteusm
       \scriptfont\eusmfam=\sixeusm \scriptscriptfont\eusmfam=\fiveeusm}%
   \fi
 \fi
 \ifx\undefined\eusbfam
 \else \font@\eighteusb=eusb8 \font@\sixeusb=eusb6
   \ifsyntax@\else \addto\tenpoint{\textfont\eusbfam=\teneusb
       \scriptfont\eusbfam=\seveneusb \scriptscriptfont\eusbfam=\fiveeusb}%
     \addto\eightpoint{\textfont\eusbfam=\eighteusb
       \scriptfont\eusbfam=\sixeusb \scriptscriptfont\eusbfam=\fiveeusb}%
   \fi
 \fi
 \ifx\undefined\eurmfam
 \else \font@\eighteurm=eurm8 \font@\sixeurm=eurm6
   \ifsyntax@\else \addto\tenpoint{\textfont\eurmfam=\teneurm
       \scriptfont\eurmfam=\seveneurm \scriptscriptfont\eurmfam=\fiveeurm}%
     \addto\eightpoint{\textfont\eurmfam=\eighteurm
       \scriptfont\eurmfam=\sixeurm \scriptscriptfont\eurmfam=\fiveeurm}%
   \fi
 \fi
 \ifx\undefined\eurbfam
 \else \font@\eighteurb=eurb8 \font@\sixeurb=eurb6
   \ifsyntax@\else \addto\tenpoint{\textfont\eurbfam=\teneurb
       \scriptfont\eurbfam=\seveneurb \scriptscriptfont\eurbfam=\fiveeurb}%
    \addto\eightpoint{\textfont\eurbfam=\eighteurb
       \scriptfont\eurbfam=\sixeurb \scriptscriptfont\eurbfam=\fiveeurb}%
   \fi
 \fi
 \ifx\undefined\cmmibfam
 \else \font@\eightcmmib=cmmib8 \font@\sixcmmib=cmmib6
   \ifsyntax@\else \addto\tenpoint{\textfont\cmmibfam=\tencmmib
       \scriptfont\cmmibfam=\sevencmmib \scriptscriptfont\cmmibfam=\fivecmmib}%
    \addto\eightpoint{\textfont\cmmibfam=\eightcmmib
       \scriptfont\cmmibfam=\sixcmmib \scriptscriptfont\cmmibfam=\fivecmmib}%
   \fi
 \fi
 \ifx\undefined\cmbsyfam
 \else \font@\eightcmbsy=cmbsy8 \font@\sixcmbsy=cmbsy6
   \ifsyntax@\else \addto\tenpoint{\textfont\cmbsyfam=\tencmbsy
      \scriptfont\cmbsyfam=\sevencmbsy \scriptscriptfont\cmbsyfam=\fivecmbsy}%
    \addto\eightpoint{\textfont\cmbsyfam=\eightcmbsy
      \scriptfont\cmbsyfam=\sixcmbsy \scriptscriptfont\cmbsyfam=\fivecmbsy}%
   \fi
 \fi
 \let\topmatter\relax}
\def\chapterno@{\uppercase\expandafter{\romannumeral\chaptercount@}}
\newcount\chaptercount@
\def\chapter{\nofrills@{\afterassignment\chapterno@
                        CHAPTER \global\chaptercount@=}\chapter@
 \DNii@##1{\leavevmode\hskip-\leftskip
   \rlap{\vbox to\z@{\vss\centerline{\eightpoint
   \chapter@##1\unskip}\baselineskip2pc\null}}\hskip\leftskip
   \nofrills@false}%
 \FN@\next@}
\newbox\titlebox@

\def\title{\nofrills@{\relax}\title@%
 \DNii@##1\endtitle{\global\setbox\titlebox@\vtop{\tenpoint\bf
 \raggedcenter@\ignorespaces
 \baselineskip1.3\baselineskip\title@{##1}\endgraf}%
 \ifmonograph@ \edef\next{\the\leftheadtoks}\ifx\next\empty
    \leftheadtext{##1}\fi
 \fi
 \edef\next{\the\rightheadtoks}\ifx\next\empty \rightheadtext{##1}\fi
 }\FN@\next@}
\newbox\authorbox@
\def\author#1\endauthor{\global\setbox\authorbox@
 \vbox{\tenpoint\smc\raggedcenter@\ignorespaces
 #1\endgraf}\relaxnext@ \edef\next{\the\leftheadtoks}%
 \ifx\next\empty\leftheadtext{#1}\fi}
\newbox\affilbox@
\def\affil#1\endaffil{\global\setbox\affilbox@
 \vbox{\tenpoint\raggedcenter@\ignorespaces#1\endgraf}}
\newcount\addresscount@
\addresscount@\z@
\def\address#1\endaddress{\global\advance\addresscount@\@ne
  \expandafter\gdef\csname address\number\addresscount@\endcsname
  {\vskip12\p@ minus6\p@\noindent\eightpoint\smc\ignorespaces#1\par}}
\def\email{\nofrills@{\eightpoint{\it E-mail\/}:\enspace}\email@
  \DNii@##1\endemail{%
  \expandafter\gdef\csname email\number\addresscount@\endcsname
  {\def\usualspace{{\it\enspace}}\smallskip\noindent\eightpoint\email@
  \ignorespaces##1\par}}%
 \FN@\next@}
\def\thedate@{}
\def\date#1\enddate{\gdef\thedate@{\tenpoint\ignorespaces#1\unskip}}
\def\thethanks@{}
\def\thanks#1\endthanks{\gdef\thethanks@{\eightpoint\ignorespaces#1.\unskip}}
\def\thekeywords@{}
\def\keywords{\nofrills@{{\it Key words and phrases.\enspace}}\keywords@
 \DNii@##1\endkeywords{\def\thekeywords@{\def\usualspace{{\it\enspace}}%
 \eightpoint\keywords@\ignorespaces##1\unskip.}}%
 \FN@\next@}
\def\thesubjclass@{}
\def\subjclass{\nofrills@{{\rm2000 {\it Mathematics Subject
   Classification\/}.\enspace}}\subjclass@
 \DNii@##1\endsubjclass{\def\thesubjclass@{\def\usualspace
  {{\rm\enspace}}\eightpoint\subjclass@\ignorespaces##1\unskip.}}%
 \FN@\next@}
\newbox\abstractbox@
\def\abstract{\nofrills@{{\smc Abstract.\enspace}}\abstract@
 \DNii@{\setbox\abstractbox@\vbox\bgroup\noindent$$\vbox\bgroup
  \def\envir@{abstract}\advance\hsize-2\indenti
  \usualspace@{{\enspace}}\eightpoint \noindent\abstract@\ignorespaces}%
 \FN@\next@}
\def\endabstract{\par\unskip\egroup$$\egroup}
\def\widestnumber#1#2{\begingroup\let\head\null\let\subhead\empty
   \let\subsubhead\subhead
   \ifx#1\head\global\setbox\tocheadbox@\hbox{#2.\enspace}%
   \else\ifx#1\subhead\global\setbox\tocsubheadbox@\hbox{#2.\enspace}%
   \else\ifx#1\key\bgroup\let\endrefitem@\egroup
        \key#2\endrefitem@\global\refindentwd\wd\keybox@
   \else\ifx#1\no\bgroup\let\endrefitem@\egroup
        \no#2\endrefitem@\global\refindentwd\wd\nobox@
   \else\ifx#1\page\global\setbox\pagesbox@\hbox{\quad\bf#2}%
   \else\ifx#1\item\setboxz@h{#2}\global\rosteritemwd\wdz@
        \global\advance\rosteritemwd by.5\parindent
   \else\message{\string\widestnumber is not defined for this option
   (\string#1)}%
\fi\fi\fi\fi\fi\fi\endgroup}
\newif\ifmonograph@
\def\Monograph{\monograph@true \let\headmark\rightheadtext
  \let\varindent@\indent \def\headfont@{\bf}\def\proclaimheadfont@{\smc}%
  \def\demofont@{\smc}}
\let\varindent@\indent

\newbox\tocheadbox@    \newbox\tocsubheadbox@
\newbox\tocbox@
\def\toc{\toc@{Contents}}
\def\newtocdefs{%
   \def \title##1\endtitle
       {\penaltyandskip@\z@\smallskipamount
        \hangindent\wd\tocheadbox@\noindent{\bf##1}}%
   \def \chapter##1{%
        Chapter \uppercase\expandafter{\romannumeral##1.\unskip}\enspace}%
   \def \specialhead##1\endspecialhead
       {\par\hangindent\wd\tocheadbox@ \noindent##1\par}%
   \def \head##1 ##2\endhead
       {\par\hangindent\wd\tocheadbox@ \noindent
        \if\notempty{##1}\hbox to\wd\tocheadbox@{\hfil##1\enspace}\fi
        ##2\par}%
   \def \subhead##1 ##2\endsubhead
       {\par\vskip-\parskip {\normalbaselines
        \advance\leftskip\wd\tocheadbox@
        \hangindent\wd\tocsubheadbox@ \noindent
        \if\notempty{##1}\hbox to\wd\tocsubheadbox@{##1\unskip\hfil}\fi
         ##2\par}}%
   \def \subsubhead##1 ##2\endsubsubhead
       {\par\vskip-\parskip {\normalbaselines
        \advance\leftskip\wd\tocheadbox@
        \hangindent\wd\tocsubheadbox@ \noindent
        \if\notempty{##1}\hbox to\wd\tocsubheadbox@{##1\unskip\hfil}\fi
        ##2\par}}}
\def\toc@#1{\relaxnext@
   \def\page##1%
       {\unskip\penalty0\null\hfil
        \rlap{\hbox to\wd\pagesbox@{\quad\hfil##1}}\hfilneg\penalty\@M}%
 \DN@{\ifx\next\nofrills\DN@\nofrills{\nextii@}%
      \else\DN@{\nextii@{{#1}}}\fi
      \next@}%
 \DNii@##1{%
\ifmonograph@\bgroup\else\setbox\tocbox@\vbox\bgroup
   \centerline{\headfont@\ignorespaces##1\unskip}\nobreak
   \vskip\belowheadskip \fi
   \setbox\tocheadbox@\hbox{0.\enspace}%
   \setbox\tocsubheadbox@\hbox{0.0.\enspace}%
   \leftskip\indenti \rightskip\leftskip
   \setbox\pagesbox@\hbox{\bf\quad000}\advance\rightskip\wd\pagesbox@
   \newtocdefs
 }%
 \FN@\next@}
\def\endtoc{\par\egroup}
\let\pretitle\relax
\let\preauthor\relax
\let\preaffil\relax
\let\predate\relax
\let\preabstract\relax
\let\prepaper\relax
\def\dedicatory #1\enddedicatory{\def\preabstract{{\medskip
  \eightpoint\it \raggedcenter@#1\endgraf}}}
\def\thetranslator@{}
\def\translator#1\endtranslator{\def\thetranslator@{\nobreak\medskip
 \line{\eightpoint\hfil Translated by \uppercase{#1}\qquad\qquad}\nobreak}}
\outer\def\endtopmatter{\runaway@{abstract}%
 \edef\next{\the\leftheadtoks}\ifx\next\empty
  \expandafter\leftheadtext\expandafter{\the\rightheadtoks}\fi
 \ifmonograph@\else
   \ifx\thesubjclass@\empty\else \makefootnote@{}{\thesubjclass@}\fi
   \ifx\thekeywords@\empty\else \makefootnote@{}{\thekeywords@}\fi
   \ifx\thethanks@\empty\else \makefootnote@{}{\thethanks@}\fi
 \fi
  \pretitle
  \ifmonograph@ \topskip7pc \else \topskip4pc \fi
  \box\titlebox@
  \topskip10pt
  \preauthor
  \ifvoid\authorbox@\else \vskip2.5pc plus1pc \unvbox\authorbox@\fi
  \preaffil
  \ifvoid\affilbox@\else \vskip1pc plus.5pc \unvbox\affilbox@\fi
  \predate
  \ifx\thedate@\empty\else \vskip1pc plus.5pc \line{\hfil\thedate@\hfil}\fi
  \preabstract
  \ifvoid\abstractbox@\else \vskip1.5pc plus.5pc \unvbox\abstractbox@ \fi
  \ifvoid\tocbox@\else\vskip1.5pc plus.5pc \unvbox\tocbox@\fi
  \prepaper
  \vskip2pc plus1pc
}
\def\document{\let\fontlist@\relax\let\alloclist@\relax
  \tenpoint}

\newskip\aboveheadskip       \aboveheadskip1.8\bigskipamount
\newdimen\belowheadskip      \belowheadskip1.8\medskipamount

\def\headfont@{\smc}
\def\penaltyandskip@#1#2{\relax\ifdim\lastskip<#2\relax\removelastskip
      \ifnum#1=\z@\else\penalty@#1\relax\fi\vskip#2%
  \else\ifnum#1=\z@\else\penalty@#1\relax\fi\fi}
\def\nobreak{\penalty\@M
  \ifvmode\def\penalty@{\let\penalty@\penalty\count@@@}%
  \everypar{\let\penalty@\penalty\everypar{}}\fi}
\let\penalty@\penalty
\def\heading#1\endheading{\head#1\endhead}

\def\specialheadfont@{\bf}
\outer\def\specialhead{\par\penaltyandskip@{-200}\aboveheadskip
  \begingroup\interlinepenalty\@M\rightskip\z@ plus\hsize \let\\\linebreak
  \specialheadfont@\noindent\ignorespaces}
\def\endspecialhead{\par\endgroup\nobreak\vskip\belowheadskip}
\let\headmark\eat@
\newskip\subheadskip       \subheadskip\medskipamount
\def\subheadfont@{\bf}
\outer\def\subhead{\nofrills@{.\enspace}\subhead@
 \DNii@##1\endsubhead{\par\penaltyandskip@{-100}\subheadskip
  \varindent@{\usualspace@{{\subheadfont@\enspace}}%
 \subheadfont@\ignorespaces##1\unskip\subhead@}\ignorespaces}%
 \FN@\next@}
\outer\def\subsubhead{\nofrills@{.\enspace}\subsubhead@
 \DNii@##1\endsubsubhead{\par\penaltyandskip@{-50}\medskipamount
      {\usualspace@{{\it\enspace}}%
  \it\ignorespaces##1\unskip\subsubhead@}\ignorespaces}%
 \FN@\next@}
\def\proclaimheadfont@{\bf}
\outer\def\proclaim{\runaway@{proclaim}\def\envir@{proclaim}%
  \nofrills@{.\enspace}\proclaim@
 \DNii@##1{\penaltyandskip@{-100}\medskipamount\varindent@
   \usualspace@{{\proclaimheadfont@\enspace}}\proclaimheadfont@
   \ignorespaces##1\unskip\proclaim@
  \sl\ignorespaces}%
 \FN@\next@}
\outer\def\endproclaim{\let\envir@\relax\par\rm
  \penaltyandskip@{55}\medskipamount}
\def\demoheadfont@{\it}
\def\demo{\runaway@{proclaim}\nofrills@{.\enspace}\demo@
     \DNii@##1{\par\penaltyandskip@\z@\medskipamount
  {\usualspace@{{\demoheadfont@\enspace}}%
  \varindent@\demoheadfont@\ignorespaces##1\unskip\demo@}\rm
  \ignorespaces}\FN@\next@}
\def\enddemo{\par\medskip}
\def\qed{\ifhmode\unskip\nobreak\fi\quad\ifmmode\square\else$\m@th\square$\fi}
\let\remark\demo
\let\endremark\enddemo
\def\definition{\runaway@{proclaim}%
  \nofrills@{.\demoheadfont@\enspace}\definition@
        \DNii@##1{\penaltyandskip@{-100}\medskipamount
        {\usualspace@{{\demoheadfont@\enspace}}%
        \varindent@\demoheadfont@\ignorespaces##1\unskip\definition@}%
        \rm \ignorespaces}\FN@\next@}


\newdimen\rosteritemwd
\newcount\rostercount@
\newif\iffirstitem@
\let\plainitem@\item
\newtoks\everypartoks@
\def\par@{\everypartoks@\expandafter{\the\everypar}\everypar{}}
\def\roster{\edef\leftskip@{\leftskip\the\leftskip}%
 \relaxnext@
 \rostercount@\z@  
 \def\item{\FN@\rosteritem@}%
 \DN@{\ifx\next\runinitem\let\next@\nextii@\else
  \let\next@\nextiii@\fi\next@}%
 \DNii@\runinitem  
  {\unskip  
   \DN@{\ifx\next[\let\next@\nextii@\else
    \ifx\next"\let\next@\nextiii@\else\let\next@\nextiv@\fi\fi\next@}%
   \DNii@[####1]{\rostercount@####1\relax
    \enspace{\rm(\number\rostercount@)}~\ignorespaces}%
   \def\nextiii@"####1"{\enspace{\rm####1}~\ignorespaces}%
   \def\nextiv@{\enspace{\rm(1)}\rostercount@\@ne~}%
   \par@\firstitem@false  
   \FN@\next@}%
 \def\nextiii@{\par\par@  
  \penalty\@m\smallskip\vskip-\parskip
  \firstitem@true}%
 \FN@\next@}
\def\rosteritem@{\iffirstitem@\firstitem@false\else\par\vskip-\parskip\fi
 \leftskip3\parindent\noindent  
 \DNii@[##1]{\rostercount@##1\relax
  \llap{\hbox to2.5\parindent{\hss\rm(\number\rostercount@)}%
   \hskip.5\parindent}\ignorespaces}%
 \def\nextiii@"##1"{%
  \llap{\hbox to2.5\parindent{\hss\rm##1}\hskip.5\parindent}\ignorespaces}%
 \def\nextiv@{\advance\rostercount@\@ne
  \llap{\hbox to2.5\parindent{\hss\rm(\number\rostercount@)}%
   \hskip.5\parindent}}%
 \ifx\next[\let\next@\nextii@\else\ifx\next"\let\next@\nextiii@\else
  \let\next@\nextiv@\fi\fi\next@}

\newif\ifnextRunin@
\def\endroster{\relaxnext@
 \par\leftskip@  
 \penalty-50 \vskip-\parskip\smallskip  
 \DN@{\ifx\next\Runinitem\let\next@\relax
  \else\nextRunin@false\let\item\plainitem@  
   \ifx\next\par 
    \DN@\par{\everypar\expandafter{\the\everypartoks@}}%
   \else  
    \DN@{\noindent\everypar\expandafter{\the\everypartoks@}}%
  \fi\fi\next@}%
 \FN@\next@}
\newcount\rosterhangafter@
\def\Runinitem#1\roster\runinitem{\relaxnext@
 \rostercount@\z@ 
 \def\item{\FN@\rosteritem@}%
 \def\runinitem@{#1}%
 \DN@{\ifx\next[\let\next\nextii@\else\ifx\next"\let\next\nextiii@
  \else\let\next\nextiv@\fi\fi\next}%
 \DNii@[##1]{\rostercount@##1\relax
  \def\item@{{\rm(\number\rostercount@)}}\nextv@}%
 \def\nextiii@"##1"{\def\item@{{\rm##1}}\nextv@}%
 \def\nextiv@{\advance\rostercount@\@ne
  \def\item@{{\rm(\number\rostercount@)}}\nextv@}%
 \def\nextv@{\setbox\z@\vbox  
  {\ifnextRunin@\noindent\fi  
  \runinitem@\unskip\enspace\item@~\par  
  \global\rosterhangafter@\prevgraf}%
  \firstitem@false  
  \ifnextRunin@\else\par\fi  
  \hangafter\rosterhangafter@\hangindent3\parindent
  \ifnextRunin@\noindent\fi  
  \runinitem@\unskip\enspace 
  \item@~\ifnextRunin@\else\par@\fi  
  \nextRunin@true\ignorespaces}%
 \FN@\next@}
\def\footmarkform@#1{$\m@th^{#1}$}
\let\thefootnotemark\footmarkform@
\def\makefootnote@#1#2{\insert\footins
 {\interlinepenalty\interfootnotelinepenalty
 \eightpoint\splittopskip\ht\strutbox\splitmaxdepth\dp\strutbox
 \floatingpenalty\@MM\leftskip\z@\rightskip\z@\spaceskip\z@\xspaceskip\z@
 \leavevmode{#1}\footstrut\ignorespaces#2\unskip\lower\dp\strutbox
 \vbox to\dp\strutbox{}}}
\newcount\footmarkcount@
\footmarkcount@\z@
\def\footnotemark{\let\@sf\empty\relaxnext@
 \ifhmode\edef\@sf{\spacefactor\the\spacefactor}\/\fi
 \DN@{\ifx[\next\let\next@\nextii@\else
  \ifx"\next\let\next@\nextiii@\else
  \let\next@\nextiv@\fi\fi\next@}%
 \DNii@[##1]{\footmarkform@{##1}\@sf}%
 \def\nextiii@"##1"{{##1}\@sf}%
 \def\nextiv@{\iffirstchoice@\global\advance\footmarkcount@\@ne\fi
  \footmarkform@{\number\footmarkcount@}\@sf}%
 \FN@\next@}
\def\footnotetext{\relaxnext@
 \DN@{\ifx[\next\let\next@\nextii@\else
  \ifx"\next\let\next@\nextiii@\else
  \let\next@\nextiv@\fi\fi\next@}%
 \DNii@[##1]##2{\makefootnote@{\footmarkform@{##1}}{##2}}%
 \def\nextiii@"##1"##2{\makefootnote@{##1}{##2}}%
 \def\nextiv@##1{\makefootnote@{\footmarkform@{\number\footmarkcount@}}{##1}}%
 \FN@\next@}
\def\footnote{\let\@sf\empty\relaxnext@
 \ifhmode\edef\@sf{\spacefactor\the\spacefactor}\/\fi
 \DN@{\ifx[\next\let\next@\nextii@\else
  \ifx"\next\let\next@\nextiii@\else
  \let\next@\nextiv@\fi\fi\next@}%
 \DNii@[##1]##2{\footnotemark[##1]\footnotetext[##1]{##2}}%
 \def\nextiii@"##1"##2{\footnotemark"##1"\footnotetext"##1"{##2}}%
 \def\nextiv@##1{\footnotemark\footnotetext{##1}}%
 \FN@\next@}
\def\adjustfootnotemark#1{\advance\footmarkcount@#1\relax}
\def\footnoterule{\kern-3\p@
  \hrule width 5pc\kern 2.6\p@} 
\def\captionfont@{\smc}
\def\topcaption#1#2\endcaption{%
  {\dimen@\hsize \advance\dimen@-\captionwidth@
   \rm\raggedcenter@ \advance\leftskip.5\dimen@ \rightskip\leftskip
  {\captionfont@#1}%
  \if\notempty{#2}.\enspace\ignorespaces#2\fi
  \endgraf}\nobreak\bigskip}
\def\botcaption#1#2\endcaption{%
  \nobreak\bigskip
  \setboxz@h{\captionfont@#1\if\notempty{#2}.\enspace\rm#2\fi}%
  {\dimen@\hsize \advance\dimen@-\captionwidth@
   \leftskip.5\dimen@ \rightskip\leftskip
   \noindent \ifdim\wdz@>\captionwidth@ 
   \else\hfil\fi 
  {\captionfont@#1}\if\notempty{#2}.\enspace\rm#2\fi\endgraf}}
\def\@ins{\par\begingroup\def\vspace##1{\vskip##1\relax}%
  \def\captionwidth##1{\captionwidth@##1\relax}%
  \setbox\z@\vbox\bgroup} 
\def\block{\RIfMIfI@\nondmatherr@\block\fi
       \else\ifvmode\vskip\abovedisplayskip\noindent\fi
        $$\def\endblock{\par\egroup$$}\fi
  \vbox\bgroup\advance\hsize-2\indenti\noindent}
\def\endblock{\par\egroup}
\def\cite#1{{\rm[{\citefont@\m@th#1}]}}
\def\citefont@{\rm}
\def\refsfont@{\eightpoint}
\outer\def\Refs{\runaway@{proclaim}%
 \relaxnext@ \DN@{\ifx\next\nofrills\DN@\nofrills{\nextii@}\else
  \DN@{\nextii@{References}}\fi\next@}%
 \DNii@##1{\penaltyandskip@{-200}\aboveheadskip
  \line{\hfil\headfont@\ignorespaces##1\unskip\hfil}\nobreak
  \vskip\belowheadskip
  \begingroup\refsfont@\sfcode`.=\@m}%
 \FN@\next@}
\def\endRefs{\par\endgroup}
\newbox\nobox@            \newbox\keybox@           \newbox\bybox@
\newbox\paperbox@         \newbox\paperinfobox@     \newbox\jourbox@
\newbox\volbox@           \newbox\issuebox@         \newbox\yrbox@
\newbox\pagesbox@         \newbox\bookbox@          \newbox\bookinfobox@
\newbox\publbox@          \newbox\publaddrbox@      \newbox\finalinfobox@
\newbox\edsbox@           \newbox\langbox@
\newif\iffirstref@        \newif\iflastref@
\newif\ifprevjour@        \newif\ifbook@            \newif\ifprevinbook@
\newif\ifquotes@          \newif\ifbookquotes@      \newif\ifpaperquotes@
\newdimen\bysamerulewd@
\setboxz@h{\refsfont@\kern3em}
\bysamerulewd@\wdz@
\newdimen\refindentwd
\setboxz@h{\refsfont@ 00. }
\refindentwd\wdz@
\outer\def\ref{\begingroup \noindent\hangindent\refindentwd
 \firstref@true \def\nofrills{\def\refkern@{\kern3sp}}%
 \ref@}
\def\ref@{\book@false \bgroup\let\endrefitem@\egroup \ignorespaces}
\def\moreref{\endrefitem@\endref@\firstref@false\ref@}%
\def\transl{\endrefitem@\endref@\firstref@false
  \book@false
  \prepunct@
  \setboxz@h\bgroup \aftergroup\unhbox\aftergroup\z@
    \def\endrefitem@{\unskip\refkern@\egroup}\ignorespaces}%
\def\emptyifempty@{\dimen@\wd\currbox@
  \advance\dimen@-\wd\z@ \advance\dimen@-.1\p@
  \ifdim\dimen@<\z@ \setbox\currbox@\copy\voidb@x \fi}
\let\refkern@\relax
\def\endrefitem@{\unskip\refkern@\egroup
  \setboxz@h{\refkern@}\emptyifempty@}\ignorespaces
\def\refdef@#1#2#3{\edef\next@{\noexpand\endrefitem@
  \let\noexpand\currbox@\csname\expandafter\eat@\string#1box@\endcsname
    \noexpand\setbox\noexpand\currbox@\hbox\bgroup}%
  \toks@\expandafter{\next@}%
  \if\notempty{#2#3}\toks@\expandafter{\the\toks@
  \def\endrefitem@{\unskip#3\refkern@\egroup
  \setboxz@h{#2#3\refkern@}\emptyifempty@}#2}\fi
  \toks@\expandafter{\the\toks@\ignorespaces}%
  \edef#1{\the\toks@}}
\refdef@\no{}{. }
\refdef@\key{[\m@th}{] }
\refdef@\by{}{}
\def\bysame{\by\hbox to\bysamerulewd@{\hrulefill}\thinspace
   \kern0sp}
\def\manyby{\message{\string\manyby is no longer necessary; \string\by
  can be used instead, starting with version 2.0 of \styname.STY}\by}
\refdef@\paper{\ifpaperquotes@``\fi\it}{}
\refdef@\paperinfo{}{}
\def\jour{\endrefitem@\let\currbox@\jourbox@
  \setbox\currbox@\hbox\bgroup
  \def\endrefitem@{\unskip\refkern@\egroup
    \setboxz@h{\refkern@}\emptyifempty@
    \ifvoid\jourbox@\else\prevjour@true\fi}%
\ignorespaces}
\refdef@\vol{\ifbook@\else\bf\fi}{}
\refdef@\issue{no. }{}
\refdef@\yr{}{}
\refdef@\pages{}{}
\def\page{\endrefitem@\def\pp@{\def\pp@{pp.~}p.~}\let\currbox@\pagesbox@
  \setbox\currbox@\hbox\bgroup\ignorespaces}
\def\pp@{pp.~}
\def\book{\endrefitem@ \let\currbox@\bookbox@
 \setbox\currbox@\hbox\bgroup\def\endrefitem@{\unskip\refkern@\egroup
  \setboxz@h{\ifbookquotes@``\fi}\emptyifempty@
  \ifvoid\bookbox@\else\book@true\fi}%
  \ifbookquotes@``\fi\it\ignorespaces}
\def\inbook{\endrefitem@
  \let\currbox@\bookbox@\setbox\currbox@\hbox\bgroup
  \def\endrefitem@{\unskip\refkern@\egroup
  \setboxz@h{\ifbookquotes@``\fi}\emptyifempty@
  \ifvoid\bookbox@\else\book@true\previnbook@true\fi}%
  \ifbookquotes@``\fi\ignorespaces}
\refdef@\eds{(}{, eds.)}
\def\ed{\endrefitem@\let\currbox@\edsbox@
 \setbox\currbox@\hbox\bgroup
 \def\endrefitem@{\unskip, ed.)\refkern@\egroup
  \setboxz@h{(, ed.)}\emptyifempty@}(\ignorespaces}
\refdef@\bookinfo{}{}
\refdef@\publ{}{}
\refdef@\publaddr{}{}
\refdef@\finalinfo{}{}
\refdef@\lang{(}{)}

\let\refdef@\relax 
\def\ppunbox@#1{\ifvoid#1\else\prepunct@\unhbox#1\fi}
\def\nocomma@#1{\ifvoid#1\else\changepunct@3\prepunct@\unhbox#1\fi}
\def\changepunct@#1{\ifnum\lastkern<3 \unkern\kern#1sp\fi}
\def\prepunct@{\count@\lastkern\unkern
  \ifnum\lastpenalty=0
    \let\penalty@\relax
  \else
    \edef\penalty@{\penalty\the\lastpenalty\relax}%
  \fi
  \unpenalty
  \let\refspace@\ \ifcase\count@,
\or;\or.\or 
  \or\let\refspace@\relax
  \else,\fi
  \ifquotes@''\quotes@false\fi \penalty@ \refspace@
}
\def\transferpenalty@#1{\dimen@\lastkern\unkern
  \ifnum\lastpenalty=0\unpenalty\let\penalty@\relax
  \else\edef\penalty@{\penalty\the\lastpenalty\relax}\unpenalty\fi
  #1\penalty@\kern\dimen@}
\def\endref{\endrefitem@\lastref@true\endref@
  \par\endgroup \prevjour@false \previnbook@false }
\def\endref@{%
\iffirstref@
  \ifvoid\nobox@\ifvoid\keybox@\indent\fi
  \else\hbox to\refindentwd{\hss\unhbox\nobox@}\fi
  \ifvoid\keybox@
  \else\ifdim\wd\keybox@>\refindentwd
         \box\keybox@
       \else\hbox to\refindentwd{\unhbox\keybox@\hfil}\fi\fi
  \kern4sp\ppunbox@\bybox@
\fi 
  \ifvoid\paperbox@
  \else\prepunct@\unhbox\paperbox@
    \ifpaperquotes@\quotes@true\fi\fi
  \ppunbox@\paperinfobox@
  \ifvoid\jourbox@
    \ifprevjour@ \nocomma@\volbox@
      \nocomma@\issuebox@
      \ifvoid\yrbox@\else\changepunct@3\prepunct@(\unhbox\yrbox@
        \transferpenalty@)\fi
      \ppunbox@\pagesbox@
    \fi 
  \else \prepunct@\unhbox\jourbox@
    \nocomma@\volbox@
    \nocomma@\issuebox@
    \ifvoid\yrbox@\else\changepunct@3\prepunct@(\unhbox\yrbox@
      \transferpenalty@)\fi
    \ppunbox@\pagesbox@
  \fi 
  \ifbook@\prepunct@\unhbox\bookbox@ \ifbookquotes@\quotes@true\fi \fi
  \nocomma@\edsbox@
  \ppunbox@\bookinfobox@
  \ifbook@\ifvoid\volbox@\else\prepunct@ vol.~\unhbox\volbox@
  \fi\fi
  \ppunbox@\publbox@ \ppunbox@\publaddrbox@
  \ifbook@ \ppunbox@\yrbox@
    \ifvoid\pagesbox@
    \else\prepunct@\pp@\unhbox\pagesbox@\fi
  \else
    \ifprevinbook@ \ppunbox@\yrbox@
      \ifvoid\pagesbox@\else\prepunct@\pp@\unhbox\pagesbox@\fi
    \fi \fi
  \ppunbox@\finalinfobox@
  \iflastref@
    \ifvoid\langbox@.\ifquotes@''\fi
    \else\changepunct@2\prepunct@\unhbox\langbox@\fi
  \else
    \ifvoid\langbox@\changepunct@1%
    \else\changepunct@3\prepunct@\unhbox\langbox@
      \changepunct@1\fi
  \fi
}
\outer\def\enddocument{%
 \runaway@{proclaim}%
\ifmonograph@ 
\else
 \nobreak
 \thetranslator@
 \count@\z@ \loop\ifnum\count@<\addresscount@\advance\count@\@ne
 \csname address\number\count@\endcsname
 \csname email\number\count@\endcsname
 \repeat
\fi
 \vfill\supereject\end}

\def\headfont@{\headfonts}
\def\proclaimheadfont@{\bf}
\def\specialheadfont@{\bf}
\def\subheadfont@{\bf}
\def\demoheadfont@{\smc}

\newif\ifThisToToc \ThisToTocfalse
\newif\iftocloaded \tocloadedfalse

\def\C@L{\noexpand\Cal}\def\B@B{\noexpand\Bbb}\def\fR@K{\noexpand\frak}
\def\S@{\noexpand\S}\def\P@P{\noexpand\"}
\def\xpar{\\}

\def\writetoc#1{\iftocloaded\ifThisToToc\begingroup\def\totoc{}
  \def\Cal{\noexpand\C@L}\def\Bbb{\noexpand\B@B}
  \def\frak{\noexpand\fR@K}\def\goth{\frak}\def\S{\noexpand\S@}
  \def\"{\noexpand\P@P}
  \def\xpar{\par\penalty100000 }\def\idx##1{##1}\def\\{\xpar}
  \edef\next@{\write\toc{\noindent#1\leaderfill\noexpand\folio\par}}%
  \next@\endgroup\global\ThisToTocfalse\fi\fi}
\def\leaderfill{\leaders\hbox to 1em{\hss.\hss}\hfill}

\newif\ifindexloaded \indexloadedfalse
\def\idx#1{\ifindexloaded\begingroup\def\ign{}\def\it{}\def\/{}%
 \def\smc{}\def\bf{}\def\tt{}%
 \def\Cal{\noexpand\C@L}\def\Bbb{\noexpand\B@B}%
 \def\frak{\noexpand\fR@K}\def\goth{\frak}\def\S{\noexpand\S@}%
  \def\"{\noexpand\P@P}%
 {\edef\next@{\write\index{#1, \noexpand\folio}}\next@}%
 \endgroup\fi{#1}}
\def\ign#1{}

\def\totoc{\global\ThisToToctrue}

\outer\def\head#1\endhead{\par\penaltyandskip@{-200}\aboveheadskip
 {\headfont@\raggedcenter@\interlinepenalty\@M
 \ignorespaces#1\endgraf}\nobreak
 \vskip\belowheadskip
 \headmark{#1}\writetoc{#1}}

\outer\def\chaphead#1\endchaphead{\par\penaltyandskip@{-200}\aboveheadskip
 {\chapheadfonts\raggedcenter@\interlinepenalty\@M
 \ignorespaces#1\endgraf}\nobreak
 \vskip3\belowheadskip
 \headmark{#1}\writetoc{#1}}

\def\folio{{\foliofont@\ifnum\pageno<\z@ \romannumeral-\pageno
 \else\number\pageno \fi}}
\newtoks\leftheadtoks
\newtoks\rightheadtoks

\def\leftheadtext{\nofrills@{\relax}\lht@
  \DNii@##1{\leftheadtoks\expandafter{\lht@{##1}}%
    \mark{\the\leftheadtoks\noexpand\else\the\rightheadtoks}
    \ifsyntax@\setboxz@h{\def\\{\unskip\space\ignorespaces}%
        \headlinefont@##1}\fi}%
  \FN@\next@}
\def\rightheadtext{\nofrills@{\relax}\rht@
  \DNii@##1{\rightheadtoks\expandafter{\rht@{##1}}%
    \mark{\the\leftheadtoks\noexpand\else\the\rightheadtoks}%
    \ifsyntax@\setboxz@h{\def\\{\unskip\space\ignorespaces}%
        \headlinefont@##1}\fi}%
  \FN@\next@}
\def\NoRunningHeads{\global\runheads@false\global\let\headmark\eat@}

\newif\iffirstpage@     \firstpage@true
\newif\ifrunheads@      \runheads@true

\newdimen\fullhsize \fullhsize=\hsize
\newdimen\fullvsize \fullvsize=\vsize
\def\fullline{\hbox to\fullhsize}

\def\pagenumbers{\gdef\folio{\folio@}}

\let\norunningheads\NoRunningHeads
\def\userunningheads{\global\runheads@true}
\norunningheads

\headline={\def\chapter#1{\chapterno@. }%
  \def\\{\unskip\space\ignorespaces}\ifrunheads@\headlinefont@
    \ifodd\pageno\rightheadline \else\leftheadline\fi
   \else\hfil\fi\ifNoRunHeadline\global\NoRunHeadlinefalse\fi}
\let\folio@\folio
\def\foliofont@{\foliofont}
\def\foliofont{\eightrm}
\def\headlinefont@{\headlinefont}
\def\headlinefont{\eightpoint\smc}
\def\leftheadline{\rlap{\folio}\hfill
   \ifNoRunHeadline\else\iftrue\topmark\fi\fi \hfill}
\def\rightheadline{\hfill\ifNoRunHeadline
   \else \expandafter\fi
  \hfill \llap{\folio}}
\footline={{\eightpoint\bottremark}%
   \ifrunheads@\else\hfil{\let\foliofont\tenrm\folio}\fi\hfil}
\def\bottremark{}
 
\newif\ifNoRunHeadline      
\def\norunninghead{\global\NoRunHeadlinetrue}
\norunninghead

\output={\output@}
%
\newif\ifoffset\offsetfalse
\output={\output@}
\def\output@{%
 \ifoffset 
  \ifodd\count0\advance\hoffset by0.5truecm
   \else\advance\hoffset by-0.5truecm\fi\fi
 \shipout\vbox{%
  \makeheadline \pagebody \makefootline }%
 \advancepageno \ifnum\outputpenalty>-\@MM\else\dosupereject\fi}

\def\indexoutput#1{%
  \ifoffset 
   \ifodd\count0\advance\hoffset by0.5truecm
    \else\advance\hoffset by-0.5truecm\fi\fi
  \shipout\vbox{\makeheadline
  \vbox to\fullvsize{\boxmaxdepth\maxdepth%
  \ifvoid\topins\else\unvbox\topins\fi%
  #1 %
  \ifvoid\footins\else 
    \vskip\skip\footins
    \footnoterule
    \unvbox\footins\fi
  \ifr@ggedbottom \kern-\dimen@ \vfil \fi}%
  \baselineskip2pc
  \makefootline}%
 \global\advance\pageno\@ne
 \ifnum\outputpenalty>-\@MM\else\dosupereject\fi}
 
 \newbox\partialpage \newdimen\halfsize \halfsize=0.5\fullhsize
 \advance\halfsize by-0.5em

 \def\begindoublecolumns{\output={\indexoutput{\unvbox255}}%
   \begingroup \def\line{\fullline}
   \output={\global\setbox\partialpage=\vbox{\unvbox255\bigskip}}\eject
   \output={\doublecolumnout}\hsize=\halfsize \vsize=2\fullvsize}
 \def\enddoublecolumns{\output={\balancecolumns}\eject
  \endgroup \pagegoal=\fullvsize%
  \output={\output@}}
\def\doublecolumnout{\splittopskip=\topskip \splitmaxdepth=\maxdepth
  \dimen@=\fullvsize \advance\dimen@ by-\ht\partialpage
  \setbox0=\vsplit255 to \dimen@ \setbox2=\vsplit255 to \dimen@
  \indexoutput{\pagesofar} \unvbox255 \penalty\outputpenalty}
\def\pagesofar{\unvbox\partialpage
  \wd0=\hsize \wd2=\hsize \hbox to\fullhsize{\box0\hfil\box2}}
\def\balancecolumns{\setbox0=\vbox{\unvbox255} \dimen@=\ht0
  \advance\dimen@ by\topskip \advance\dimen@ by-\baselineskip
  \divide\dimen@ by2 \splittopskip=\topskip
  {\vbadness=10000 \loop \global\setbox3=\copy0
    \global\setbox1=\vsplit3 to\dimen@
    \ifdim\ht3>\dimen@ \global\advance\dimen@ by1pt \repeat}
  \setbox0=\vbox to\dimen@{\unvbox1} \setbox2=\vbox to\dimen@{\unvbox3}
  \pagesofar}

\tenpoint
\catcode`\@=\active

\def\smallheadings{\let\chapheadfonts\tenpoint\let\headfonts\tenpoint}

\tenpoint
\catcode`\@=\active

\def\LL{\leavevmode\setbox0=\hbox{L}\hbox to\wd0{\hss\char'40L}}

\def\ep{\varepsilon}

\def\ka{\kappa}

\def\si{\sigma}

\def\P{{\Bbb P}}

\def\today{\ifcase\month\or
 January\or February\or March\or April\or May\or June\or
 July\or August\or September\or October\or November\or December\fi
 \space\number\day, \number\year}

\def\({\left(}
\def\){\right)}
\def\[{\left[}
\def\]{\right]}

\def\Re{\operatorname{Re}}

\def\sgn{\operatorname{sgn}}

\def\3{\ss}
\catcode`\@=11
\def\dddot#1{\vbox{\ialign{##\crcr
      .\hskip-.5pt.\hskip-.5pt.\crcr\noalign{\kern1.5\p@\nointerlineskip}
      $\hfil\displaystyle{#1}\hfil$\crcr}}}

\newif\iftab@\tab@false
\newif\ifvtab@\vtab@false
\def\tab{\bgroup\tab@true\vtab@false\vst@bfalse\Strich@false%
   \def\\{\global\hline@@false%
     \ifhline@\global\hline@false\global\hline@@true\fi\cr}
   \edef\l@{\the\leftskip}\ialign\bgroup\hskip\l@##\hfil&&##\hfil\cr}
\def\endtab{\cr\egroup\egroup}
\def\vtab{\vtop\bgroup\vst@bfalse\vtab@true\tab@true\Strich@false%
   \bgroup\def\\{\cr}\ialign\bgroup&##\hfil\cr}
\def\endvtab{\cr\egroup\egroup\egroup}
\def\stab{\D@cke0.5pt\null 
 \bgroup\tab@true\vtab@false\vst@bfalse\Strich@true\Let@@\vspace@
 \normalbaselines\offinterlineskip
  \openup\spreadmlines@
 \edef\l@{\the\leftskip}\ialign
 \bgroup\hskip\l@##\hfil&&##\hfil\crcr}
\def\endstab{\crcr\egroup
 \egroup}
\newif\ifvst@b\vst@bfalse
\def\vstab{\D@cke0.5pt\null
 \vtop\bgroup\tab@true\vtab@false\vst@btrue\Strich@true\bgroup\Let@@\vspace@
 \normalbaselines\offinterlineskip
  \openup\spreadmlines@\bgroup}
\def\endvstab{\crcr\egroup\egroup
 \egroup\tab@false\Strich@false}

\newdimen\htstrut@
\htstrut@8.5\p@
\newdimen\htStrut@
\htStrut@12\p@
\newdimen\dpstrut@
\dpstrut@3.5\p@
\newdimen\dpStrut@
\dpStrut@3.5\p@
\def\openup{\afterassignment\@penup\dimen@=}
\def\@penup{\advance\lineskip\dimen@
  \advance\baselineskip\dimen@
  \advance\lineskiplimit\dimen@
  \divide\dimen@ by2
  \advance\htstrut@\dimen@
  \advance\htStrut@\dimen@
  \advance\dpstrut@\dimen@
  \advance\dpStrut@\dimen@}
\def\Let@@{\relax%
    \def\\{\global\hline@@false%
     \ifhline@\global\hline@false\global\hline@@true\fi\cr}%
    \iffalse}\fi}
\def\matrix{\null\,\vcenter\bgroup
 \tab@false\vtab@false\vst@bfalse\Strich@false\Let@@\vspace@
 \normalbaselines\openup\spreadmlines@\ialign
 \bgroup\hfil$\m@th##$\hfil&&\quad\hfil$\m@th##$\hfil\crcr
 \Mathstrut@\crcr\noalign{\kern-\baselineskip}}
\def\endmatrix{\crcr\Mathstrut@\crcr\noalign{\kern-\baselineskip}\egroup
 \egroup\,}
\def\smatrix{\D@cke0.5pt\null\,
 \vcenter\bgroup\tab@false\vtab@false\vst@bfalse\Strich@true\Let@@\vspace@
 \normalbaselines\offinterlineskip
  \openup\spreadmlines@\ialign
 \bgroup\hfil$\m@th##$\hfil&&\quad\hfil$\m@th##$\hfil\crcr}
\def\endsmatrix{\crcr\egroup
 \egroup\,\Strich@false}
\newdimen\D@cke
\def\Dicke#1{\global\D@cke#1}
\newtoks\tabs@\tabs@{&}
\newif\ifStrich@\Strich@false
\newif\iff@rst

\def\Stricherr@{\iftab@\ifvtab@\errmessage{\noexpand\s not allowed
     here. Use \noexpand\vstab!}%
  \else\errmessage{\noexpand\s not allowed here. Use \noexpand\stab!}%
  \fi\else\errmessage{\noexpand\s not allowed
     here. Use \noexpand\smatrix!}\fi}
\def\format{\ifvst@b\else\crcr\fi\egroup\iffalse{\fi\ifnum`}=0 \fi\format@}
\def\format@#1\\{\def\preamble@{#1}%
 \def\Str@chfehlt##1{\ifx##1\s\Stricherr@\fi\ifx##1\\\let\Next\relax%
   \else\let\Next\Str@chfehlt\fi\Next}%
 \def\c{\hfil\noexpand\ifhline@@\hbox{\vrule height\htStrut@%
   depth\dpstrut@ width\z@}\noexpand\fi%
   \ifStrich@\hbox{\vrule height\htstrut@ depth\dpstrut@ width\z@}%
   \fi\iftab@\else$\m@th\fi\the\hashtoks@\iftab@\else$\fi\hfil}%
 \def\r{\hfil\noexpand\ifhline@@\hbox{\vrule height\htStrut@%
   depth\dpstrut@ width\z@}\noexpand\fi%
   \ifStrich@\hbox{\vrule height\htstrut@ depth\dpstrut@ width\z@}%
   \fi\iftab@\else$\m@th\fi\the\hashtoks@\iftab@\else$\fi}%
 \def\l{\noexpand\ifhline@@\hbox{\vrule height\htStrut@%
   depth\dpstrut@ width\z@}\noexpand\fi%
   \ifStrich@\hbox{\vrule height\htstrut@ depth\dpstrut@ width\z@}%
   \fi\iftab@\else$\m@th\fi\the\hashtoks@\iftab@\else$\fi\hfil}%
 \def\s{\ifStrich@\ \the\tabs@\vrule width\D@cke\the\hashtoks@%
          \fi\the\tabs@\ }%
 \def\sa{\ifStrich@\vrule width\D@cke\the\hashtoks@%
            \the\tabs@\ %
            \fi}%
 \def\se{\ifStrich@\ \the\tabs@\vrule width\D@cke\the\hashtoks@\fi}%
 \def\cd{\hfil\noexpand\ifhline@@\hbox{\vrule height\htStrut@%
   depth\dpstrut@ width\z@}\noexpand\fi%
   \ifStrich@\hbox{\vrule height\htstrut@ depth\dpstrut@ width\z@}%
   \fi$\dsize\m@th\the\hashtoks@$\hfil}%
 \def\rd{\hfil\noexpand\ifhline@@\hbox{\vrule height\htStrut@%
   depth\dpstrut@ width\z@}\noexpand\fi%
   \ifStrich@\hbox{\vrule height\htstrut@ depth\dpstrut@ width\z@}%
   \fi$\dsize\m@th\the\hashtoks@$}%
 \def\ld{\noexpand\ifhline@@\hbox{\vrule height\htStrut@%
   depth\dpstrut@ width\z@}\noexpand\fi%
   \ifStrich@\hbox{\vrule height\htstrut@ depth\dpstrut@ width\z@}%
   \fi$\dsize\m@th\the\hashtoks@$\hfil}%
 \ifStrich@\else\Str@chfehlt#1\\\fi%
 \setbox\z@\hbox{\xdef\Preamble@{\preamble@}}\ifnum`{=0 \fi\iffalse}\fi
 \ialign\bgroup\span\Preamble@\crcr}
\newif\ifhline@\hline@false
\newif\ifhline@@\hline@@false
\def\hlinefor#1{\multispan@{\strip@#1 }\leaders\hrule height\D@cke\hfill%
    \global\hline@true\ignorespaces}
\def\Item "#1"{\par\noindent\hangindent2\parindent%
  \hangafter1\setbox0\hbox{\rm#1\enspace}\ifdim\wd0>2\parindent%
  \box0\else\hbox to 2\parindent{\rm#1\hfil}\fi\ignorespaces}
\def\ITEM #1"#2"{\par\noindent\hangafter1\hangindent#1%
  \setbox0\hbox{\rm#2\enspace}\ifdim\wd0>#1%
  \box0\else\hbox to 0pt{\rm#2\hss}\hskip#1\fi\ignorespaces}
\def\item"#1"{\par\noindent\hang%
  \setbox0=\hbox{\rm#1\enspace}\ifdim\wd0>\the\parindent%
  \box0\else\hbox to \parindent{\rm#1\hfil}\enspace\fi\ignorespaces}
\let\plainitem@\item
\catcode`\@=13

\magnification1200
\hsize13cm
\vsize19cm
\newdimen\fullhsize
\newdimen\fullvsize
\newdimen\halfsize
\fullhsize13cm
\fullvsize19cm
\halfsize=0.5\fullhsize
\advance\halfsize by-0.5em

\catcode`\@=11
\font\tenln    = line10
\font\tenlnw   = linew10

\newskip\Einheit \Einheit=0.5cm
\newcount\xcoord \newcount\ycoord
\newdimen\xdim \newdimen\ydim \newdimen\PfadD@cke \newdimen\Pfadd@cke

\newcount\@tempcnta
\newcount\@tempcntb

\newdimen\@tempdima
\newdimen\@tempdimb

\newdimen\@wholewidth
\newdimen\@halfwidth

\newcount\@xarg
\newcount\@yarg
\newcount\@yyarg
\newbox\@linechar
\newbox\@tempboxa
\newdimen\@linelen
\newdimen\@clnwd
\newdimen\@clnht

\newif\if@negarg

\def\@whilenoop#1{}
\def\@whiledim#1\do #2{\ifdim #1\relax#2\@iwhiledim{#1\relax#2}\fi}
\def\@iwhiledim#1{\ifdim #1\let\@nextwhile=\@iwhiledim
        \else\let\@nextwhile=\@whilenoop\fi\@nextwhile{#1}}

\def\@whileswnoop#1\fi{}
\def\@whilesw#1\fi#2{#1#2\@iwhilesw{#1#2}\fi\fi}
\def\@iwhilesw#1\fi{#1\let\@nextwhile=\@iwhilesw
         \else\let\@nextwhile=\@whileswnoop\fi\@nextwhile{#1}\fi}

\def\thinlines{\let\@linefnt\tenln \let\@circlefnt\tencirc
  \@wholewidth\fontdimen8\tenln \@halfwidth .5\@wholewidth}
\def\thicklines{\let\@linefnt\tenlnw \let\@circlefnt\tencircw
  \@wholewidth\fontdimen8\tenlnw \@halfwidth .5\@wholewidth}
\thinlines

\PfadD@cke1pt \Pfadd@cke0.5pt
\def\PfadDicke#1{\PfadD@cke#1 \divide\PfadD@cke by2 \Pfadd@cke\PfadD@cke \multiply\PfadD@cke by2}
\long\def\LOOP#1\REPEAT{\def\BODY{#1}\ITERATE}
\def\ITERATE{\BODY \let\next\ITERATE \else\let\next\relax\fi \next}
\let\REPEAT=\fi
\def\Punkt{\hbox{\raise-2pt\hbox to0pt{\hss$\ssize\bullet$\hss}}}
\def\DuennPunkt(#1,#2){\unskip
  \raise#2 \Einheit\hbox to0pt{\hskip#1 \Einheit
          \raise-2.5pt\hbox to0pt{\hss$\bullet$\hss}\hss}}
\def\NormalPunkt(#1,#2){\unskip
  \raise#2 \Einheit\hbox to0pt{\hskip#1 \Einheit
          \raise-3pt\hbox to0pt{\hss\twelvepoint$\bullet$\hss}\hss}}
\def\DickPunkt(#1,#2){\unskip
  \raise#2 \Einheit\hbox to0pt{\hskip#1 \Einheit
          \raise-4pt\hbox to0pt{\hss\fourteenpoint$\bullet$\hss}\hss}}
\def\Kreis(#1,#2){\unskip
  \raise#2 \Einheit\hbox to0pt{\hskip#1 \Einheit
          \raise-4pt\hbox to0pt{\hss\fourteenpoint$\circ$\hss}\hss}}

\def\Line@(#1,#2)#3{\@xarg #1\relax \@yarg #2\relax
\@linelen=#3\Einheit
\ifnum\@xarg =0 \@vline
  \else \ifnum\@yarg =0 \@hline \else \@sline\fi
\fi}

\def\@sline{\ifnum\@xarg< 0 \@negargtrue \@xarg -\@xarg \@yyarg -\@yarg
  \else \@negargfalse \@yyarg \@yarg \fi
\ifnum \@yyarg >0 \@tempcnta\@yyarg \else \@tempcnta -\@yyarg \fi
\ifnum\@tempcnta>6 \@badlinearg\@tempcnta0 \fi
\ifnum\@xarg>6 \@badlinearg\@xarg 1 \fi
\setbox\@linechar\hbox{\@linefnt\@getlinechar(\@xarg,\@yyarg)}%
\ifnum \@yarg >0 \let\@upordown\raise \@clnht\z@
   \else\let\@upordown\lower \@clnht \ht\@linechar\fi
\@clnwd=\wd\@linechar
\if@negarg \hskip -\wd\@linechar \def\@tempa{\hskip -2\wd\@linechar}\else
     \let\@tempa\relax \fi
\@whiledim \@clnwd <\@linelen \do
  {\@upordown\@clnht\copy\@linechar
   \@tempa
   \advance\@clnht \ht\@linechar
   \advance\@clnwd \wd\@linechar}%
\advance\@clnht -\ht\@linechar
\advance\@clnwd -\wd\@linechar
\@tempdima\@linelen\advance\@tempdima -\@clnwd
\@tempdimb\@tempdima\advance\@tempdimb -\wd\@linechar
\if@negarg \hskip -\@tempdimb \else \hskip \@tempdimb \fi
\multiply\@tempdima \@m
\@tempcnta \@tempdima \@tempdima \wd\@linechar \divide\@tempcnta \@tempdima
\@tempdima \ht\@linechar \multiply\@tempdima \@tempcnta
\divide\@tempdima \@m
\advance\@clnht \@tempdima
\ifdim \@linelen <\wd\@linechar
   \hskip \wd\@linechar
  \else\@upordown\@clnht\copy\@linechar\fi}

\def\@hline{\ifnum \@xarg <0 \hskip -\@linelen \fi
\vrule height\Pfadd@cke width \@linelen depth\Pfadd@cke
\ifnum \@xarg <0 \hskip -\@linelen \fi}

\def\@getlinechar(#1,#2){\@tempcnta#1\relax\multiply\@tempcnta 8
\advance\@tempcnta -9 \ifnum #2>0 \advance\@tempcnta #2\relax\else
\advance\@tempcnta -#2\relax\advance\@tempcnta 64 \fi
\char\@tempcnta}

\def\Vektor(#1,#2)#3(#4,#5){\unskip\leavevmode
  \xcoord#4\relax \ycoord#5\relax
      \raise\ycoord \Einheit\hbox to0pt{\hskip\xcoord \Einheit
         \Vector@(#1,#2){#3}\hss}}

\def\Vector@(#1,#2)#3{\@xarg #1\relax \@yarg #2\relax
\@tempcnta \ifnum\@xarg<0 -\@xarg\else\@xarg\fi
\ifnum\@tempcnta<5\relax
\@linelen=#3\Einheit
\ifnum\@xarg =0 \@vvector
  \else \ifnum\@yarg =0 \@hvector \else \@svector\fi
\fi
\else\@badlinearg\fi}

\def\@hvector{\@hline\hbox to 0pt{\@linefnt
\ifnum \@xarg <0 \@getlarrow(1,0)\hss\else
    \hss\@getrarrow(1,0)\fi}}

\def\@vvector{\ifnum \@yarg <0 \@downvector \else \@upvector \fi}

\def\@svector{\@sline
\@tempcnta\@yarg \ifnum\@tempcnta <0 \@tempcnta=-\@tempcnta\fi
\ifnum\@tempcnta <5
  \hskip -\wd\@linechar
  \@upordown\@clnht \hbox{\@linefnt  \if@negarg
  \@getlarrow(\@xarg,\@yyarg) \else \@getrarrow(\@xarg,\@yyarg) \fi}%
\else\@badlinearg\fi}

\def\@upline{\hbox to \z@{\hskip -.5\Pfadd@cke \vrule width \Pfadd@cke
   height \@linelen depth \z@\hss}}

\def\@downline{\hbox to \z@{\hskip -.5\Pfadd@cke \vrule width \Pfadd@cke
   height \z@ depth \@linelen \hss}}

\def\@upvector{\@upline\setbox\@tempboxa\hbox{\@linefnt\char'66}\raise
     \@linelen \hbox to\z@{\lower \ht\@tempboxa\box\@tempboxa\hss}}

\def\@downvector{\@downline\lower \@linelen
      \hbox to \z@{\@linefnt\char'77\hss}}

\def\@getlarrow(#1,#2){\ifnum #2 =\z@ \@tempcnta='33\else
\@tempcnta=#1\relax\multiply\@tempcnta \sixt@@n \advance\@tempcnta
-9 \@tempcntb=#2\relax\multiply\@tempcntb \tw@
\ifnum \@tempcntb >0 \advance\@tempcnta \@tempcntb\relax
\else\advance\@tempcnta -\@tempcntb\advance\@tempcnta 64
\fi\fi\char\@tempcnta}

\def\@getrarrow(#1,#2){\@tempcntb=#2\relax
\ifnum\@tempcntb < 0 \@tempcntb=-\@tempcntb\relax\fi
\ifcase \@tempcntb\relax \@tempcnta='55 \or
\ifnum #1<3 \@tempcnta=#1\relax\multiply\@tempcnta
24 \advance\@tempcnta -6 \else \ifnum #1=3 \@tempcnta=49
\else\@tempcnta=58 \fi\fi\or
\ifnum #1<3 \@tempcnta=#1\relax\multiply\@tempcnta
24 \advance\@tempcnta -3 \else \@tempcnta=51\fi\or
\@tempcnta=#1\relax\multiply\@tempcnta
\sixt@@n \advance\@tempcnta -\tw@ \else
\@tempcnta=#1\relax\multiply\@tempcnta
\sixt@@n \advance\@tempcnta 7 \fi\ifnum #2<0 \advance\@tempcnta 64 \fi
\char\@tempcnta}

\def\Diagonale(#1,#2)#3{\unskip\leavevmode
  \xcoord#1\relax \ycoord#2\relax
      \raise\ycoord \Einheit\hbox to0pt{\hskip\xcoord \Einheit
         \Line@(1,1){#3}\hss}}
\def\AntiDiagonale(#1,#2)#3{\unskip\leavevmode
  \xcoord#1\relax \ycoord#2\relax 
      \raise\ycoord \Einheit\hbox to0pt{\hskip\xcoord \Einheit
         \Line@(1,-1){#3}\hss}}
\def\Pfad(#1,#2),#3\endPfad{\unskip\leavevmode
  \xcoord#1 \ycoord#2 \thicklines\ZeichnePfad#3\endPfad\thinlines}
\def\ZeichnePfad#1{\ifx#1\endPfad\let\next\relax
  \else\let\next\ZeichnePfad
    \ifnum#1=1
      \raise\ycoord \Einheit\hbox to0pt{\hskip\xcoord \Einheit
         \vrule height\Pfadd@cke width1 \Einheit depth\Pfadd@cke\hss}%
      \advance\xcoord by 1
    \else\ifnum#1=2
      \raise\ycoord \Einheit\hbox to0pt{\hskip\xcoord \Einheit
        \hbox{\hskip-\PfadD@cke\vrule height1 \Einheit width\PfadD@cke depth0pt}\hss}%
      \advance\ycoord by 1
    \else\ifnum#1=3
      \raise\ycoord \Einheit\hbox to0pt{\hskip\xcoord \Einheit
         \Line@(1,1){1}\hss}
      \advance\xcoord by 1
      \advance\ycoord by 1
    \else\ifnum#1=4
      \raise\ycoord \Einheit\hbox to0pt{\hskip\xcoord \Einheit
         \Line@(1,-1){1}\hss}
      \advance\xcoord by 1
      \advance\ycoord by -1
    \else\ifnum#1=5
      \advance\xcoord by -1
      \raise\ycoord \Einheit\hbox to0pt{\hskip\xcoord \Einheit
         \vrule height\Pfadd@cke width1 \Einheit depth\Pfadd@cke\hss}%
    \else\ifnum#1=6
      \advance\ycoord by -1
      \raise\ycoord \Einheit\hbox to0pt{\hskip\xcoord \Einheit
        \hbox{\hskip-\PfadD@cke\vrule height1 \Einheit width\PfadD@cke depth0pt}\hss}%
    \else\ifnum#1=7
      \advance\xcoord by -1
      \advance\ycoord by -1
      \raise\ycoord \Einheit\hbox to0pt{\hskip\xcoord \Einheit
         \Line@(1,1){1}\hss}
    \else\ifnum#1=8
      \advance\xcoord by -1
      \advance\ycoord by +1
      \raise\ycoord \Einheit\hbox to0pt{\hskip\xcoord \Einheit
         \Line@(1,-1){1}\hss}
    \fi\fi\fi\fi
    \fi\fi\fi\fi
  \fi\next}
\def\hSSchritt{\leavevmode\raise-.4pt\hbox to0pt{\hss.\hss}\hskip.2\Einheit
  \raise-.4pt\hbox to0pt{\hss.\hss}\hskip.2\Einheit
  \raise-.4pt\hbox to0pt{\hss.\hss}\hskip.2\Einheit
  \raise-.4pt\hbox to0pt{\hss.\hss}\hskip.2\Einheit
  \raise-.4pt\hbox to0pt{\hss.\hss}\hskip.2\Einheit}
\def\vSSchritt{\vbox{\baselineskip.2\Einheit\lineskiplimit0pt
\hbox{.}\hbox{.}\hbox{.}\hbox{.}\hbox{.}}}
\def\DSSchritt{\leavevmode\raise-.4pt\hbox to0pt{%
  \hbox to0pt{\hss.\hss}\hskip.2\Einheit
  \raise.2\Einheit\hbox to0pt{\hss.\hss}\hskip.2\Einheit
  \raise.4\Einheit\hbox to0pt{\hss.\hss}\hskip.2\Einheit
  \raise.6\Einheit\hbox to0pt{\hss.\hss}\hskip.2\Einheit
  \raise.8\Einheit\hbox to0pt{\hss.\hss}\hss}}
\def\dSSchritt{\leavevmode\raise-.4pt\hbox to0pt{%
  \hbox to0pt{\hss.\hss}\hskip.2\Einheit
  \raise-.2\Einheit\hbox to0pt{\hss.\hss}\hskip.2\Einheit
  \raise-.4\Einheit\hbox to0pt{\hss.\hss}\hskip.2\Einheit
  \raise-.6\Einheit\hbox to0pt{\hss.\hss}\hskip.2\Einheit
  \raise-.8\Einheit\hbox to0pt{\hss.\hss}\hss}}
\def\SPfad(#1,#2),#3\endSPfad{\unskip\leavevmode
  \xcoord#1 \ycoord#2 \ZeichneSPfad#3\endSPfad}
\def\ZeichneSPfad#1{\ifx#1\endSPfad\let\next\relax
  \else\let\next\ZeichneSPfad
    \ifnum#1=1
      \raise\ycoord \Einheit\hbox to0pt{\hskip\xcoord \Einheit
         \hSSchritt\hss}%
      \advance\xcoord by 1
    \else\ifnum#1=2
      \raise\ycoord \Einheit\hbox to0pt{\hskip\xcoord \Einheit
        \hbox{\hskip-2pt \vSSchritt}\hss}%
      \advance\ycoord by 1
    \else\ifnum#1=3
      \raise\ycoord \Einheit\hbox to0pt{\hskip\xcoord \Einheit
         \DSSchritt\hss}
      \advance\xcoord by 1
      \advance\ycoord by 1
    \else\ifnum#1=4
      \raise\ycoord \Einheit\hbox to0pt{\hskip\xcoord \Einheit
         \dSSchritt\hss}
      \advance\xcoord by 1
      \advance\ycoord by -1
    \else\ifnum#1=5
      \advance\xcoord by -1
      \raise\ycoord \Einheit\hbox to0pt{\hskip\xcoord \Einheit
         \hSSchritt\hss}%
    \else\ifnum#1=6
      \advance\ycoord by -1
      \raise\ycoord \Einheit\hbox to0pt{\hskip\xcoord \Einheit
        \hbox{\hskip-2pt \vSSchritt}\hss}%
    \else\ifnum#1=7
      \advance\xcoord by -1
      \advance\ycoord by -1
      \raise\ycoord \Einheit\hbox to0pt{\hskip\xcoord \Einheit
         \DSSchritt\hss}
    \else\ifnum#1=8
      \advance\xcoord by -1
      \advance\ycoord by 1
      \raise\ycoord \Einheit\hbox to0pt{\hskip\xcoord \Einheit
         \dSSchritt\hss}
    \fi\fi\fi\fi
    \fi\fi\fi\fi
  \fi\next}
\def\Koordinatenachsen(#1,#2){\unskip
 \hbox to0pt{\hskip-.5pt\vrule height#2 \Einheit width.5pt depth1 \Einheit}%
 \hbox to0pt{\hskip-1 \Einheit \xcoord#1 \advance\xcoord by1
    \vrule height0.25pt width\xcoord \Einheit depth0.25pt\hss}}
\def\Koordinatenachsen(#1,#2)(#3,#4){\unskip
 \hbox to0pt{\hskip-.5pt \ycoord-#4 \advance\ycoord by1
    \vrule height#2 \Einheit width.5pt depth\ycoord \Einheit}%
 \hbox to0pt{\hskip-1 \Einheit \hskip#3\Einheit 
    \xcoord#1 \advance\xcoord by1 \advance\xcoord by-#3 
    \vrule height0.25pt width\xcoord \Einheit depth0.25pt\hss}}
\def\Gitter(#1,#2){\unskip \xcoord0 \ycoord0 \leavevmode
  \LOOP\ifnum\ycoord<#2
    \loop\ifnum\xcoord<#1
      \raise\ycoord \Einheit\hbox to0pt{\hskip\xcoord \Einheit\Punkt\hss}%
      \advance\xcoord by1
    \repeat
    \xcoord0
    \advance\ycoord by1
  \REPEAT}
\def\Gitter(#1,#2)(#3,#4){\unskip \xcoord#3 \ycoord#4 \leavevmode
  \LOOP\ifnum\ycoord<#2
    \loop\ifnum\xcoord<#1
      \raise\ycoord \Einheit\hbox to0pt{\hskip\xcoord \Einheit\Punkt\hss}%
      \advance\xcoord by1
    \repeat
    \xcoord#3
    \advance\ycoord by1
  \REPEAT}
\def\Label#1#2(#3,#4){\unskip \xdim#3 \Einheit \ydim#4 \Einheit
  \def\lo{\advance\xdim by-.5 \Einheit \advance\ydim by.5 \Einheit}%
  \def\llo{\advance\xdim by-.25cm \advance\ydim by.5 \Einheit}%
  \def\loo{\advance\xdim by-.5 \Einheit \advance\ydim by.25cm}%
  \def\o{\advance\ydim by.25cm}%
  \def\ro{\advance\xdim by.5 \Einheit \advance\ydim by.5 \Einheit}%
  \def\rro{\advance\xdim by.25cm \advance\ydim by.5 \Einheit}%
  \def\roo{\advance\xdim by.5 \Einheit \advance\ydim by.25cm}%
  \def\l{\advance\xdim by-.30cm}%
  \def\r{\advance\xdim by.30cm}%
  \def\lu{\advance\xdim by-.5 \Einheit \advance\ydim by-.6 \Einheit}%
  \def\llu{\advance\xdim by-.25cm \advance\ydim by-.6 \Einheit}%
  \def\luu{\advance\xdim by-.5 \Einheit \advance\ydim by-.30cm}%
  \def\u{\advance\ydim by-.30cm}%
  \def\ru{\advance\xdim by.5 \Einheit \advance\ydim by-.6 \Einheit}%
  \def\rru{\advance\xdim by.25cm \advance\ydim by-.6 \Einheit}%
  \def\ruu{\advance\xdim by.5 \Einheit \advance\ydim by-.30cm}%
  #1\raise\ydim\hbox to0pt{\hskip\xdim
     \vbox to0pt{\vss\hbox to0pt{\hss$#2$\hss}\vss}\hss}%
}
\catcode`\@=13

\TagsOnRight

\def\ZeilAV{55}
\def\ZeilAM{54}
\def\WJZeAA{53}
\def\WhWaAA{52}
\def\SulaAC{51}
\def\SlatAC{50}
\def\SlatZZ{49}
\def\SlatZY{48}
\def\SagaAL{47}
\def\RubMAG{46}
\def\RubMAC{45}
\def\RaVeAA{44}
\def\PeWZAA{43}
\def\PaScAA{42}
\def\OwEBAA{41}
\def\MohaAE{40}
\def\MohaAB{39}
\def\LindAA{38}
\def\KrGVAA{37}
\def\KratBN{36}
\def\KratBK{35}
\def\KratAP{34}
\def\KaMGAC{33}
\def\KaMGAB{32}
\def\KaJGAA{31}
\def\JoSaAB{30}
\def\GuOVAA{29}
\def\GrJSAA{28}
\def\GrKPAA{27}
\def\GospAB{26}
\def\GeViAB{25}
\def\GeViAA{24}
\def\GaRaAA{23}
\def\ForrAD{22}
\def\ForrAC{21}
\def\ForrAB{20}
\def\FlSeAA{19}
\def\FlOdAA{18}
\def\FishAA{17}
\def\FiscAA{16}
\def\EsGuAA{15}
\def\EngbAB{14}
\def\CsViAA{13}
\def\ComtAA{12}
\def\CiKrAB{11}
\def\CiEKAA{10}
\def\BrEOAD{9}
\def\BrEOAC{8}
\def\BrEsAC{7}
\def\BrEsAB{6}
\def\BrEsAA{5}
\def\BailAA{4}
\def\BiTrAA{3}
\def\BirkAA{2}
\def\ArMEAA{1}

\def\PP{1.1}
\def\ZY{2.1}
\def\ZZ{2.2}
\def\ZZa{2.3}
\def\ZT{2.4}
\def\ZU{2.5}
\def\ZV{2.6}
\def\ZS{2.7}
\def\ZI{2.8}
\def\ZW{2.9}
\def\ZWa{2.10}
\def\ZX{2.11}
\def\Za{3.1}
\def\Zb{3.2}
\def\Zc{3.3}
\def\ZB{3.4}
\def\ZA{3.5}
\def\ZBa{3.6}
\def\ZAa{3.7}
\def\ZG{4.1}
\def\ZD{4.2}
\def\ZC{4.3}
\def\ZE{4.4}
\def\ZF{4.5}
\def\ZN{4.6}
\def\ZNa{4.7}
\def\Zd{4.8}
\def\Zdd{4.9}
\def\Zddd{4.10}
\def\Zg{5.1}
\def\Zi{5.2}
\def\Zh{5.3}
\def\Zl{5.4}
\def\Zk{5.5}
\def\Zm{6.1}
\def\Zn{6.2}
\def\Zo{6.3}
\def\Zp{6.4}
\def\Zpa{6.5}
\def\Zq{6.6}
\def\ZP{6.7}
\def\ZPa{6.8}
\def\ZL{6.9}
\def\ZM{6.10}
\def\Zr{6.11}
\def\Zra{6.12}
\def\Zs{6.13}
\def\Zt{6.14}
\def\Zu{6.15}
\def\AB{7.1}
\def\AC{7.2}
\def\AD{7.3}
\def\ZK{7.4}
\def\ZKa{7.5}
\def\ZH{7.6}
\def\ZHa{7.7}
\def\ZHc{7.8}
\def\ZJ{7.9}
\def\Zv{7.10}
\def\Zw{7.11}
\def\ZKb{7.12}
\def\Zff{8.1}
\def\Zf{8.2}
\def\ZO{8.3}
\def\ZOa{8.4}
\def\ZQa{8.5}
\def\ZQb{8.6}
\def\ZQc{8.7}
\def\ZR{8.8}
\def\Zx{8.9}

\def\TA{1}
\def\TH{2}
\def\TB{3}
\def\TC{4}
\def\TCa{5}
\def\TI{6}
\def\TJ{7}
\def\TK{8}
\def\TL{9}
\def\TD{10}
\def\TN{11}
\def\TE{12}
\def\TF{13}
\def\TFa{14}
\def\TG{15}

\def\FA{1}
\def\FB{2}
\def\FC{3}
\def\FD{4}
\def\FE{5}
\def\FF{6}
\def\FG{7}

\def\vv#1{\left\vert#1\right\vert}
\def\po#1#2{(#1)_#2}
\def\P{{\Cal P}}
\def\GF{\operatorname{GF}}

\topmatter 
\title Watermelon configurations with wall interaction: 
exact and asymptotic results
\endtitle 
\author C.~Krattenthaler$^\dagger$
\endauthor 
\affil 
Institut Girard 
Desargues, Universit\'e Claude Bernard Lyon-I,\\
21, avenue Claude Bernard, F-69622 Villeurbanne Cedex, France.\\
e-mail: kratt\@euler.univ-lyon1.fr\\
WWW: \tt http://igd.univ-lyon1.fr/\~{}kratt
\endaffil
\address Institut Girard 
Desargues, Universit\'e Claude Bernard Lyon-I, 
21, avenue Claude Bernard, F-69622 Villeurbanne Cedex, France.
\endaddress
\thanks{$^\dagger$Research partially supported by the Austrian
Science Foundation FWF, grants P12094-MAT and P13190-MAT, 
and by EC's IHRP Programme,
grant HPRN-CT-2001-00272, ``Algebraic Combinatorics in Europe"}\endthanks
\subjclass Primary 82B20;
 Secondary 05A15 05A16 60G50 82B23 82B26 82B41 82B43
\endsubjclass
\keywords vicious walkers, watermelon configurations, non-intersecting
lattice\break paths, touchings, reflection principle, hypergeometric
series, singularity analysis
\endkeywords
\abstract 
We perform an exact and asymptotic analysis of the model of $n$
vicious walkers interacting with a wall via contact potentials, a
model introduced by Brak, Essam and Owczarek. More specifically, we
study the partition function of watermelon configurations 
which start on the wall, but may end at arbitrary height, 
and their mean number of contacts with the
wall. We improve and extend the earlier (partially non-rigorous) results by
Brak, Essam and Owczarek, providing new exact results, and more
precise and more general asymptotic results, 
in particular full asymptotic expansions
for the partition function and the mean number of contacts.
Furthermore, we relate this circle of problems to earlier results in
the combinatorial and statistical literature.
\endabstract
\endtopmatter
\document

\subhead 1. Introduction\endsubhead
The problem of vicious walkers was introduced by Fisher \cite{\FishAA}, 
who also gave a number of physical applications of the model, such as, for
example, to modelling wetting and melting. 
The general model is one of $n$ random walkers on a 
$d$-dimensional lattice who at 
regular time intervals simultaneously take one step 
in the direction of one of the allowed lattice vectors
such that at no time two walkers occupy the same lattice site. 

Numerous papers have been written on the subject since then.
Most of them analyse the model of vicious walkers in a continuum
limit (such as for example 
\cite{\FishAA, \ForrAB, \ForrAC, \ForrAD}). It has been
realized only recently that in fact there are many interesting cases in
which even exact results in form of nice closed product formulas are available
(see for example \cite{\ArMEAA, \BrEsAC, 
\BrEOAC, \BrEOAD, \EsGuAA, \GuOVAA, 
\KrGVAA, \OwEBAA, \RubMAC}), and that asymptotic analysis can be performed
directly on the model, without taking recourse to continuum limits,
thus obtaining more precise estimates (see for example 
\cite{\BrEOAC, \EsGuAA, \KrGVAA, \OwEBAA, \RubMAC}).

In this paper
we consider $n$ vicious walkers 
(the ``vicious"
constraint demanding that at no time two walkers occupy the same site)
in the plane integer lattice with
allowed steps of the form $(1,1)$ ({\it up-step}) and $(1,-1)$ ({\it
down-step}),
which in addition do not run below the $x$-axis (the {\it wall}). 
Given integer vectors $\bold a=(a_1,a_2,\dots,a_n)$ and $\bold
e=(e_1,e_2,\dots,e_n)$, where all $a_i$'s are of the same parity and
all $e_i$'s are of the same parity, we study the partition function
$Z^{(n)}_t(\bold a\to\bold e;\ka):=\sum _{\bold P} ^{}\ka^{c(\bold
P)}$, where the sum is over all families $\bold P=(P_1,P_2,\dots,P_n)$
of $n$ vicious walkers as above, the $i$-th walker $P_i$ starting at
$(0,a_i)$ and ending at $(t,e_i)$, $i=1,2,\dots,n$, and with $c(\bold
P)$ denoting the total number of contacts of the walkers with 
the wall. (Because of the ``vicious" constraint, it is only the
lowest of the walkers who can have contacts with the wall.) 
An example of such a family of vicious walkers with $n=4$,
$\bold a=(0,2,4,6)$, $\bold e=(0,2,4,6)$, $t=8$ is shown in
Figure~\FA.a, whereas Figure~\FA.b shows a 
family of vicious walkers with $n=5$,
$\bold a=(0,2,4,6,8)$, $\bold e=(2,4,6,8,10)$, $t=12$.
The number of contacts with the wall (the $x$-axis) is 3
in Figure~\FA.a, and it is 4 in Figure~\FA.b.

\midinsert
\vskip10pt
\vbox{
$$
\Einheit.4cm
        \Gitter(9,12)(0,-1)
        \Koordinatenachsen(9,12)(0,-1)
        \Pfad(0,0),33434434\endPfad
        \Pfad(0,2),33344344\endPfad
        \Pfad(0,4),33343444\endPfad
        \Pfad(0,6),33334444\endPfad
\PfadDicke{.5pt}
        \Pfad(8,-2),22222222222222\endPfad
        \DickPunkt(0,0)
        \DickPunkt(0,2)
        \DickPunkt(0,4)
        \DickPunkt(0,6)
        \DickPunkt(8,0)
        \DickPunkt(8,2)
        \DickPunkt(8,4)
        \DickPunkt(8,6)
        \Label\l{u_1}(0,0)
        \Label\l{u_2}(0,2)
        \Label\l{u_3}(0,4)
        \Label\l{u_4}(0,6)
        \Label\r{v_1}(8,0)
        \Label\r{v_2}(8,2)
        \Label\r{v_3}(8,4)
        \Label\r{v_4}(8,6)
\hbox{\hskip6cm}
\Einheit.4cm
        \Gitter(13,14)(0,-1)
        \Koordinatenachsen(13,14)(0,-1)
\Pfad(0,0),343343443433\endPfad
\Pfad(0,2),334333444433\endPfad
\Pfad(0,4),334333443434\endPfad
\Pfad(0,6),333343344344\endPfad
\Pfad(0,8),333343343444\endPfad
\DickPunkt(0,0)
\DickPunkt(0,2)
\DickPunkt(0,4)
\DickPunkt(0,6)
\DickPunkt(0,8)
\DickPunkt(12,2)
\DickPunkt(12,4)
\DickPunkt(12,6)
\DickPunkt(12,8)
\DickPunkt(12,10)
        \Label\l{u_1}(0,0)
        \Label\l{u_2}(0,2)
        \Label\l{u_3}(0,4)
        \Label\l{u_4}(0,6)
        \Label\l{u_5}(0,8)
        \Label\r{v_1}(12,2)
        \Label\r{v_2}(12,4)
        \Label\r{v_3}(12,6)
        \Label\r{v_4}(12,8)
        \Label\r{v_5}(12,10)
\hbox{\hskip4.8cm}
$$
\centerline{\eightpoint a.\hskip6cm b.}
\centerline{\eightpoint Vicious walkers}
\vskip7pt
\centerline{\eightpoint Figure \FA}
}
\vskip10pt
\endinsert

The special topology of vicious walkers that we shall be mostly
concerned with in this paper is {\it watermelon configurations} ({\it
watermelons} for short), which are families of vicious walkers in
which the $y$-coordinates of neighbouring starting points differ by
2, the same being true for the end points. The $y$-coordinate of the
end point of the lowest walker is called the {\it deviation} of the
watermelon configuration. Thus, both vicious walker families in
Figure~\FA\ are in fact watermelon configurations, the one in
Figure~\FA.a has deviation 0, the one in Figure~\FA.b has deviation
2.

The above described vicious walker model with contact interaction
has been introduced by Owczarek, Essam and Brak in
\cite{\OwEBAA}. They obtained a determinantal formula for the
partition function, and undertook a 
(partially non-rigorous\footnote{Their approach could probably
be made rigorous by applying the Birkhoff--Trjitzinsky theory
\cite{\BirkAA, \BiTrAA} of
determining the asymptotics of solutions to difference equations,
surveyed in \cite{\WJZeAA}.}) study of
its asymptotic properties as the length $t$ of the walks tends to
infinity in the case of watermelon configurations starting and ending
on the wall. (By definition, 
these arise if $\bold a=\bold e=(0,2,4,\dots,2n-2)$.) 
This study was crucially based on recurrence relations for the
partition function. It revealed a critical value of
the parameter $\ka$ at $\ka=2$, and they also provided a scaling
analysis around this critical value. However, the use of recurrence
relations intrinsically makes it only possible to obtain the order of
magnitude of the partition function, but not the multiplicative
constant, not to mention any error terms. This paper was followed by
the articles \cite{\BrEsAB, \BrEsAC} by Brak and
Essam, in which exact results for the
partition function for watermelons which start on the wall
and end at some fixed deviation $y$ not necessarily equal to 0 are
obtained. The latter three papers do on the other hand not contain any
asymptotic results for this more general situation.

The purpose of our paper is three-fold. First, we relate 
this circle of problems to earlier results in
the combinatorial and statistical literature. In particular, we show
that the determinantal formulas in \cite{\BrEOAD} and \cite{\OwEBAA}
follow directly from a (now) classical result on non-intersecting
path families in acyclic graphs due to Lindstr\"om \cite{\LindAA}, 
which was
rediscovered by Gessel and Viennot \cite{\GeViAA, \GeViAB}. (In fact, {\it
all\/} determinantal formulas for vicious walkers follow from that
result.) Furthermore, we
outline that the results in \cite{\BrEOAC} and \cite{\OwEBAA} on the
partition function of a single walker are well-known in nonparametric
statistics. Second, we derive a new exact formula for the partition
function 
$$Z^{(n)}_t(y;\ka):=
Z^{(n)}_t\big((0,2,\dots,2n-2)\to(y,y+2,\dots,y+2n-2);\ka\big)\tag\PP$$
for watermelon configurations which start on the wall and
end at deviation $y$ (see Theorem~\TK), which expresses
$Z^{(n)}_t(y;\ka)$ in terms of a double sum. This is the central
result of our paper, from which all other results are derived. Not
only does it allow us to rederive, in a uniform manner, all
previously obtained exact results \cite{\BrEsAB, \BrEsAC, \BrEOAC, 
\CiKrAB, \KrGVAA, \OwEBAA} on the partition function 
$Z^{(n)}_t(y;\ka)$ in the literature (see Corollary~\TCa,
Theorem~\TL, Corollary~\TD, Theorem~\TN), it also enables us to
extend these in two cases to significantly larger domains of the
parameter $\ka$ (see Corollary~\TCa\ and Theorem~\TN).
Third, we improve and extend the asymptotic results by 
Owczarek, Essam and Brak \cite{\OwEBAA}. Not only do we show how to
rigorously find the asymptotic form of the partition function 
for watermelon
configurations which start on the wall and end at an arbitrary
(fixed, or non-fixed) deviation as
the length of the walks tends to infinity, our
approach allows to even derive full asymptotic expansions. In
particular, our results confirm the phase transition at $\ka=2$
predicted by Owczarek, Essam and Brak.
Moreover, we provide new exact and asymptotic results for the
(normalized) mean
number of contacts of the watermelons with the wall. 
These results show that, for the length of the walks being large, 
the normalized mean number of contacts is (asymptotically)
proportional to a constant if $\ka$ is less than the critical value
$2$, it is proportional to the
square root of the length of the walks if $\ka$ is equal to the
critical value, and it is
proportional to the length of the walks if $\ka$ is greater than $2$.
Again, all these
results are made by possible by our double sum formula for the
partition function given in Theorem~\TK.

The techniques that we employ to prove our results are 
(1)~combinatorial: path manipulation, Lindstr\"om--Gessel--Viennot
theorem on non-intersecting lattice paths, tableau combinatorics (in
particular: jeu de taquin); (2)~algebraic-manipulatory: 
hypergeometric series identities;\footnote{All the hypergeometric
calculations in this paper were carried out using the author's 
{\sl Mathematica\/} package HYP, which is designed for a
convenient handling of hypergeometric series, and is available from
{\tt http://www.mat.univie.ac.at/\~{}kratt}.} and for
the asymptotic calculations (3) analytic: singularity analysis.
In particular, our approach clearly shows how the various 
results that have been obtained earlier are connected.

Our paper is organised as follows. In the next section we address the
analysis of a single walker with wall interaction, and relate it to
the statistical literature, in particular to papers by Engelberg 
\cite{\EngbAB} and Mohanty \cite{\MohaAB}. Then, in
Section~3, we address the case of several walkers.
We recall the Lindstr\"om--Gessel--Viennot theorem and
demonstrate how it directly implies the determinantal formula 
\cite{\BrEOAD} by Brak, Essam and Owczarek. In addition,
in Proposition~\TB,
we also deduce a slightly different determinantal formula for the
partition function for watermelon configurations which start  
on the wall, which will be more convenient for our subsequent computations. 
In Section~4 we restrict our attention to the exact enumeration of
watermelon configurations which start {\it and end\/} on the wall.
We provide a new, short proof for a result by Brak and Essam \cite{\BrEsAB,
Theorem~6}, which expresses the partition
function in form of a single sum. This proof shows in particular that 
{\it guessing} the result is sufficient to prove it, the details being
filled in by the computer (see Theorem~\TC\ and its proof). Moreover,
we reprove and extend an alternative expression for the partition function found
earlier by Owczarek, Essam and Brak \cite{\OwEBAA, Eq.~(4.65)} (see
Corollary~\TCa).
Section~5 begins the analysis of
watermelon configurations which start on
the wall and end at some fixed deviation $y$ {\it not necessarily equal
to} 0.
The main purpose of the section is to find a manageable expression for
the number of these watermelons with a fixed number of contacts with
the wall. This is accomplished by showing that these are equinumerous
with another set of vicious walker families (see Proposition~\TI), a
result that has been previously obtained by Brak and Essam 
\cite{\BrEsAB, Cor.~1}. Our proof however is bijective, being based
on a modified jeu de taquin, thus solving a problem posed in
\cite{\BrEsAB}. (We remark that since the first version of the present
article was written, Rubey \cite{\RubMAG} has found a
completely different bijection,
which, in fact, works in a much more general setting.) 
This result is then used in Section~6
to find an exact expression
in form of a double sum for the partition function $Z^{(n)}_t(y;\ka)$ 
for watermelons which start on the wall and end at
some arbitrary fixed deviation (see Theorem~\TK), which is new and
the central result of this paper. By
the use of hypergeometric summation and transformation formulas we
are able to rederive, respectively extend, previously obtained
alternative expressions due to Brak and Essam (private communication)
(see Theorems~\TL\ and \TN). 
By the same techniques we are also able to rederive earlier closed form
results by the author, Guttmann and Viennot \cite{\KrGVAA} in the
case that $\ka=1$, and by
Owczarek, Essam and Brak \cite{\OwEBAA} and Ciucu
and the author \cite{\CiKrAB} in the case that $\ka=2$, 
see Corollary~\TD. Section~7 is devoted to the asymptotic analysis of
the partition function $Z^{(n)}_t(y;\ka)$. Our starting point is the
double sum formula for $Z^{(n)}_t(y;\ka)$ provided by Theorem~\TK,
which, when combined with the technique of singularity
analysis, yields full asymptotic expansions for
$Z^{(n)}_t(y;\ka)$
as the length $t$ of the walks tends to infinity (see
Theorem~\TE). Thus we
confirm the predictions by Owczarek, Essam and Brak \cite{\OwEBAA},
making them more precise at the same time. The final Section~8 is
devoted to the study of the mean number of contacts. In the same
special cases as above, namely for $\ka=1$ and $\ka=2$, we are able
to obtain simple closed formulas, see Theorems~\TF\ and \TFa. 
Singularity analysis again allows
us to obtain full asymptotic expansions for the mean number of
contacts, which is the contents of Theorem~\TG.

\subhead 2. The combinatorics of a single walker\endsubhead
In this section we consider a single walker with wall
interaction. What we aim at is the computation of the partition
function $\sum _{P} ^{}\ka^{c(P)}$, where the
sum is over all walks $P$ from $(0,a)$ to $(t,e)$ which never run
below the $x$-axis (the wall), and where $c(P)$ denotes the number 
of contacts of the walk with the wall. Slightly deviating from the
notation introduced in the Introduction, we denote this partition
function by $Z^{(1)}_{t}(a\to e;\ka)$.

As we are going to develop,
this type of question is classical in the literature on nonparametric
statistics. We summarize what is known in the proposition below.
It gives an explicit formula for the partition
function $Z^{(1)}_{t}(a\to e;\ka)$, and, as a simple corollary, an
algebraic expression for the corresponding generating function. This 
is an old result due to Engelberg \cite{\EngbAB, Cor.~3.2} 
(for the case where the
starting point is on the $x$-axis) and Mohanty \cite{\MohaAB,
Cor.~1~(iv)} (for the general case). It has recently
been rediscovered by Brak, Essam and Owczarek
in \cite{\BrEOAC, Eq.~(3.11)}. There is an elegant
combinatorial proof, the origin of which is difficult to track down
(see for example \cite{\KaJGAA, proof of Lemma~1} 
for an occurrence of that idea\footnote{The original 
inspiration for that idea may
be due to Cs\'aki and Vincze \cite{\CsViAA, p.~100}. There it
appears in a slightly modified form, and is used there to solve a slightly
modified problem, namely to find 
the number of walks which {\it cross} the $x$-axis
a given number of times.}). 
Since it is apparently less well-known than it should be, we
reproduce it below.

\proclaim{Proposition \TA}Let $a$ and $e$ be non-negative integers
with $a+ e\equiv t$ mod $2$. Then
the partition function $Z^{(1)}_{t}(a\to e;\ka)$
of a single vicious walker starting at $(0,a)$ and ending at $(t,e)$
and never running below the $x$-axis is given by
$$\multline \binom {t}{(t+a-e)/2}-\binom {t}{(t-a-e)/2}\\
+\sum _{\ell\ge1} ^{}\(\binom {t-\ell}{(t+a+e-2)/2}-\binom
{t-\ell}{(t+a+e)/2}\)\ka^\ell.
\endmultline
\tag\ZY$$
The generating function, the coefficients of which are 
these partition functions, is equal to
$$\multline \underset t\equiv a+e\ (\hbox{\sevenrm mod }2)\to{\sum 
_{t\ge0} ^{}}Z^{(1)}_{t}(a\to e;\ka)z^{t}=
\frac {1} {\sqrt{1-4z^2}}\(\frac {1-\sqrt{1-4z^2}} {2z}\)^{\vert
e-a\vert}\\
-
\frac {1} {\sqrt{1-4z^2}}\(\frac {1-\sqrt{1-4z^2}} {2z}\)^{e+a}
\(1-\frac {\frac {2\ka} {2-\ka}\sqrt{1-4z^2}} {1+\frac {\ka}
{2-\ka}\sqrt{1-4z^2}}\).
\endmultline
\tag\ZZ$$
\endproclaim

\demo{Proof} Let for the moment $a$ and $e$ be at least 1.
The coefficient of $\ka^0$ in $Z^{(1)}_{t}(a\to e;\ka)$ counts the
walks from $(0,a)$ to $(t,e)$ which {\it do not touch} the $x$-axis.
As is classical, this number can be obtained from the reflection
principle (see e.g\. \cite{\ComtAA, p.~22}). It is equal to
$$\binom {t}{(t+a-e)/2}-\binom {t}{(t-a-e)/2},\tag\ZZa$$
which explains the first line in (\ZY). 

On the other hand, for
$\ell\ge1$ the coefficient of $z^\ell$ counts the
walks from $(0,a)$ to $(t,e)$ which do not cross the $x$-axis and touch
it {\it exactly} $\ell$ times. Let us consider such a walk, see
Figure~\FB\ for an example in which $t=17$, $a=3$, $e=2$, $\ell=3$.

\midinsert
\vskip10pt
\vbox{
$$
\Koordinatenachsen(18,5)(0,0)
\Gitter(18,5)(0,0)
\Pfad(0,3),4434\endPfad
\Pfad(4,1),4\endPfad
\raise1pt\hbox to0pt{\Pfad(4,1),4\endPfad\hss}%
\raise-1pt\hbox to0pt{\Pfad(4,1),4\endPfad\hss}%
\Pfad(5,0),3334434\endPfad
\Pfad(12,1),4\endPfad
\raise1pt\hbox to0pt{\Pfad(12,1),4\endPfad\hss}%
\raise-1pt\hbox to0pt{\Pfad(12,1),4\endPfad\hss}%
\Pfad(13,0),3\endPfad
\Pfad(14,1),4\endPfad
\raise1pt\hbox to0pt{\Pfad(14,1),4\endPfad\hss}%
\raise-1pt\hbox to0pt{\Pfad(14,1),4\endPfad\hss}%
\Pfad(15,0),33\endPfad
\PfadDicke{.5pt}
\Pfad(17,-1),222222\endPfad
\DickPunkt(0,3)
\DickPunkt(17,2)
\Label\lo{\hbox{$(0,a)$\hskip.5cm}}(0,3)
\Label\ro{\hbox{\hskip.5cm$(t,e)$}}(17,2)
\hskip8.5cm
$$ 
\centerline{\eightpoint A walk touching the $x$-axis $\ell$ times}
\vskip7pt
\centerline{\eightpoint Figure \FB}
}
\vskip10pt
\endinsert

Such a walk is now transformed into a walk from $(0,a)$ to
$(t-\ell,e+\ell)$ by deleting all the steps {\it immediately
preceding a touching point} from the original walk and gluing the
walk pieces together. In our example in Figure~\FB, these steps are
indicated by thick line segments. Figure~\FC\ shows the result after
deletion of these steps. The circles should be ignored at this point.

\midinsert
\vskip10pt
\vbox{
$$
\Koordinatenachsen(18,6)(0,0)
\Gitter(18,6)(0,0)
\Pfad(0,3),44343334434333\endPfad
\PfadDicke{.5pt}
\Pfad(17,-1),2222222\endPfad
\SPfad(-1,1),1111111111111111111\endSPfad
\DickPunkt(0,3)
\DickPunkt(14,5)
\Kreis(4,1)
\Kreis(11,2)
\Kreis(12,3)
\Label\lo{\hbox{$(0,a)$\hskip.5cm}}(0,3)
\Label\o{(t-\ell,e+\ell)}(14,5)
\Label\r{y=1}(18,1)
\hskip8.5cm
$$ 
\centerline{\eightpoint After deletion of steps}
\vskip7pt
\centerline{\eightpoint Figure \FC}
}
\vskip10pt
\endinsert

We claim that this mapping is indeed a bijection between
walks from $(0,a)$ to $(t,e)$ which do not cross the $x$-axis and touch
it exactly $\ell$ times and walks from $(0,a)$ to
$(t-\ell,e+\ell)$ which {\it do not run below the horizontal line $y=1$} but
{\it touch $y=1$ at least once}. For establishing the claim we have
to explain how we can reconstruct the original walk from one of the
walks of the second type. In order to accomplish this, we
just have to find the points in the walk where steps had been
deleted. Indeed, as is straight-forward to see, for each $j$, $1\le
j\le\ell$, the right-most point on the walk which is on the
horizontal line $y=j$ is such a point, and these are in fact all of
them. In our example in Figure~\FC\ these points are circled.

Hence, the number that we want to find is equal to the number of
walks from $(0,a)$ to $(t-\ell,e+\ell)$ which do not run below the
line $y=1$ {\it minus} the number of
walks from $(0,a)$ to $(t-\ell,e+\ell)$ which do not run below the
line $y=2$. By another application of the reflection principle, this
is
$$\multline
\(\binom {t-\ell}{(t+a-e-2\ell)/2}-\binom {t-\ell}{(t-a-e-2\ell)/2}\)\\
-\(\binom {t-\ell}{(t+a-e-2\ell)/2}-\binom {t-\ell}{(t-a-e-2\ell+2)/2}\),
\endmultline$$
which, after cancellation, gives exactly the expression which appears
as the coefficient of $\ka^\ell$, $\ell\ge1$, in (\ZY). 

The above considerations establish (\ZY) for $a,e\ge1$. If $a=0$ and
$e\ge1$, then we argue that the partition function $Z^{(1)}_{t}(0\to
e;\ka)$ is equal to $\ka Z^{(1)}_{t-1}(1\to e;\ka)$ because the first
step in a walk that starts at $(0,0)$ must necessarily be an up-step,
the additional factor $\ka$ taking into account the touching of the
original walk at $(0,0)$. Now we may apply formula (\ZY) with $a=1$.
It is just a matter of simple manipulations to convert the obtained
expression to (\ZY) with $a=0$. If $a\ge1$ and $e=0$, respectively
if $a=e=0$, one argues similarly. We leave the details to the reader.

\medskip
In order to establish the second assertion of the proposition,
we must compute the generating function
$$\underset t\equiv a+e\ (\hbox{\sevenrm mod }2)\to{\sum 
_{t\ge0} ^{}}\(\binom {t}{(t+a-e)/2}-\binom {t}{(t-a-e)/2}\)z^t
\tag\ZT$$
on the one hand, and
$$\underset t\equiv a+e\ (\hbox{\sevenrm mod }2)\to{\sum 
_{t\ge0} ^{}}\bigg(\sum _{\ell\ge1} ^{}\(\binom {t-\ell}{(t+a+e-2)/2}-\binom
{t-\ell}{(t+a+e)/2}\)\ka^\ell\bigg)z^t
\tag\ZU$$
on the other.

We may concentrate on the case where $e\ge a$, because it is
combinatorially obvious that an interchange of $a$ and $e$ gives the
same result. Therefore, from now on we assume $e\ge a$. 

In order to compute the series (\ZT), we write
it in a telescoping form as
$$\multline 
\sum _{j=0} ^{a-1}\underset t\equiv a+e\ (\hbox{\sevenrm mod }2)\to{\sum 
_{t\ge0} ^{}}\(\binom {t}{(t+a-e-2j)/2}-\binom
{t}{(t+a-e-2j-2)/2}\)z^t\\
=\sum _{j=0} ^{a-1}{\sum 
_{r\ge0} ^{}}\(\binom {2r+2j+e-a}{r}-\binom
{2r+2j+e-a}{r-1}\)z^{2r+2j+e-a},
\endmultline
\tag\ZV
$$
where in the next-to-last line we performed the index transformation $t\to
2r+2j+e-a$.

Now we appeal to the well-known fact (see e.g\. \cite{\MohaAE, displayed 
line before (1.21) with $\mu=1$}, in combination with \cite{\MohaAE, (1.19)}) 
that 
$$\sum _{r=0} ^{\infty}\frac {m+1} {r+m+1}\binom
{2r+m}{r}x^r=\sum _{r=0} ^{\infty}\(\binom
{2r+m}{r}-\binom {2r+m}{r-1}\)x^r=C(x)^{m+1},\tag\ZS$$
where $C(x)$ is the generating function for the Catalan numbers,
$$C(x)=\sum _{r=0} ^{\infty}\frac {1} {r+1}\binom
{2r}{r}x^r=\sum _{r=0} ^{\infty}\(\binom
{2r}{r}-\binom {2r}{r-1}\)x^r=\frac {1-\sqrt{1-4x}} {2x}.\tag\ZI$$
If we use this in (\ZV), then we find that the series (\ZT) is
equal to
$$\multline
\sum _{j=0} ^{a-1}z^{2j+e-a}C(z^2)^{2j+e-a+1}
=\frac {z^{e-a}C(z^2)^{e-a+1}-
z^{e+a}C(z^2)^{e+a+1}} {1-z^2C(z^2)^2}\\
=\frac {1} {\sqrt{1-4z^2}}\(\(\frac {1-\sqrt{1-4z^2}} {2z}\)^{
e-a}-\(\frac {1-\sqrt{1-4z^2}} {2z}\)^{e+a}\),
\endmultline
\tag\ZW
$$
the last line being due to the 
easily verified fact that 
$$1-xC(x)^2=C(x)\sqrt{1-4x}\tag\ZWa$$.

For the computation of the series (\ZU) we proceed in a similar
fashion. Here, we first interchange summations and then replace the
index $t$ of the (now) inner summation by $2r+2\ell+a+e-2$ to
transform the series to
$$\sum _{\ell\ge1} ^{}\sum _{r\ge0} ^{}\(\binom 
{2r+\ell+a+e-2}{r}-\binom
{2r+\ell+a+e-2}{r-1}\)z^{2r+2\ell+a+e-2}\ka^\ell.
$$
Now we apply (\ZS) to the inner summation and obtain
$$\multline
\sum _{\ell\ge1} ^{}z^{2\ell+a+e-2}C(z^2)^{\ell+a+e-1}\ka^\ell
=\frac {\ka\,z^{a+e}C(z^2)^{a+e}}{1-\ka z^2C(z^2)}\\
=\frac {1} {\sqrt{1-4z^2}}\(\frac {1-\sqrt{1-4z^2}} {2z}\)^{e+a}
\frac {\frac {2\ka} {2-\ka}\sqrt{1-4z^2}} {1+\frac {\ka}
{2-\ka}\sqrt{1-4z^2}}.
\endmultline
\tag\ZX$$
Summing the expressions (\ZW) and (\ZX) gives exactly (\ZZ).\quad \quad
\qed
\enddemo

\subhead 3. The combinatorics of $n$ vicious walkers\endsubhead
In this section we recall Lindstr\"om's classical result \cite{\LindAA,
Lemma~1} on the
enumeration of families of vicious walkers (non-intersecting 
paths)\footnote{Lindstr\"om used
the term ``pairwise node disjoint paths". The term ``non-intersecting,"
which is most often used nowadays in combinatorial literature,
was coined by Gessel and Viennot \cite{\GeViAA}.} 
in directed graphs (directed networks), and apply it to
obtain determinantal formulae for our partition function
$Z^{(n)}_{t}(\bold a\to \bold e;\ka)$.

Let $G=(V,E)$ be a directed acyclic graph with vertices $V$ and directed
edges $E$. 
Furthermore, we are given a function $w$ which assigns a
weight $w(x)$ to every vertex or edge $x$. Let us define the weight
$w(P)$ of a walk $P$ in the graph by $\prod _{e} ^{}w(e)\prod _{v}
^{}w(v)$, where the first product is over all edges $e$ of the walk
$P$ and the second product is over all vertices $v$ of $P$. 
We denote the set of all walks in $G$ from $u$ to $v$ by $\P(u\to
v)$, and the set of all families $(P_1,P_2,\dots,P_n)$ of walks,
where $P_i$ runs from $u_i$ to $v_i$, $i=1,2,\dots,n$, by $\P(\bold
u\to\bold v)$, with $\bold u=(u_1,u_2,\dots,u_n)$ and 
$\bold v=(v_1,v_2,\dots,v_n)$. The symbol $\P^+(\bold
u\to\bold v)$ stands for the set of all families $(P_1,P_2,\dots,P_n)$ 
in $\P(\bold u\to\bold v)$ with the additional property that no two
walks share a vertex. We call such families of walk(er)s {\it ``vicious
walkers"} or, alternatively, {\it ``non-intersecting lattice paths"}. 
The weight $w(\bold P)$ of a family $\bold
P=(P_1,P_2,\dots,P_n)$ of walks is defined as the product $\prod
_{i=1} ^{n}w(P_i)$ of all the weights of the walks in the family.
Finally, given a set $M$ with weight function $w$, we write 
$\GF(\Cal M;w)$ for the generating function $\sum
_{x\in\Cal M} ^{}w(x)$. 

We need two further notations before we are able to state the
Lindstr\"om--Gessel--Viennot theorem.\footnote{By a curious coincidence,
Lindstr\"om's result (the motivation of which was matroid theory!) 
was rediscovered in the 1980s at about the same
time in three different
communities, not knowing from each other at that time: in statistical
physics by Fisher \cite{\FishAA, Sec.~5.3} in order to apply it to
the analysis of vicious walkers as a model of wetting and melting, 
in combinatorial chemistry by John and Sachs \cite{\JoSaAB} and
Gronau, Just, Schade, Scheffler and Wojciechowski \cite{\GrJSAA}
in order to compute Pauling's bond order
in benzenoid hydrocarbon molecules, and in enumerative combinatorics
by Gessel and Viennot \cite{\GeViAA, \GeViAB} in order to count
tableaux and plane partitions. Since only Gessel and Viennot
rediscovered it in its most general form, I propose to call this theorem
the ``Lindstr\"om--Gessel--Viennot theorem." It must however be
mentioned that in fact the same idea appeared even earlier in work by
Karlin and McGregor \cite{\KaMGAB, \KaMGAC} in a probabilistic
framework, as well as that the so-called ``Slater determinant"
in quantum mechanics (cf\. \cite{\SlatZY} and \cite{\SlatZZ, Ch.~11}) 
may qualify as an ``ancestor" of
the Lindstr\"om--Gessel--Viennot determinant.} The symbol $S_n$
denotes the symmetric group of order $n$. Given a permutation $\si\in
S_n$, we write $\bold u_\si$ for
$(u_{\si(1)},u_{\si(2)},\dots,u_{\si(n)})$. Then
$$\sum _{\si\in S_n} ^{}(\sgn\si) \cdot\GF(\P^+(\bold u_\si\to \bold
v);w)
=\det_{1\le i,j\le n}\big(\GF(\P(u_{j}\to v_i);w)\big).
\tag\Za$$

Most often, this theorem is applied in the case where the only permutation
$\si$ for which vicious walks exist is the identity
permutation, so that the sum on the left-hand side reduces to a
single term which counts all families $(P_1,P_2,\dots,P_n)$ of 
vicious walks, the $i$-th walk $P_i$ running from $A_i$
to $E_i$, $i=1,2,\dots,n$. This case occurs for example if for any
pair of walks $(P,Q)$ with $P$ running from $u_a$ to $v_d$ and $Q$
running from $u_b$ to $v_c$, $a<b$ and $c<d$, it is true that $P$ and
$Q$ must have a common vertex. Explicitly, in that case we have
$$
\GF(\P^+(\bold u\to \bold v);w)=
\det_{1\le i,j\le n}\big(\GF(\P(u_{j}\to v_i);w)\big).
\tag\Zb$$
This is also the case which we encounter in this
paper.\footnote{There exist however also several 
interesting applications of the
{\it general\/} form of the Lindstr\"om--Gessel--Viennot theorem in the
literature, see \cite{\CiEKAA, \FiscAA, \SulaAC}.}

\proclaim{Proposition \TH}Let $\bold a=(a_1,a_2,\dots,a_n)$ and
$\bold e=(e_1,e_2,\dots,e_n)$ be $n$-tuples of non-negative integers
with $0\le a_1<a_2<\dots<a_n$, $0\le e_1<e_2<\dots<e_n$, all $a_i$'s of the
same parity, all $e_i$'s of the same parity, such that $a_i+e_i\equiv
t$ {\rm(mod 2)} for all $i$.
As before, let $Z^{(n)}_{t}(\bold a\to\bold e;\ka)$ denote the partition
function for families of 
$n$ vicious walkers, the $i$-th starting at $(0,a_i)$ and
ending at $(t,e_i)$, none of them running below the $x$-axis, where the
weight of a walker configuration $\bold P$ is defined as 
$\ka^{c(\bold P)}$ with $c(\bold P)$
denoting the number of contacts of the walkers with the $x$-axis. Then
$$Z^{(n)}_{t}(\bold a\to\bold e;\ka)=
\det_{1\le i,j\le n}\(Z^{(1)}_t(a_j\to e_i;\ka)\),\tag\Zc$$
where $Z^{(1)}_t(a\to e;\ka)$ is given by {\rm(\ZY)}.
\endproclaim
\demo{Proof}We apply (\Zb) for $G$ the graph with vertices the points
$(x,y)$ in the plane integer lattice with $y\ge0$ and with directed
edges $(x,y)\to(x+1,y-1)$, $y\ge1$, and $(x,y)\to(x+1,y+1)$, $y\ge0$.
As the weight $w$ we choose the function which assigns 1 to every
edge and to every vertex $(x,y)$ with $y>0$, and which assigns $\ka$
to every vertex on the $x$-axis. If we now choose $u_i=(0,a_i)$ and
$v_i=(t,e_i)$ in the
Lindstr\"om--Gessel--Viennot theorem (\Zb), then the generating
function on the left-hand side is exactly the partition function 
$Z^{(n)}_{t}(\bold a\to\bold e;\ka)$ for $n$ vicious walkers, 
whereas the generating function
which gives the $(i,j)$-entry of the determinant on the right-hand
side is the partition function $Z^{(1)}_{t}(a_j\to e_i;\ka)$ for a
single walker from $(0,a_j)$ to $(t,e_i)$, the latter being given
explicitly by Proposition~\TA.\quad \quad \qed
\enddemo

This result has also been obtained by Brak, Essam and Owczarek 
in \cite{\BrEOAD, Eqs.~(22), (23)}
(using the Bethe Ansatz method; their derivation is much more complicated,
but it is on the other hand a method which is more
widely applicable that they use). 

\medskip
In this paper we will primarily analyse the case of watermelon configurations
which start on the wall. To be precise, this is the case where the
starting point $u_i$ is $(0,2i-2)$ and the end point $v_i$ is
$(t,y+2i-2)$, $i=1,2,\dots,n$, for some non-negative integer $y$.
In our analysis of this particular case we shall in fact use two
variants of the above formula, which are also corollaries of the
Lindstr\"om--Gessel--Viennot theorem.

\proclaim{Proposition \TB}
Let $t$ and $y$ be non-negative integers with $t\equiv y$ {\rm(mod 2)}.
As in the Introduction, let $Z^{(n)}_{t}(y;\ka)$ 
denote the partition function for families of 
$n$ vicious walkers, the $i$-th starting at $(0,2i-2)$ and
ending at $(t,y+2i-2)$, $i=1,2,\dots,n$, 
none of them running below the $x$-axis, where the
weight of a walker configuration $\bold P$ is defined as 
$\ka^{c(\bold P)}$ with $c(\bold P)$
denoting the number of contacts of the walkers with the $x$-axis. Then
$$Z^{(n)}_{t}(y;\ka)=\frac {1} {\ka^{n-1}}
\det_{0\le i,j\le n-1}(B_{i,j}(t,y;\ka)),\tag\ZB$$
where
$$\align
B_{i,j}(t,y;\ka)&=\sum _{\ell\ge1} ^{}\(\binom {t+2j-\ell}{\frac
{t+y}2+i+j-1}-\binom
{t+2j-\ell}{\frac {t+y}2+i+j}\)\ka^\ell\\
&=
\sum _{\ell\ge1} ^{}\frac {y+2i+\ell-1} {\frac {t+y} {2}+i+j}\binom
{t+2j-\ell}{\frac {t+y} {2}+i+j-1}\ka^\ell.
\tag\ZA\endalign$$
Furthermore, let $Z^{(n)}_{2r}(\ka)$ denote the partition
function for families of 
$n$ vicious walkers, the $i$-th starting at $(0,2i-2)$ and
ending at $(2r,2i-2)$, $i=1,2,\dots,n$, 
none of them running below the $x$-axis, with the same
weight of walker configurations. Then
$$Z^{(n)}_{2r}(\ka)=\frac {1} {\ka^{2n-2}}
\det_{0\le i,j\le n-1}(C(r+i+j;\ka)),\tag\ZBa$$
where
$$C(r;\ka)=\sum _{\ell=2} ^{r+1}\(\binom {2r-\ell}{r-1}-\binom
{2r-\ell}r\)\ka^\ell=
\sum _{\ell=2} ^{r+1}\frac {\ell-1} {r}\binom
{2r-\ell}{r-1}\ka^\ell.
\tag\ZAa$$
\endproclaim

\demo{Proof} We start by proving the first claim. 
Given a family of vicious walkers as in the first statement of
the proposition, we may freely attach $2i-2$ up-steps at the
beginning of the $i$-th walk, $i=1,2,\dots,n$, 
see Figure~\FD\ for the resulting walks if we do this with the
watermelon configuration in Figure~\FA.b.

\midinsert
\vskip10pt
\vbox{
$$
\Einheit.4cm
        \Gitter(13,14)(-8,-1)
        \Koordinatenachsen(13,14)(-8,-1)
\Pfad(0,0),343343443433\endPfad
\Pfad(-2,0),33334333444433\endPfad
\Pfad(-4,0),3333334333443434\endPfad
\Pfad(-6,0),333333333343344344\endPfad
\Pfad(-8,0),33333333333343343444\endPfad
\DickPunkt(0,0)
\DickPunkt(0,2)
\DickPunkt(0,4)
\DickPunkt(0,6)
\DickPunkt(0,8)
\DickPunkt(12,2)
\DickPunkt(12,4)
\DickPunkt(12,6)
\DickPunkt(12,8)
\DickPunkt(12,10)
        \Label\lu{u'_1}(0,0)
        \Label\lu{u'_2}(-2,0)
        \Label\lu{u'_3}(-4,0)
        \Label\lu{u'_4}(-6,0)
        \Label\lu{u'_5}(-8,0)
        \Label\r{v_1}(12,2)
        \Label\r{v_2}(12,4)
        \Label\r{v_3}(12,6)
        \Label\r{v_4}(12,8)
        \Label\r{v_5}(12,10)
\hbox{\hskip1.5cm}
$$
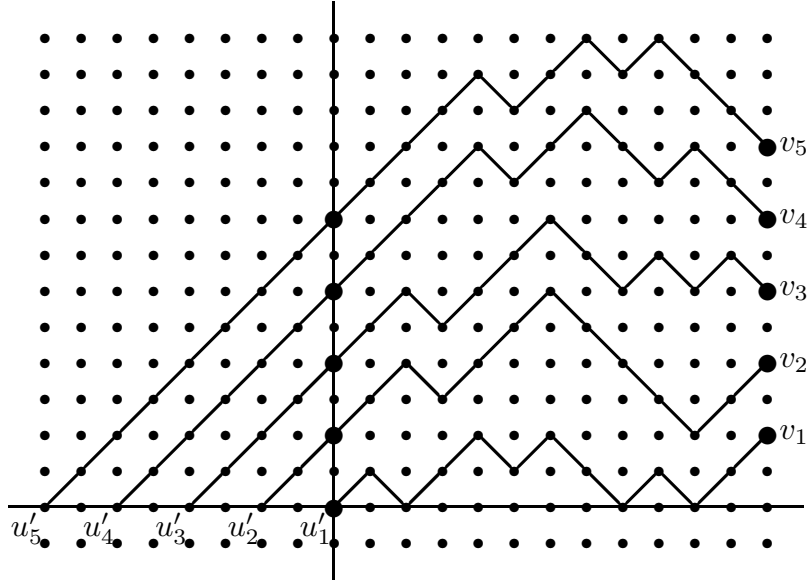
\centerline{\eightpoint Watermelon configuration with attached walk
pieces at the beginning}
\vskip7pt
\centerline{\eightpoint Figure \FD}
}
\vskip10pt
\endinsert

It is obvious that the number of families of vicious walkers 
with starting points
$u'_i=(-2i+2,0)$ (instead of $u_i=(0,2i-2)$) and end points 
$v_i=(t,y+2i-2)$, $i=1,2,\dots,n$,
each walker not running below the $x$-axis (see
Figure~\FD),
is exactly the same as the number of families of vicious walkers
in the first statement of the proposition (compare Figure~\FA). Clearly, by
attaching these walk pieces at the beginning, we
introduced $n-1$ further contacts of the first walk with the $x$-axis
(namely in the points $u'_n,u'_{n-1},\dots, u'_2$).

If we now apply (\Zb) with $u_i$ replaced by $u'_i$, 
then we obtain that the partition function for the
latter families of walkers is given by
$$\det_{1\le i,j\le n}\big(\GF(\P((-2j+2,0)\to (t,y+2i-2));w)\big),$$
where $w$ is our contact weight. By definition, the generating
function 
$$\GF(\P((-2j+2,0)\to (t,y+2i-2));w)$$ 
is equal to
$Z^{(1)}_{t+2j-2}(0\to y+2i-2;\ka)$, which by Proposition~\TA\ is exactly the
same as $B_{i-1,j-1}(t,y;\ka)$. Since our operation of adding walk pieces
introduced $n-1$ additional contacts with the $x$-axis, we must divide
the above determinant by $\ka^{n-1}$, and we obtain the final result
(\ZB) after replacing $i$ by $i+1$ and $j$ by $j+1$.

The proof of formula (\ZBa) is completely analogous. Here, one
has to also attach $2i-2$ {\it down}-steps at the {\it end\/} of the $i$-th
walk. We leave the details to the reader.\quad \quad \qed
\enddemo

\subhead 4. Exact results for the partition function
for watermelons of deviation 0\endsubhead
In this section we address
the exact evaluation of the partition function
of watermelon configurations which start {\it and end\/} on the wall. To be
precise, we consider the partition function 
$Z^{(n)}_{2r}(\ka)=\sum _{\bold P} ^{}\ka^{c(\bold
P)}$, where the sum is over all families $\bold P=(P_1,P_2,\dots,P_n)$
of $n$ vicious walkers never running below the $x$-axis, 
the $i$-th walker $P_i$ starting at
$(0,2i-2)$ and ending at $(2r,2i-2)$, $i=1,2,\dots,n$, and with $c(\bold
P)$ denoting the total number of contacts of the walkers with 
the wall. It should be noted that
$Z^{(n)}_{2r}(\ka)=Z^{(n)}_{2r}((0,2,4,\dots,2n-2)\to(0,2,4,\dots,2n-2);\ka)$
in the notation of the Introduction and of Proposition~\TH.

The theorem below provides an expression for the partition function
$Z^{(n)}_{2r}(\ka)$ in terms of a single sum, by finding a product
formula for the number of the above watermelon configurations {\it
with exactly $\ell$ contacts} with the wall. This result was first
obtained by Brak and Essam \cite{\BrEsAB, Theorem~6} by clever determinant
manipulations, which allowed them to find a recurrence for the
partition function. Our proof is completely different. It
demonstrates that, once we have a guess for the formula (which can be
found completely automatically using the {\sl Mathematica} program
{\tt Rate}\footnote{written by the author;
available from {\tt
http://www.mat.univie.ac.at/\~{}kratt}.}, 
respectively its {\sl Maple} equivalent {\tt GUESS}\footnote{written by 
Fran\c cois B\'eraud and Bruno Gauthier; available from\break 
{\tt http://www-igm.univ-mlv.fr/\~{}gauthier}.}), it is already
proved, as the validity of the formula 
can be established in a completely automatic fashion, by making
use of what is known under the name of ``Dodgson's
condensation me\-thod" and the
Gosper--Zeilberger algorithm. A more
general result will be obtained later in Theorem~\TK, however using
completely different (non-automatic) methods.

\proclaim{Theorem \TC} The partition function $Z^{(n)}_{2r}(\ka)$ is
equal to 
$$\frac {(r-1)!\prod _{i=0} ^{n-1}(2i+1)!\prod _{i=0} ^{n-2}(2r+2i)!} 
{\prod _{i=0} ^{2n-2}(r+i)!}\sum _{\ell=0} ^{r-1}\binom
{2r-\ell-2}{r-1}\binom {\ell+2n-1}\ell\ka^{\ell+2}.\tag\ZG$$
In particular, the number of walker configurations with
exactly $\ell$ contacts with the $x$-axis is equal to
$$\frac {(r-1)!\prod _{i=0} ^{n-1}(2i+1)!\prod _{i=0} ^{n-2}(2r+2i)!} 
{\prod _{i=0} ^{2n-2}(r+i)!}\binom
{2r-\ell}{r-1}\binom {\ell+2n-3}{\ell-2}.$$
\endproclaim

\demo{Proof} ``Dodgson's
condensation me\-thod" (cf\. \cite{\KratBN, Sec.~2.3}) is based on
the following determinant identity due to Desnanot and Jacobi.
Let $A$ be an $n\times n$ matrix. Denote the submatrix of $A$ in which
rows $i_1,i_2,\dots,i_k$ and columns $j_1,j_2,\dots,j_k$ are 
omitted by $A_{i_1,i_2,\dots,i_k}^{j_1,j_2,\dots,j_k}$. Then there holds
$$\det A\cdot \det A_{1,n}^{1,n}=\det A_{1}^{1}\cdot \det A_n^n-
\det A_1^n\cdot \det A_n^1.
\tag\ZD$$
If $n=2$, then $\det A^{1,n}_{1,n}$, the determinant of a
$0\times 0$ matrix, has to be interpreted as 1.

Because of Proposition~\TB, $Z^{(n)}_{2r}(\ka)$ is essentially equal
to a determinant. We now apply the above determinant identity with $A$
equal to the determinant in (\ZBa). Thus we obtain
$$Z^{(n)}_{2r}(\ka)Z^{(n-2)}_{2r+4}(\ka)=
Z^{(n-1)}_{2r+4}(\ka)Z^{(n-1)}_{2r}(\ka)-\(Z^{(n-1)}_{2r+2}(\ka)\)^2.
\tag\ZC$$
This relation enables us to prove the theorem by induction on
$n$. It is clearly true for $n=1$, thanks to Proposition~\TA.
Equation~(\ZC) allows to perform the induction step. We only have to
check that (\ZC) is true with $Z^{(n)}_{2r}(\ka)$ as asserted in the
statement of the theorem. 

Let first $n=2$. Then the term $Z^{(0)}_{2r+4}(\ka)$ appears in
(\ZC), which by virtue of (\ZBa) is $\ka^2$ times the determinant of a
$0\times0$ matrix. Since, according to our earlier convention, the
latter determinant should be interpreted as 1, we
should let $Z^{(0)}_{2r}(\ka):=\ka^2$. Thus, for $n=2$ Equation~(\ZC) becomes
$$\multline \frac {(2r)!} {r\,(r+1)!\,(r+2)!}\sum _{\ell=0} ^{r-1}
\binom {2r-\ell-2}{r-1}(\ell+1)(\ell+2)(\ell+3)\ka^{\ell+4}\\
=\(\sum _{\ell=0} ^{r+1}\frac {\ell+1} {r+2}\binom
{2r-\ell+2}{r+1}\ka^{\ell+2}\)
\(\sum _{\ell=0} ^{r-1}\frac {\ell+1} {r}\binom
{2r-\ell-2}{r-1}\ka^{\ell+2}\)
\\-
\(\sum _{\ell=0} ^{r}\frac {\ell+1} {r+1}\binom
{2r-\ell}{r}\ka^{\ell+2}\)^2.
\endmultline$$
By comparison of coefficients of $\ka^{e+4}$ on both sides, this is seen
to be equivalent to verifying that
$$\multline 
(e+1)(e+2)(e+3)\binom {2r-e-2}{r-1}\\
=\sum _{\ell=0} ^{e}\Bigg(\frac {\ell+1} {r+2}\binom
{2r-\ell+2}{r+1}
\frac {e-\ell+1} {r}\binom
{2r-e+\ell-2}{r-1}\\
-\frac {\ell+1} {r+1}\binom
{2r-\ell}{r}\frac {e-\ell+1} {r+1}\binom
{2r-e+\ell}{r}\Bigg).
\endmultline$$
This identity is easily proved once one observes that the summand of
the sum on the right-hand side is equal to
$G(e,\ell+1)-G(e,\ell)$, where
$$\multline G(e,\ell)={\frac{  
     \left( 2 r - \ell +1 \right) !\, \left(  2 r- e + \ell -2 \right) !} 
   {\left( r+1 \right) ! \,\left( r+2 \right) ! \,
     \left( r - \ell + 1 \right) !\, \left( r - e + \ell -1 \right) !}}
\\
\times
( -6 - 8 e - 2 {e^2} + 9 \ell + 6 e \ell + {e^2} \ell - {\ell^2} + 
       2 e {\ell^2} + {e^2} {\ell^2} - 3 {\ell^3} - 2 e {\ell^3} \\
+ {\ell^4} -        6 e r 
- 2 {e^2} r - 3 \ell r - 4 e \ell r - {e^2} \ell r + 
       3 {\ell^2} r + e {\ell^2} r + 6 {r^2} + 2 e {r^2} ),
\endmultline$$
since then the sum on the right-hand side is equal to
$G(e,e+1)-G(e,0)$, which can be seen to be equal to the left-hand
side upon some simplification. (Clearly, Gosper's algorithm
\cite{\GospAB}, \cite{\GrKPAA, \S5.7}, \cite{\PeWZAA, \S II.5} was used
to find $G(e,\ell)$. The particular implementation that we used is
the {\sl Mathematica} implementation by Paule and Schorn
\cite{\PaScAA}.)

Now let $n>2$. In this case, substitution of the claimed expression
for $Z^{(n)}_{2r}(\ka)$ in (\ZC) and comparison of coefficients of 
$\ka^{e+4}$ on both sides yields, after some manipulation, 
that we have to verify the summation
$$\multline \sum _{\ell=0} ^{e}
\Bigg(  {{\left( n-1 \right)  \left( 2 n-1 \right)  \left( r+1 \right)  }
     \over {\left( 2 n + r-2 \right)  \left( 2 r+1 \right)}}
      {\binom {2r - \ell - 2} {r-1}} {\binom {2n + \ell -1} \ell}\\
\cdot
{\binom {2r - e + \ell + 2} {
r+1}}      {\binom { 2 n + e - \ell -5} {e - \ell}} 
\\
-{{\left( r+1 \right)  \left( n + r-1 \right)  
       \left(  2 n + 2 r-3 \right)  }
      \over {\left( 2 n + r-2 \right)  \left(  2 r+1 \right) }} 
{\binom {2r - \ell + 2} {r+1}} {\binom {2n + \ell -3} \ell} 
\\
\cdot
        {\binom {2r - e + \ell - 2} {
r-1}}       {\binom {2n + e - \ell -3} {e - \ell}}
\\+ 
  r {\binom { 2 r-\ell} r} {\binom {2n + \ell -3} \ell} 
    {\binom {2r-e + \ell } r} 
{\binom {2n + e - \ell -3} {e - \ell}}\Bigg)
=0.
\endmultline\tag\ZE$$
The most straight-forward way to do this is to feed the sum into the
Gosper--Zeilberger algorithm \cite{\PeWZAA, \ZeilAM, \ZeilAV} 
(again we used the {\sl Mathematica}
implementation by Paule and Schorn \cite{\PaScAA}). If $S(n)$ denotes
the sum on the left-hand side of (\ZE), then the output of the
algorithm is the recurrence
$$\multline 
   \left(  n-1 \right)  \left( 2 n-1 \right)  
    \left(  e + 4 n-5 \right)  \left( 2 n + r-1 \right)  
    \left( 2 n + r \right)  {S}(n+1) \\
-\left(   e + 4 n-1 \right)  \left(  e + 2 n - r-2 \right)  
    \left(  e + 2 n - r-1 \right)  \left( n + r \right)  
    \left(  2 n + 2 r-1 \right)  { S}(n) = 0.
\endmultline\tag\ZF$$
Now, it is straight-forward to check that $S(1)=0$, thus proving (\ZE)
for $n=1$. The validity of (\ZE) for arbitrary $n$ then follows upon
induction on $n$ from the recurrence (\ZF).

Thus, the proof of the theorem is complete.\quad \quad \qed
\enddemo

If we reverse the order of summation in the sum in
(\ZG) (i.e., if we replace $\ell$ by $r-1-\ell$) and then rewrite 
it using the standard hypergeometric notation
$${}_p F_q\!\left[\matrix a_1,\dots,a_p\\ b_1,\dots,b_q\endmatrix; 
z\right]=\sum _{m=0} ^{\infty}\frac {\po{a_1}{m}\cdots\po{a_p}{m}}
{m!\,\po{b_1}{m}\cdots\po{b_q}{m}} z^m\ ,$$
where the Pochhammer symbol 
$(\alpha)_m$ is defined by $(\alpha)_m:=\alpha(\alpha+1)\cdots(\alpha+m-1)$,
$m\ge1$, $(\alpha)_0:=1$, then we obtain
$$
{ \ka^{1 + r} 
  {{({ \textstyle 2 n}) _{r-1}} \over  { (r-1)!}}} 
  {} _{2} F _{1} \!\left [ \matrix {r, 1 - r}\\ { 2 - 2 n - r}\endmatrix 
; {\displaystyle \frac{1}{\ka}}\right ]  ,
$$
or, more precisely 
(the above $_2F_1$-series is actually ill-defined because of
the denominator parameter $2-2n-r$ which is a negative integer),
$$\lim_{\ep\to0}
{ \ka^{1 + r} 
  {{({ \textstyle 2 n}) _{r-1}} \over  { (r-1)!}}} 
  {} _{2} F _{1} \!\left [ \matrix {r, 1 - r}\\ 
{ 2 - 2 n - r+\ep}\endmatrix 
; {\displaystyle \frac{1}{\ka}}\right ]  .
\tag\ZN$$
In view of the fact that there are numerous $_2F_1$-transformation
formulae, we can obtain
numerous equivalent expressions for this sum and, hence, for (\ZG)
itself. 

A particular such expression has been given by
Owczarek, Essam and Brak in \cite{\OwEBAA, Eq.~(4.65)},
with proof provided by Brak and Essam in \cite{\BrEsAC}, and was the
starting point for their scaling analysis of the partition
function $Z^{(n)}_{2r}(\ka)$ in \cite{\OwEBAA}. This expression provides
an expansion of (\ZG) as an (infinite) power series in $(\ka-1)/\ka^2$. It is
however only valid for $\ka\le2$. In the corollary below, we show 
how the expression of Owczarek, Essam and Brak can be readily derived 
by applying appropriate $_2F_1$-transformation formulas to (\ZN).
Moreover, this technique yields also an expression which is valid
for $\ka>2$. In the case $n=1$ the latter expansion was, 
for example, previously observed by Brak and Essam \cite{\BrEsAA, Eq.~(25)}.

\proclaim{Corollary \TCa}
Let $\ka>2(\sqrt2-1)$.
Then the partition function $Z^{(n)}_{2r}(\ka)$ is
equal to
$$\multline 
\frac {\prod _{i=0} ^{n-1}(2i+1)!\,(2r+2i)!} 
{\prod _{i=0} ^{2n-1}(r+i)!}\Bigg(\frac {1} {\ka^{2n-2}}
\sum _{h=0} ^{\infty}\binom {n+h-1}h \frac
{(2r+2n-1)_{2h}} {(r+n)_h\,(r+2n)_h}\(\frac {\ka-1}
{\ka^2}\)^h\\
+\chi(\ka>2)
\frac {(\ka-2)\ka^{2r+2n-1}} {(\ka-1)^{r+2n-1}}\frac {(r+2n-1)!} 
{(2n-1)!\,(r+n-1)!}\kern3cm\\
\cdot
\sum _{h=0} ^{n-1}\binom {n-1}h\frac
{(r+n-h-1)!\,(r+2n-h-2)!} {(2r+2n-2h-2)!}\(\frac {1-\ka}
{\ka^2}\)^h  
\Bigg),
\endmultline\tag\ZNa$$
where $\chi(A)=1$ if $A$ is true and $\chi(A)=0$ otherwise. 
\endproclaim

\demo{Proof}
First let $\ka\le2$.
We write the sum in (\ZG) in the form (\ZN) and
apply the transformation formula (cf\. \cite{\SlatAC, (1.8.10),
terminating form})
$$
{} _{2} F _{1} \!\left [ \matrix { a, -N}\\ { c}\endmatrix ; {\displaystyle
   z}\right ]  =
{{    ({ \textstyle c-a}) _{N} }\over {({ \textstyle c}) _{N} }}
{} _{2} F _{1} \!\left [ \matrix { -N, a}\\ { 1 + a - c -
       N}\endmatrix ; {\displaystyle 1 - z}\right ]  ,
$$
where $N$ is a non-negative integer. Now the limit 
$\ep\to0$ can be safely performed, and we obtain
$$
\ka^{1 + r}\frac{ 
    ({ \textstyle 2 n + r}) _{r-1} }{({ \textstyle 1}) _{r-1} }
{} _{2} F _{1} \!\left [ \matrix { 1 - r, r}\\ { 2 n + 
r}\endmatrix ; {\displaystyle 1 - \frac{1}{\ka}}\right ]  
$$
for the sum in (\ZG). Next we apply the quadratic transformation formula
(see \cite{\RaVeAA, (3.31), reversed}) 
$$
{} _{2} F _{1} \!\left [ \matrix { a, 1 - a}\\ { c}\endmatrix ; {\displaystyle
   z}\right ]  = {{\left( 1 - z \right) }^{c-1}} 
   {} _{2} F _{1} \!\left [ \matrix { {c\over 2}-{{a}\over 2},{a\over
   2}+{c\over 2}- {{1}\over 2}}\\ 
  { c}\endmatrix ; {\displaystyle 4 z\left( 1 - z \right)  }\right ],
\tag\Zd$$
which holds when $\Re z<1/2$, as well as if $z=1/2$,
and substitute the result back
in (\ZG). After some manipulation, 
we obtain the expression in the first line of (\ZNa),
which, taking into account the condition under which (\Zd) is true
and the latter series converges,
holds for $2(\sqrt 2-1)<\ka\le2$.

If on the other hand $\ka>2$, then we proceed differently.
We start with the expression
$$\multline
\lim_{\ep\to0}
{ \ka^{1 + r} 
  {{({ \textstyle 2 n}) _{r-1}} \over  { (r-1)!}}} 
  {} _{2} F _{1} \!\left [ \matrix {r-\ep, 1 - r+\ep}\\ 
{ 2 - 2 n - r+\ep}\endmatrix 
; {\displaystyle \frac{1}{\ka}}\right ] \\
=
\lim_{\ep\to0}
{ \ka^{1 + r} 
  {{({ \textstyle 2 n}) _{r-1}} \over  { (r-1)!}}} 
\sum _{k=0} ^{\infty}\frac  {(r-\ep)_k\,( 1 - r+\ep)_k} 
{( 2 - 2 n - r+\ep)_k\,k!}\ka^{-k},
\endmultline\tag\Zdd$$
which is clearly equivalent to (\ZN).
{The reader should however note the differences to (\ZN) in placing the $\ep$.}
We split the sum over $k$ into the three ranges $0\le k\le r-1$,
$r\le k\le r+2n-2$, and $r+2n-1\le k$, and then we perform the
limit $\ep\to0$. Thereby each summand in the range $r\le k\le r+2n-2$
vanishes identically. Thus, upon writing the sum over the range $r+2n-1\le
k$ in hypergeometric notation, the expression (\Zdd) becomes
$$
{ \ka^{1 + r} 
  {{({ \textstyle 2 n}) _{r-1}} \over  { (r-1)!}}} 
\sum _{k=0} ^{r-1}\frac  {(r)_k\,( 1 - r)_k} 
{( 2 - 2 n - r)_k\,k!}{\ka^{-k}}
-
\frac{      ({ \textstyle r}) _{r+2 n-1} }
{\ka^{2 n-2}\,
(r + 2 n -1)! }
  {} _{2} F _{1} \!\left [ \matrix { 2 n, 2r + 2 n-1}\\ { 2 n + 
r}\endmatrix ; {\displaystyle \frac{1}{\ka}}\right ] .
\tag\Zddd$$
The first sum in this expression is exactly the sum in (\ZG) (which is
most easily seen from the equivalent expression (\ZN)). Therefore, if
we equate (\Zdd) and (\Zddd), we see that the sum in (\ZG) is equal to
$$\multline
\frac{      ({ \textstyle r}) _{r+2 n-1} }
{\ka^{2 n-2}\,
(r + 2 n -1)! }
  {} _{2} F _{1} \!\left [ \matrix { 2 n, 2r + 2 n-1}\\ { 2 n + 
r}\endmatrix ; {\displaystyle \frac{1}{\ka}}\right ] \\
+
\lim_{\ep\to0}
{ \ka^{1 + r} 
  {{({ \textstyle 2 n}) _{r-1}} \over  { (r-1)!}}} 
  {} _{2} F _{1} \!\left [ \matrix {r-\ep, 1 - r+\ep}\\ 
{ 2 - 2 n - r+\ep}\endmatrix 
; {\displaystyle \frac{1}{\ka}}\right ] .
\endmultline
$$
To the first $_2F_1$-series we apply the quadratic transformation
formula (cf\. \cite{\BailAA, Ex. 4.(iii), p. 97, reversed})
$$
{} _{2} F _{1} \!\left [ \matrix { a, b}\\ { {1\over 2} + {a\over 2} + {b\over
   2}}\endmatrix ; {\displaystyle z}\right ] =
  {} _{2} F _{1} \!\left [ \matrix { {a\over 2}, {b\over 2}}\\ { {1\over 2} +
   {a\over 2} + {b\over 2}}\endmatrix ; {\displaystyle 4 z\left( 1 - z \right)
    }\right ] ,
$$
which is valid for $\Re z<1/2$, while to the second $_2F_1$-series 
we apply the quadratic transformation formula (\Zd), and subsequently
the transformation (cf\. \cite{\SlatAC, (1.3.15)})
$$
{} _{2} F _{1} \!\left [ \matrix { a, b}\\ { c}\endmatrix ; {\displaystyle
   z}\right ]  = {{\left( 1 - z \right) }^{c-a-b}} 
   {} _{2} F _{1} \!\left [ \matrix { c-a, c-b}\\ { c}\endmatrix ;
    {\displaystyle z}\right ] .
$$
If the resulting
expression is simplified and then substituted back in (\ZG), 
we obtain eventually the expression (\ZNa) for $\ka>2$.\quad \quad \qed
\enddemo

\remark{Remark}
(1) Other alternative expression for (\ZG) will be obtained later
(in greater generality) in Theorems~\TL\ and \TN.

(2) The sum in (\ZG) does not simplify in general. However, there are two
cases where it does, namely if $\ka=1$ and if $\ka=2$. We state these
results (again in greater generality) in Corollary~\TD.
\endremark

\subhead 5. Vicious walkers, semistandard tableaux, and jeu de taquin
\endsubhead
The purpose of this section is to show that the number of watermelon
configurations which start on the wall, end at some deviation $y$,
and have a given number of contacts with the $x$-axis is equinumerous
to another family of vicious walkers which do not run below the
$x$-axis, but have {\it no} restriction on the number of contacts
with the $x$-axis. The upshot of this is that the latter number
of watermelons can now be written
in form of a {\it single} determinant by applying the
Lindstr\"om-Gessel--Viennot determinant (\Zb). 
The result has been originally
found by Brak and Essam \cite{\BrEsAB, Cor.~1}. They derived it by
manipulating (a variant of) 
the Lindstr\"om--Gessel--Viennot determinant (\ZB).
Since the result itself is purely combinatorial,
they posed the problem of finding a bijective proof. We resolve
this problem here. The bijective proof that we shall give below
is based on the fact that vicious walker configurations are in
bijection with semistandard tableaux (see \cite{\GuOVAA, \KrGVAA})
and on a modified jeu de taquin for semistandard tableaux 
that occurred earlier 
in a bijective proof \cite{\KratBK} of Stanley's hook-content
formula for the number of semistandard tableaux of a given shape with
bounded entries. 

\proclaim{Proposition~\TI}
The number of families of $n$ vicious walkers, the $i$-th 
starting at $(0,2i-2)$ and
ending at $(t,y+2i-2)$, $i=1,2,\dots,n$, 
none of them running below the $x$-axis, and where the first walk
has $\ell+1$ contacts with the $x$-axis {\rm(}including its starting point
$(0,0)${\rm)} is the same as 
the number of families of $n$ vicious walkers, the $i$-th 
starting at $(0,2i-2)$ and
ending at $(t,y+2i-2)$, $i=1,2,\dots,n-1$, the $n$-th walker
running from $(0,2n-2)$ to $(t-\ell-1,y+2n+\ell-3)$, 
none of them running below the $x$-axis.
\endproclaim
\remark{Remark}
The statement of the proposition deviates slightly from the one in 
\cite{\BrEsAB, Cor.~1}. 
The equivalence of both statements can be seen by observing that 
in the formulation of Brak and Essam there are forced walk portions
at the beginning and at the end of the walks, which can
hence be omitted. 
There is also a typo there: the definition of ${v'}_n^{f}$ must be
$(t-m,y+2(n-1)+m)$.
\endremark
\demo{Proof of Proposition}
It is well-known that configurations of vicious walkers are in
bijection with so-called {\it semistandard tableaux}. By definition,
a semistandard tableau
is an array of integers in which entries are weakly increasing along rows
and strictly increasing along columns. Figure~\FE\ shows three examples.

\midinsert
\vskip10pt
\vbox{
$$
\smatrix\format\sa\c\s\c\s\c\s\c\se\\
\hlinefor9\\
&\hbox to10pt{\hss 1\hss}&&\hbox to10pt{\hss 2\hss}&&\hbox to10pt{\hss 2\hss}&&
\hbox to10pt{\hss 6\hss}&\\
\hlinefor9\\
&4&&6&&6&&7&\\
\hlinefor9\\
&6&&7&&7&&8&\\
\hlinefor9\\
&7&&8&&9&&10&\\
\hlinefor9\\
&9&&9&&11&&11&\\
\hlinefor9
\endsmatrix
\hskip2cm
\smatrix \format\sa\c\s\c\s\c\s\c\se\\
\hlinefor9\\
&\hbox to10pt{\hss 3\hss}&&\hbox to10pt{\hss 3\hss}&&\hbox to10pt{\hss
5\hss}&&
\hbox to10pt{\hss 6\hss}&\\
\hlinefor9\\
&4&&6&&6&&7&\\
\hlinefor9\\
&6&&7&&9&\\
\hlinefor7\\
&9&&10&&11&\\
\hlinefor7\\
&10&&12&&12&\\
\hlinefor7
\endsmatrix
\hskip2cm
\smatrix \format\sa\c\s\c\s\c\s\c\se\\
\hlinefor9\\
&\hbox to10pt{\hss 2\hss}&&\hbox to10pt{\hss 2\hss}&&\hbox to10pt{\hss 4\hss}&&
\hbox to10pt{\hss 5\hss}&\\
\hlinefor9\\
&3&&5&&5&&6&\\
\hlinefor9\\
&5&&6&&8&\\
\hlinefor7\\
&8&&9&&10&\\
\hlinefor7\\
&9&&11&&11&\\
\hlinefor7
\endsmatrix
$$
\centerline{\eightpoint a.\hskip4.2cm b.\hskip4.2cm c.}
\vskip4pt
\centerline{\eightpoint Semistandard tableaux}
\vskip7pt
\centerline{\eightpoint Figure \FE}
}
\vskip10pt
\endinsert

To explain the bijection, let us first consider the watermelon
configurations in the statement of the proposition (this is the first
set of families of vicious walkers in the statement of the proposition).
An example for $n=4$, $t=12$, $y=2$, and $\ell=3$ is shown in the
Figure~\FF.a. Now label down-steps by the $x$-coordinate of their
starting point, so that a step from $(a,b)$ to $(a+1,b-1)$ is labelled by
$a$, see Figure~\FF. Then, out of the labels of the $i$-th
walk, form the $i$-th column of the corresponding tableau. If this
procedure is applied to the family of vicious walkers in 
Figure~\FF.a then the tableau in Figure~\FE.a
is obtained. It is not difficult to see that
the resulting array of numbers is always a semistandard tableau.
In fact, the entries are
trivially strictly increasing along columns, and they are weakly
increasing along rows because the walks do not touch each other.
Moreover, the restriction that no walker may run below the $x$-axis
translates into the condition that the entries in the $i$-th row of
the corresponding tableau are at least $2i-1$. Finally, any contact of
the first walker with the $x$-axis (except the one in $(0,0)$) 
produces an entry in the first
column which is exactly equal to its lower bound. (In the tableau in
the Figure~\FE.a these are the entries 1, 7 and 9.)

\midinsert
\vskip10pt
\vbox{
$$
\Einheit.4cm
        \Gitter(13,13)(0,-1)
        \Koordinatenachsen(13,13)(0,-1)
\Pfad(0,0),343343443433\endPfad
\Pfad(0,2),334333444433\endPfad
\Pfad(0,4),334333443434\endPfad
\Pfad(0,6),333333444344\endPfad
\DickPunkt(0,0)
\DickPunkt(0,2)
\DickPunkt(0,4)
\DickPunkt(0,6)
\DickPunkt(12,2)
\DickPunkt(12,4)
\DickPunkt(12,6)
\DickPunkt(12,8)
\Label\r{\italic 1}(1,1)
\Label\r{\italic 4}(4,2)
\Label\r{\italic 6}(6,2)
\Label\r{\italic 7}(7,1)
\Label\r{\italic 9}(9,1)
\Label\r{\italic 2}(2,4)
\Label\r{\italic 6}(6,6)
\Label\r{\italic 7}(7,5)
\Label\r{\italic 8}(8,4)
\Label\r{\italic 9}(9,3)
\Label\r{\italic 2}(2,6)
\Label\r{\italic 6}(6,8)
\Label\r{\italic 7}(7,7)
\Label\r{\italic 9}(9,7)
\Label\r{\italic 1\italic1}(11,7)
\Label\r{\italic 6}(6,12)
\Label\r{\italic 7}(7,11)
\Label\r{\italic 8}(8,10)
\Label\r{\italic 1\italic0}(10,10)
\Label\r{\italic 1\italic1}(11,9)
\hbox{\hskip6cm}
        \Gitter(13,13)(0,-1)
        \Koordinatenachsen(13,13)(0,-1)
\Pfad(0,0),334434334433\endPfad
\Pfad(0,2),334334433434\endPfad
\Pfad(0,4),333344334344\endPfad
\Pfad(0,6),33333443\endPfad
\DickPunkt(0,0)
\DickPunkt(0,2)
\DickPunkt(0,4)
\DickPunkt(0,6)
\DickPunkt(12,2)
\DickPunkt(12,4)
\DickPunkt(12,6)
\DickPunkt(8,10)
\Label\r{\italic 2}(2,2)
\Label\r{\italic 3}(3,1)
\Label\r{\italic 5}(5,1)
\Label\r{\italic 8}(8,2)
\Label\r{\italic 9}(9,1)
\Label\r{\italic 2}(2,4)
\Label\r{\italic 5}(5,5)
\Label\r{\italic 6}(6,4)
\Label\r{\italic 9}(9,5)
\Label\r{\italic 1\italic1}(11,5)
\Label\r{\italic 4}(4,8)
\Label\r{\italic 5}(5,7)
\Label\r{\italic 8}(8,8)
\Label\r{\italic 1\italic0}(10,8)
\Label\r{\italic 1\italic1}(11,7)
\Label\r{\italic 5}(5,11)
\Label\r{\italic 6}(6,10)
\hskip4.8cm
$$
\vskip-4pt
\centerline{\eightpoint a.\hskip5cm b.}
\vskip4pt
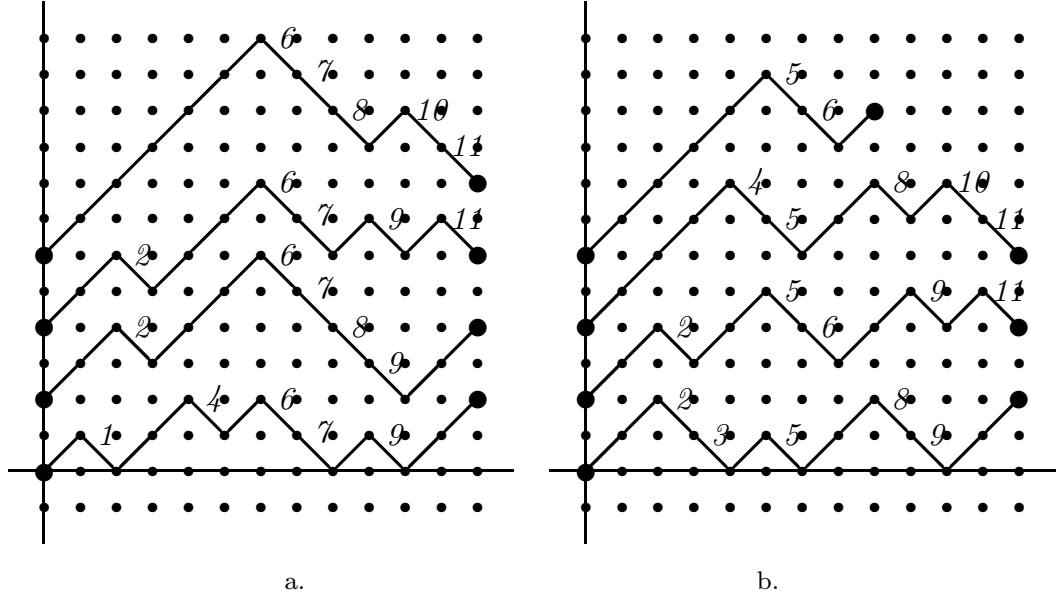
\centerline{\eightpoint Corresponding vicious walkers}
\vskip7pt
\centerline{\eightpoint Figure \FF}
}
\vskip10pt
\endinsert

Thus, the first set of families of vicious walkers in the statement of
the proposition is in bijection with semistandard tableaux with $n$
columns of length $(t-y)/2$, with the entries in row $i$ at least
$2i-1$ and at most $t-1$, $i=1,2,\dots,(t-y)/2$, and with exactly $\ell$
entries in the first column which are equal to their respective lower bound.

If we apply this same translation to the second set 
of families of vicious walkers in the statement of
the proposition, then we see that these are in bijection with 
semistandard tableaux with $n$ columns, the first $n-1$ of them
of length $(t-y)/2$, the {\it last\/} of length $(t-y-2\ell)/2$, 
with the entries in row $i$ being at least
$2i-1$ and at most $t-1$, $i=1,2,\dots,(t-y)/2$, and with the entries 
in the {\it last\/} column being at most $t-\ell-2$.
Figure~\FE.c shows the result when this translation
is applied to the family of vicious walkers in 
Figure~\FF.b. 

It is now our task to construct a bijection between the latter two
sets of tableaux. We start with a tableau from the first set (of which
an example is shown in Figure~\FE.a).
The construction makes use of a modified
form of jeu de taquin from \cite{\KratBK, (2.2)--(2.4)}.
We start with the bottom-most entry in the first column which is equal to
its lower bound. (In Figure~\FE.a this is the entry 9
in the first column.) Let us call it the ``special entry," and denote
it by $s$. Now we move $s$ to the right/bottom by the following
procedure, denoted by (JT):

\smallskip
(JT) Compare the special entry $s$ with its right neighbour,
$x$ say (if there is no right neighbour, then, by convention, 
we set $x=\infty$), and its bottom neighbour, $y$ say 
(if there is no bottom neighbour, then, by convention, we set $y=\infty$), 
see (\Zg). If there is no right or bottom neighbour, then stop.

If not, then we have the following
situation,
$$\smatrix \format \sa\c\s\c\se\\
\hlinefor5\\
&s&&x&\\
\hlinefor5\\
&y&\\
\hlinefor3
\endsmatrix\ ,
\tag\Zg$$
where at least on of $x$ and $y$ is not $\infty$.
If $x\ge y-1$ then do the move
$$\smatrix \format \sa\c\s\c\se\\
\hlinefor5\\
&y-1&&x&\\
\hlinefor5\\
&s&\\
\hlinefor3
\endsmatrix\ .\tag\Zi
$$
If $x+1<y$ then do the move
$$\smatrix \format \sa\c\s\c\se\\
\hlinefor5\\
&x+1&&s&\\
\hlinefor5\\
&y&\\
\hlinefor3
\endsmatrix\ .\tag\Zh
$$
Repeat (JT). 

\smallskip

For example, if we apply this algorithm to the tableau in
Figure~\FE.a, the special entry being the 9 in the first column, 
then we obtain the tableau in Figure~\FG.
(In this case it is only movements of type (\Zh) that are applied.)

\midinsert
\vskip10pt
\vbox{
$$
\smatrix\format\sa\c\s\c\s\c\s\c\se\\
\hlinefor9\\
&\hbox to10pt{\hss 1\hss}&&\hbox to10pt{\hss 2\hss}&&\hbox to10pt{\hss 2\hss}&&
\hbox to10pt{\hss 6\hss}&\\
\hlinefor9\\
&4&&6&&6&&7&\\
\hlinefor9\\
&6&&7&&7&&8&\\
\hlinefor9\\
&7&&8&&9&&10&\\
\hlinefor9\\
&10&&12&&12&&9&\\
\hlinefor9
\endsmatrix
$$
\vskip7pt
\centerline{\eightpoint Figure \FG}
}
\vskip10pt
\endinsert

\noindent
Clearly, after this algorithm has terminated, the special entry will
have ended up in the bottom-right corner (as it does in this example).
We now remove the special entry (the 9 in the bottom-right corner in
our example). 

Next, the same algorithm is applied to the (now) bottom-most entry in
the first column which is equal to its lower bound. 
(In Figure~\FG\ this is the entry 7
in the first column.) Since successive jeu de taquin paths cannot
cross each other (see \cite{\KratBK, Fig.~9} and the accompanying
explanations for a detailed argument), this entry must again end up in the
bottom of the last column. Subsequently, we remove it.
This procedure is repeated until there is no entry in the first
column which is equal to its lower bound. (The final result when this
procedure is applied to the tableau in Figure~\FE.a is the tableau in
Figure~\FE.b.)

Now, what do we obtain in the end? It is easy to see that after each
movement of type (\Zi) or (\Zh) the entries are weakly increasing
along rows and strictly increasing along columns, if we ignore the
special entry. Thus, we will have obtained a semistandard tableau with
$n$ columns, out of which all of them have length $(t-y)/2$ except for
the last, which will be by $\ell$ shorter. Furthermore, all entries will
be at most $t$ (this is due to the fact that at most once 1 is added
to an entry by performing a movement of type (\Zh)), and all entries
in the $i$-th row will now be {\it strictly larger} than $2i-1$. Finally,
since originally the entry in row $i$ in the last column is at most 
$(t+y+2i-2)/2$ due to strict increase of entries along columns, 
and since whenever an entry (other than the special
entry) is moved upward (which is only possible through a movement of
type (\Zi)) it shrinks by 1, the entries in the last column will in
the end be at most $t-\ell-1$. Clearly, if we subtract 1 from all the
entries of such a tableau, we obtain a tableau from the second set
of tableaux. (In our running example, after subtracting 1 from all
the entries in the tableau in Figure~\FE.b we obtain the tableau in
Figure~\FE.c.)

We claim that this is a bijection between the sets of tableaux under
consideration. In order to see this, we have to construct the inverse
mapping. In fact, every single step can be inverted. We would start
with a tableau from the second set, add 1 to all the entries,
place an entry $s$ (whose value
will be determined later) immediately below the bottom-most entry in
the last column of the tableau, and then apply the following
``inverse" jeu de taquin, denoted by (JT*):

\smallskip

(JT*) 
Call $s$ the special entry.

If the special entry $s$ is located in the first row and column, then stop.

If not, then we have the following
situation,
$$\smatrix \format \sa\c\s\c\se\\
\omit&\omit&\hlinefor3\\
\omit&&&y&\\
\hlinefor5\\
&x&&s&\\
\hlinefor5
\endsmatrix\ ,
$$
where one of $x$ or $y$ could also be absent.
(If $y$ is actually not there, then by convention
we set $y=-\infty$.)
We denote the number of the row in which $s$ is located by $i$.

If $x$ is actually not there and $y+1\le 2i-1$, then stop. 
Otherwise, and also if 
$x\le y+1$, then do the move
$$\smatrix \format \sa\c\s\c\se\\
\omit&\omit&\hlinefor3\\
\omit&&&s&\\
\hlinefor5\\
&x&&y+1&\\
\hlinefor5
\endsmatrix\ .\tag\Zl
$$
If $x-1>y$ then do the move
$$\smatrix \format \sa\c\s\c\se\\
\omit&\omit&\hlinefor3\\
\omit&&&y&\\
\hlinefor5\\
&s&&x-1&\\
\hlinefor5
\endsmatrix\ .\tag\Zk
$$
Repeat (JT*).

\smallskip
Suppose that the special entry ended up in row $i$. Then we now set
$s=2i-1$. Subsequently, this procedure is repeated until the complete 
$((t-y)/2)\times n$ rectangle is filled. We leave it to the reader to
check that this yields indeed the desired inverse mapping.
\quad \quad \qed
\enddemo

\remark{Remark}
Brak and Essam have also a variant of Proposition~\TI\ in
\cite{\BrEsAB, Cor.~2}. It states that 
the number of families of $n$ vicious walkers, the $i$-th 
starting at $(0,2i-2)$ and
ending at $(t,y+2i-2)$, $i=1,2,\dots,n$, 
none of them running below the $x$-axis, and where $\ell$
of the contacts of the first walk with the $x$-axis (excluding 
$(0,0)$) are marked, is the same as 
the number of families of $n$ vicious walkers, the $i$-th 
starting at $(0,2i-2)$ and
ending at $(t,y+2i-2)$, $i=1,2,\dots,n-1$, the $n$-th walk
running from $(0,2n-2)$ to $(t,y+2n+2\ell-2)$, 
none of them running below the $x$-axis. Again, they prove it by
determinant manipulations and pose the problem of finding a bijective
proof. Such a bijective proof can be provided in the same manner as
the one that we found for Proposition~\TI. First, one translates both
sets of walk families into sets of semistandard tableaux in the same
manner as before. The first set translates into semistandard tableaux
with $n$
columns of length $(t-y)/2$, with the entries in row $i$ at least
$2i-1$ and at most $t-1$, $i=1,2,\dots,(t-y)/2$, and with exactly $\ell$
entries in the first column which are equal to its lower bound and
which are marked.
The second set translates into 
semistandard tableaux with $n$ columns, the first $n-1$ of them
of length $(t-y)/2$, the {\it last\/} of length $(t-y-2\ell)/2$, 
with the entries in row $i$ being at least
$2i-1$ and at most $t-1$, $i=1,2,\dots,(t-y)/2$. The bijection of the
proof of Proposition~\TI\ then becomes a bijection between the latter two
sets of tableaux when one replaces the modified jeu de taquin (JT) and
(JT*) by {\it ordinary} jeu de taquin (as for example explained in
\cite{\SagaAL}; i.e., nothing is subtracted from or added to $x$ and
$y$ in (\Zi), (\Zh), (\Zl), (\Zk); furthermore the condition $x\le
y+1$ has to be replaced by $x\le y$ everywhere, 
and $x-1>y$ has to be replaced by
$x>y$ everywhere, 
as well as the stopping rule in (JT*) has to be modified to the extent
that $y+1\le 2i-1$ has to be replaced by $y<2i-1$).
\endremark

\subhead 6. Exact formulas for the partition function 
for watermelons of arbitrary deviation\endsubhead
This section is devoted to the derivation of exact formulas for the
partition function $Z^{(n)}_t(y;\ka)$ for 
families of $n$ vicious walkers, the $i$-th 
starting at $(0,2i-2)$ and
ending at $(t,y+2i-2)$, $i=1,2,\dots,n$, 
none of them running below the $x$-axis, as defined in (\PP). 
The main result is Theorem~\TK, which expresses the partition function
in terms of a double sum. It extends Theorem~\TC\ to watermelon
configuration of an
arbitrary (but fixed) deviation. By applying hypergeometric summation and
transformation formulas to this double sum, 
we are also able to provide alternative expressions for
$Z^{(n)}_t(y;\ka)$, special cases of which having been given earlier 
by Brak and Essam (private communication).

Our starting point is Proposition~\TI.
We apply the Lindstr\"om--Gessel--Viennot formula (\Zb) to
the second set of vicious walkers in the
proposition in order to obtain a determinant for the number of
watermelons of deviation $y$ and exactly $\ell+1$ contacts with the
$x$-axis. The most useful determinant is obtained if we first prepend 
$2i-2$ up-steps to the $i$-th walk (analogously to what we did in the
proof of Proposition~\TB), so that we obtain families of $n$
vicious walkers, the $i$-th starting at $(-2i+2,0)$ and
ending at $(t,y+2i-2)$, $i=1,2,\dots,n-1$, the $n$-th walker
running from $(-2n+2,0)$ to $(t-\ell-1,y+2n+\ell-3)$, 
none of them running below the $x$-axis. Again, it is clear that the
prepended walk portions are forced, and, hence, the number of the
latter families of walkers are the same as the number of families in
the second set of walkers in Proposition~\TI. If we apply the
Lindstr\"om--Gessel--Viennot formula to the modified families of
walkers, then we obtain
$$\det_{1\le i,j\le n}\pmatrix 
\frac {y+2i-1} {\frac {t+y} {2}+i+j-1}\binom {t+2j-2}{\frac {t+y}
{2}+i+j-2}&\quad 1\le i\le n-1\\
\frac {y+2n+\ell-2} {\frac {t+y} {2}+n+j-2}\binom {t+2j-\ell-3}{\frac {t+y}
{2}+n+j-3}&\quad i=n
\endpmatrix,\tag\Zm$$
where we used (\ZZa) to compute the number of walks with given
starting and end point, which do not pass below the $x$-axis.

As it turns out, this determinant can be simplified to a single
sum. We state the corresponding result in the lemma below.

\proclaim{Lemma~\TJ}
Let $t$ and $y$ be non-negative integers with $t\equiv y$ {\rm(mod 2)}.
Then the number of families of $n$ vicious walkers, the $i$-th 
starting at $(0,2i-2)$ and
ending at $(t,y+2i-2)$, $i=1,2,\dots,n$, 
none of them running below the $x$-axis, and where the first walk
has $\ell+1$ contacts with the $x$-axis {\rm(}including its starting point
$(0,0)${\rm)} is equal to
$$\multline \prod _{i=0} ^{n-1}\frac {(t+2i)!\,i!} {(\frac {t+y}
{2}+n+i-1)!\, (\frac {t-y} {2}+i)!}\prod _{i=0} ^{n-2}\frac
{(y+n+i-1)!} {(y+2i)!}\\
\times \sum _{k=0} ^{n-1}(-1)^{n-k-1}\frac {(y+2n+\ell-2)\,(t+2k-\ell-1)!\,
(\frac {t-y} {2}+k)!} {(t+2k)!\,(\frac {t-y} {2}+k-n-\ell+1)!\, k!\,(n-k-1)!}.
\endmultline\tag\Zn$$
\endproclaim

\demo{Proof}
We expand the determinant (\Zm) along the last row. Thus we obtain
$$\sum _{k=1} ^{n}(-1)^{n+k}\frac {y+2n+\ell-2} {\frac {t+y} {2}+n+k-2}\binom
{t+2k-\ell-3}{\frac {t+y}
{2}+n+k-3}M_k,$$
where $M_k$ is the minor of (\Zm) where the last row and the $k$-th
column has been deleted. This minor is easily evaluated by means of
\cite{\KratBN, Theorem~26, (3.13)} with $A=1-y$, $B=2$, 
and $L_i=i+\chi(i\ge k)+\frac {t+y}
{2}-1$, where $\chi(A)=1$ if $A$ is true and $\chi(A)=0$ otherwise. 
After some simplification, and after replacing
$k$ by $k+1$, one obtains (\Zn).\quad \quad \qed
\enddemo

The preceding lemma provides us with a workable expression for
the number of watermelon configurations with a given number of
contacts with the $x$-axis. We will use it now to derive several
expressions for the partition function 
$Z^{(n)}_t(y;\ka)$. The first of those will
be the one which is most suited for studying the asymptotic behaviour of
the partition function, and for the mean number of contacts, which we
will do in Sections~7 and 8, respectively.

\proclaim{Theorem \TK}
Let $t$ and $y$ be non-negative integers with $t\equiv y$ {\rm(mod 2)}.
The partition function $Z^{(n)}_t(y;\ka)$ for families of 
$n$ vicious walkers, the $i$-th starting at $(0,2i-2)$ and
ending at $(t,y+2i-2)$, $i=1,2,\dots,n$, 
none of them running below the $x$-axis, where the
weight of a walker configuration $\bold P$ is defined as 
$\ka^{c(\bold P)}$ with $c(\bold P)$
denoting the number of contacts of the walkers with the $x$-axis, is
given by
$$\multline 
\frac {(\frac {t+y}2-1)!}{(\frac {t+y}2+2n-2)!}
\prod _{i=0} ^{n-2}\frac {(t+2i)!\,i!\,(y+n+i-1)!}
{(\frac {t+y}2+n+i-1)!\,
   (\frac {t-y}2+i)!\,(y+2i)!}\\
\times \sum _{\ell=0} ^{(t-y)/2}
\sum _{k=0} ^{n-1}(-1)^k
   \binom {n-1}{k}\binom {y+2k-2}{2k}
\binom {t-\ell-1}{\frac {t-y}2-\ell-k}\hskip2cm\\
\cdot
   \binom {y+\ell+2n-2}{2n-2k-1}
\frac {\big(\frac {t+y}2\big)_k\, (2n-2k-1)!\,(2k)!}
{\big(\frac {t-y}2+n-k\big)_k}\ka^{\ell+1}.
\endmultline\tag\Zo$$
\endproclaim

\remark{Remark}
In the case $n=1$ the product in (\Zo) is empty and has to be
interpreted as $1$. Moreover, because of the binomial coefficient
$\binom {n-1}k$ in the summand, the sum over $k$ collapses to just
one term, the term for $k=0$.
\endremark

\demo{Proof}
Clearly, the partition function $Z^{(n)}_t(y;\ka)$ is equal
to 
$$\sum _{\ell=0} ^{(t-y)/2}N^{(n)}_t(y;\ell)\ka^{\ell+1},\tag\Zp$$ 
where $N^{(n)}_t(y;\ell)$ is the expression in (\Zn). 
We write the sum in (\Zn) in hypergeometric notation, thus obtaining
$$\multline 
 {\left( -1 \right) }^{n+1}\frac{\left(  y+ 2 n + \ell -2\right)  
      (t-\ell-1)! \,
      (\frac{t-y}{2} )! }{
(n-1)! \,t! \,
      (\frac{t-y}{2} - n -\ell + 1)! }
\\ 
\times
      {} _{4} F _{3} \!\left [ \matrix { 1 + \frac{t}{2} - \frac{y}{2}, 
 -\frac{\ell}{2}+ \frac{t}{2}, \frac{1}{2} - \frac{\ell}{2} + \frac{t}{2}, 1 - 
n}\\ { 1 + \frac{t}{2}, \frac{1}{2} + \frac{t}{2}, 2 - \ell - n + \frac{t}{2} 
- \frac{y}{2}}\endmatrix ; {\displaystyle 1}\right ]  .
\endmultline$$
To the $_4F_3$-series we apply one of the balanced 
(which, by definition, means that
the sum of the denominator parameters exceeds the sum of the numerator
parameters by 1) $_4F_3$ 
transformation formulas due to Whipple
(cf\. \cite{\GaRaAA, Appendix (III.16), $q\uparrow1$})
$$\multline
{} _{4} F _{3} \!\left [ \matrix { a, b, c, -N}\\ { e, f, 1 + a + b + c - e -
   f - N}\endmatrix ; {\displaystyle 1}\right ] \\
=
{{( a)_N\,( -a - b + e + f)_N\,( -a - c + e + f) _{N}}\over
     {( e)_N\, (f)_N\, (-a - b - c + e + f) _{N}}}\\
\times
  {} _{4} F _{3} \!\left [ \matrix { -N, -a + e, -a + f, -a - b - c + e + f}\\
    { -a - b + e + f, -a - c + e + f, 1 - a - N}\endmatrix ; {\displaystyle
    1}\right ]   ,
\endmultline\tag\Zpa$$
where $N$ is a non-negative integer. Thus we obtain
$$\multline 
\frac{\left( y+ 2 n +\ell-2 \right) \,
    ({ \textstyle n}) _{t-n-\ell} \,
    ({ \textstyle  t+2 n -1}) _{1 - n - \frac{t}{2} - \frac{y}{2}} \,
    ({ \textstyle y+\ell}) _{ 2 n-2} }{(\frac{t-y}{2} -\ell)! }\\
\times
    {} _{4} F _{3} \!\left [ \matrix { 1 - n, \frac{y}{2}, 
   \frac{y}{2}-\frac{1}{2} , \ell - \frac{t}{2} + \frac{y}{2}}\\ {
\frac{1}{2} + \frac{\ell}{2} + \frac{y}{2}, \frac{\ell}{2} + \frac{y}{2}, 1 - n - 
\frac{t}{2} + \frac{y}{2}}\endmatrix ; {\displaystyle 1}\right ] 
\endmultline$$
for the sum in (\Zn). If we substitute this back in (\Zp), then we
obtain the claimed expression after some manipulation.\quad \quad \qed
\enddemo

\remark{Remark}
This theorem generalizes Theorem~\TC\ to watermelon configurations of
arbitrary (fixed) deviation. Indeed,
if we set $y=0$ in Theorem~\TK, then, because of 
the binomial coefficient $\binom {y+2k-2}{2k}$ appearing
in the summand in (\Zo), the only summands which are nonzero in (\Zo)
are the ones with $k=0$.
\endremark

As corollaries of Lemma~\TJ\
(which via (\Zpa) is equivalent to Theorem~\TK),
we are able to derive remarkable simple alternative 
expressions for the partition
function $Z^{(n)}_t(y;\ka)$ in terms of {\it single} sums. 
(These expressions however seem to be less suited
for asymptotic considerations.) We start with an expansion
of the partition function around $\ka=1$.

\proclaim{Theorem \TL}
With the same assumptions as in Theorem~\TK, the partition function
$Z^{(n)}_t(y;\ka)$ is equal to
$$\multline
 \Bigg(\prod _{i=0} ^{n-1}(t+2i)!\, i!\Bigg)
\Bigg(\prod _{i=0} ^{n-2}\frac {(y+n+i-1)!}
{(y+2i)!\,(\frac {t+y}2+n+i)!\,(\frac {t-y}2+i+1)!}\Bigg)\\
\times
\sum _{h=0} ^{(t-y)/2}\binom{n+h-1}{n-1}
\frac {(y+2n+2h-1) 
\,  (y+2n+h-2)!}{(y+n+h-1)!\,(\frac {t+y}2+2n+h-1)!\, (\frac {t-y}2-h)!} 
\ka(\ka-1)^h.
\endmultline\tag\Zq$$
\endproclaim

\demo{Proof}
Again, by Lemma~\TJ, the partition function $Z^{(n)}_t(y;\ka)$ is
equal to (\Zp), with $N^{(n)}_t(y;\ell)$ the expression in (\Zn). 
We write $\ka^{\ell+1}=\ka\sum _{h=0} ^{\ell}\binom \ell h(\ka-1)^h$,
and substitute this in (\Zp). Thus we obtain
$$\multline
\Pi^{(n)}_t(y)\sum _{\ell=0} ^{(t-y)/2}
\sum _{k=0} ^{n-1}(-1)^{n-k-1}\frac {(y+2n+\ell-2)\,(t+2k-\ell-1)!\,
(\frac {t-y} {2}+k)!} {(t+2k)!\,(\frac {t-y} {2}+k-n-\ell+1)!\,
k!\,(n-k-1)!}\\
\cdot
\ka\sum _{h=0} ^{\ell}\binom \ell h(\ka-1)^h,
\endmultline$$
where $\Pi^{(n)}_t(y)$ is the product in the first line of (\Zn).
Now we interchange summations, so that the sum over $\ell$ becomes 
the inner-most sum and the sum over $h$ becomes the outer-most sum, 
and we convert the sum over $\ell$ to
hypergeometric notation. This yields
$$\multline
\Pi^{(n)}_t(y)\frac {\ka} {(n-1)!}
\sum _{h=0} ^{(t-y)/2}\sum _{k=0} ^{n-1}
(-1)^{n-k-1}\binom {n-1}k\\
\cdot\frac {(t+2k-h-1)!\,
(\frac {t-y} {2}+k)!} {(t+2k)!\,(\frac {t-y} {2}+k-n-h+1)!}(\ka-1)^h\\
\cdot
\bigg((y+2n-3)\,{} _{2} F _{1} \!\left [ \matrix h+1,-\frac {t}
{2}+\frac {y} {2}-k+n+h-1\\-t-2k+h+1
\endmatrix ; {\displaystyle 1}\right ]\\
+
(h+1)\,{} _{2} F _{1} \!\left [ \matrix h+2,-\frac {t} {2}+\frac {y}
{2}-k+n+h-1 \\-t-2k+h+1
\endmatrix ; {\displaystyle 1}\right ]
\bigg).
\endmultline$$
Both $_2F_1$-series can be evaluated by means of the 
Chu--Vandermonde summation formula (cf\. \cite{\SlatAC, (1.7.7); Appendix
(III.4)})
$$
{} _{2} F _{1} \!\left [ \matrix { a, -N}\\ { c}\endmatrix ; {\displaystyle
   1}\right ]  = {{({ \textstyle c-a}) _{N} }\over 
    {({ \textstyle c}) _{N} }},
\tag\ZP$$
where $N$ is a non-negative integer. The result is simplified,
and the sum over $k$ is written in hypergeometric notation. This gives
the expression
$$\multline
\Pi^{(n)}_t(y)\frac {\ka} {(n-1)!}
\sum _{h=0} ^{(t-y)/2}
(-1)^{n-1}
\frac {(y+2n+2h-1)\,(\frac {t+y} {2}+n-1)!\,
(\frac {t-y} {2})!} {(\frac {t+y} {2}+n+h)!\,
(\frac {t-y} {2}-n-h+1)!}(\ka-1)^h\\
\cdot
{} _{3} F _{2} \!\left [ \matrix \frac {t} {2}+\frac {y} {2}+n,
\frac {t} {2}-\frac {y} {2}+1,1-n\\
\frac {t} {2}+\frac {y} {2}+n+h+1,\frac {t} {2}-\frac {y} {2}-n-h+2
\endmatrix ; {\displaystyle 1}\right ].
\endmultline$$
The $_3F_2$-series is balanced and can therefore be evaluated by
means of the Pfaff--Saalsch\"utz summation formula
(cf\. \cite{\SlatAC, (2.3.1.3); Appendix (III.2)})
$$
{} _{3} F _{2} \!\left [ \matrix { a, b, -N}\\ { c, 1 + a + b - c -
   N}\endmatrix ; {\displaystyle 1}\right ]  =
  {{({ \textstyle c-a}) _{N} \, ({ \textstyle c-b}) _{N} }\over 
    {({ \textstyle c}) _{N} \, ({ \textstyle c-a-b}) _{N} }},
$$
where $N$ is a non-negative integer. Subsequent simplification then
leads to the claimed expression (\Zq).\quad \quad \qed
\enddemo

\remark{Remarks}
(1) The coefficient of $\ka(\ka-1)^h$ in (\Zq) has a combinatorial
interpretation in terms of watermelon
configurations, the $i$-th walker of the configuration running
from $(0,2i-2)$ to $(t,y+2i-2)$, $i=1,2,\dots,n$, but never below
the $x$-axis, where contacts with the $x$-axis other than $(0,0)$ 
can be marked or not. It is an immediate result
of the comparison of (\Zq) with the definition of $Z^{(n)}_t(y;\ka)$.
In this definition, we replaced $\ka^{\ell+1}$ by $\ka\sum _{h=0}
^{\ell}\binom \ell h(\ka-1)^h$. Hence, the coefficient of 
$\ka(\ka-1)^h$ in (\Zq) counts watermelon configurations as above
where {\it exactly $h$ contacts} other than $(0,0)$ 
{\it are marked}.

(2) It is also possible to prove Theorem~\TL\ in the same (automated)
fashion as Theorem~\TC, i.e., by starting
with the Lindstr\"om--Gessel--Viennot determinant (\ZB), and then arguing
by induction, which one would base on
the condensation formula (\ZD), the verification of 
the binomial summations that
one has to establish on the way being accomplished by the
Gosper--Zeilberger algorithm.
\endremark

The sum in (\Zq) does not simplify in general. In fact, when
written in hypergeometric terms, it turns out to be a $_4F_3$-series,
$$
\frac{\ka\,
 (y + 2 n -1)! }{
    (\frac{t-y}{2} )! \,
    (\frac{t+y}{2}+2 n  -1)! \,
    (y+ n -1)! }
    {} _{4} F _{3} \!\left [ \matrix { y+ 2 n-1, \frac{y}{2} + n + 
\frac{1}{2} , n, -\frac{t}{2} + \frac{y}{2}}\\ {  \frac{y}{2}
 + n - \frac{1}{2}, y + n, 2 n + \frac{t}{2} + \frac{y}{2}}\endmatrix ; 
{\displaystyle 1 - \ka}\right ].
\tag\ZPa$$
This is a balanced $_4F_3$-series, 
and there are many identities known for balanced
$_4F_3$-series, but there is no {\it summation} formula available that
would apply in full generality. However, there are two
cases where summation formulas are available, 
namely if $\ka=1$ (trivially), and if $\ka=2$. 
The result
corresponding to $\ka=1$, giving the total {\it number} of the
watermelon configurations that we consider in this  
section (and the subsequent ones), has
been previously obtained by the author, Guttmann and Viennot 
in \cite{\KrGVAA, Theorem~6}. In fact, in
Theorem~6 of \cite{\KrGVAA} it is shown how to derive a more general result 
which is
valid for star configurations.\footnote{What is shown in \cite{\KrGVAA} is
that the result for star configurations is, in equivalent form, 
part of the folklore of
combinatorics and representation theory of symplectic groups. A 
weighted generalization of this result, in which each path configuration
is assigned a certain $q$-weight, had been proven already in 
\cite{\KratAP, Theorem~7}.}
The result corresponding to $\ka=2$
has been previously obtained by Ciucu and the author 
in \cite{\CiKrAB, Theorem~1.4} in an equivalent form,\footnote{The
combinatorial objects that are studied in \cite{\CiKrAB} are rhombus
tilings of regions composed of equilateral unit triangles. There is a
standard bijection between such tilings and families of
non-intersecting lattice paths, 
which is, for example, explained in
Section~2 of \cite{\CiKrAB}, see Figures~2.1, 2.2, 2.3(a) in that
paper. In the particular
case addressed by Theorem~1.4 of \cite{\CiKrAB}, a 45$^\circ$ rotation
turns these non-intersecting lattice paths into our vicious walkers.
Under these transformations, the weight for rhombus tilings defined
in \cite{\CiKrAB} becomes our contact weight for vicious walkers, up
to a multiplicative constant.} and, in the case $y=0$,
independently by Owczarek, Essam and Brak in
\cite{\OwEBAA, Eq.~(4.45)}.
The proof that we give
below is independent from the proofs in the papers mentioned. 

\proclaim{Corollary \TD}Let $t$ and $y$ be non-negative integers with
$t\equiv y$ {\rm(mod 2)}. Then
the total number $Z^{(n)}_{t}(y;1)$ of families of
$n$ vicious walkers, the $i$-th starting at $(0,2i-2)$ and
ending at $(t,y+2i-2)$, $i=1,2,\dots,n$, 
none of them running below the $x$-axis, is
equal to
$$
\prod _{i=0} ^{n-1}\frac {(y+n+i)!\,(t+2i)!\, i!}
{(y+2i)!\,(\frac {t+y}2+n+i)!\,(\frac {t-y}2+i)!}.
\tag\ZL$$
The partition function $Z^{(n)}_{t}(y;2)$ for these families of
vicious walkers is equal to
$$
2
\prod _{i=0} ^{n-1}\frac {(y+n+i-1)!\,(t+2i)!\, i!}
{(y+2i)!\,(\frac {t+y}2+n+i-1)!\,(\frac {t-y}2+i)!}.
\tag\ZM$$
\endproclaim

\demo{Proof} The claim (\ZL) is immediately obvious once we set
$\ka=1$ in (\Zq).

For proving the second claim, we set $\ka=2$ in (\Zq) and write
the sum in (\Zq) in hypergeometric notation (equivalently, we set
$\ka=2$ in (\ZPa)). Thus we see that this
sum is equal to 
$$\frac{2 
\,(y + 2 n -1)! }
{\big(\frac{t-y}{2}\big)! \,
    \big(\frac{t+y}{2}+2n-1\big)! \,
    (y + n -1)! }
    {} _{4} F _{3} \!\left [ \matrix { y + 2 n -1, \frac{y}{2} + n + 
\frac{1}{2} ,  n,-\frac{t}{2} + \frac{y}{2}}\\ {  \frac{y}{2}
 + n - \frac{1}{2},  y + n,2 n + \frac{t}{2} + \frac{y}{2}}\endmatrix ; 
{\displaystyle -1}\right ].
$$
The $_4F_3$-series can be evaluated by means of the summation formula
(cf\. \cite{\SlatAC, (2.3.4.6); Appendix
(III.10)})
$$
{} _{4} F _{3} \!\left [ \matrix { a, 1 + \frac{a}{2}, b, c}\\ { \
\frac{a}{2}, 1 + a - b, 1 + a - c}\endmatrix ; {\displaystyle -1}\right ]  = 
  \frac{\Gamma({ \textstyle 1 + a - b}) \,\Gamma({ \textstyle 1 + a - c}) }
   {\Gamma({ \textstyle 1 + a}) \,\Gamma({ \textstyle 1 + a - b - c}) }.
$$
Substitution of the result in (\Zq) then leads to the claimed
expression after some simplification.\quad \quad \qed
\enddemo

Our next result provides an alternative
expression for $Z^{(n)}_t(y;\ka)$ in form of an {\it infinite} series,
which resembles the one in Corollary~\TCa\ for the $y=0$ case.
Again, the corresponding formula is only valid for $\ka<2$. 
As we are going to show, it follows readily from the formula in
Theorem~\TL\ by applying another well-known hypergeometric
transformation formula. Moreover, by the same techniques we can also
obtain a similar expression which is valid for $\ka>2$.
We presume that, in analogy to the analysis in \cite{\OwEBAA}, these
expressions will be the right starting point in order to carry out 
a scaling analysis of
$Z^{(n)}_t(y;\ka)$ near the critical point $\ka=2$.

\proclaim{Theorem \TN}
Under the assumptions of Theorem~\TK, the partition function\linebreak
$Z^{(n)}_t(y;\ka)$ is equal to
$$\multline 
\Bigg(\prod _{i=0} ^{n-2}\frac {(y+n+i-1)!}{(y+2i)!}\Bigg)
\Bigg(\prod _{i=0} ^{n-1}\frac {(t+2i)! \,i!} 
{(\frac {t-y}2+i)!\,(\frac {t+y}2+n+i-1)!}\Bigg)
\\
\times
\frac {(2-\ka)}{\ka^{2n+y-1}}
\sum _{h=0} ^{\infty}\frac {(y+2n+2h-1)!\,(\frac {t+y}2+n+h-1)!} 
{h!\,(y+n+h-1)!\,(\frac {t+y}2+2n+h-1)!}
 \left(\frac {\ka-1}{\ka^2}\right)^h
\endmultline\tag\Zr$$
if $2(\sqrt2-1)<\ka<2$, and it is equal to
$$\multline 
\Bigg(\prod _{i=0} ^{n-1}\frac {(t+2i)!\,i!} 
{\big(\frac {t-y}2+i\big)!\,\big(\frac {t+y}2+n+i-1\big)!}\Bigg)
\Bigg(\prod _{i=0} ^{n-2}\frac {(y+n+i-1)!}{(y+2i)!}\Bigg)(\ka-2)\\
\times\Bigg(
\frac { 1}{\ka^{2n+y-1}}
 \sum _{h=0} ^{\infty}\frac {(y+2n+2h-1)!\,\big(\frac {t+y}2+n+h-1)!}
{\big(\frac {t+y}2+2n+h-1)!\, 
    (y+n+h-1)!\,h!} \(\frac {\ka-1}{\ka^2}\)^h\kern2cm\\
+\frac {\ka^{2n+t-1}}{(\ka-1)^{\frac {t+y}2+2n-1}}
 \sum _{h=0} ^{n-1}\frac {\big(\frac {t+y}2+2n-h-2\big)!\,
\big(\frac {t-y}2+n-h-1\big)!}{h!\,(n-h-1)!\,(t+2n-2h-2)!}
    \(\frac {1-\ka}{\ka^2}\)^h 
\Bigg)\\
\endmultline\tag\Zra$$
if $\ka>2$.
\endproclaim

\demo{Proof}
We work with the expression (\Zq) for the partition function
$Z^{(n)}_t(y;\ka)$.
In hypergeometric notation, the sum in (\Zq) is equal to (\ZPa).

Now let first $\ka<2$.
Then to the 
$_4F_3$-series in (\ZPa) we apply the quadratic transformation formula
(cf\. \cite{\BailAA, Ex. 6, p. 97, reversed})
$$\multline
{} _{4} F _{3} \!\left [ \matrix { a, \frac {a} {2}+1, b, c}\\
{\frac {a} {2},1+a-b,1+a-c}\endmatrix; {\displaystyle z}\right] \\
=
  {{\left( 1 + z \right)  }\over{{\left( 1 - z \right) }^{a+1}}}
   {} _{3} F _{2} \!\left [ \matrix { {1\over 2} + {a\over 2}, 1 + {a\over 2},
    1 + a - b - c}\\ { 1 + a - b, 1 + a - c}\endmatrix ; {\displaystyle
    -{{4 z}\over {{{\left( 1-z \right) }^2}}}}\right ] ,
\endmultline\tag\Zs$$
which is valid provided $\vert z\vert<1$ and $\vert 4z/(1-z)^2\vert
<1$.
Thus, as long as $2(\sqrt2-1)<\ka<2$, the sum in (\Zq) is equal to
$$\multline
\frac{\left( 2 - \ka \right)     }
{ \ka^{y+2n-1}}
\frac {  (y + 2 n -1)!} 
{(\frac{t-y}{2})! \,(\frac{t+y}{2}+ 2 n -1)! 
\,(y + n -1)! }\\
\times
    {} _{3} F _{2} \!\left [ \matrix {  \frac{y}{2}+n, 
\frac{y}{2}+n+\frac {1} {2}, n + \frac{t}{2} + \frac{y}{2}}\\ { 2 n 
+ \frac{t}{2} + \frac{y}{2}, n + y}\endmatrix ; {\displaystyle 
\frac{4 \left(  \ka-1 \right) }{\ka^ 2}}\right ]  .
\endmultline$$
This expression is now substituted back in (\Zq). Some simplification
then yields the expression (\Zr).

Next we address the case $\ka>2$. Clearly, we cannot follow the same
line of argument because we would not be allowed to apply (\Zs)
because of $\vert 1-\ka\vert>1$. Instead, we start with
the expression
$$\multline
\underset{\ep_2\to0}\to{\lim_{\ep_1\to0}}\frac{
\ka\,    {\left(  \ka-1 \right) }^{(t-y)/{2}}\,
 ({ \textstyle n}) _{(t-y)/2}\,(2n+t-1)_{(y-t)/2} }
{\big(\frac{t-y}{2}\big)! \,
    \big(\frac{t+y}{2}+n-1\big)! }\\
\times
    {} _{4} F _{3} \!\left [ \matrix { 1 - 2 n - t, \frac{3}{2} - n - 
\frac{t}{2},  1 - n - \frac{t}{2} - 
\frac{y}{2}+\ep_2,-\frac{t}{2} + \frac{y}{2}+\ep_1}
\\ { \frac{1}{2} - n - \frac{t}{2},  1 - n - \frac{t}{2} +
\frac{y}{2}-\ep_2, 2 - 2 n - \frac{t}{2} - 
\frac{y}{2}-\ep_1}\endmatrix ; {\displaystyle 
\frac{1}{1 - \ka}}\right ]\\
=
\underset{\ep_2\to0}\to{\lim_{\ep_1\to0}}\frac{
\ka\,    {\left(  \ka-1 \right) }^{(t-y)/{2}}\,
 ({ \textstyle n}) _{(t-y)/2}\,(2n+t-1)_{(y-t)/2} }
{\big(\frac{t-y}{2}\big)! \,
    \big(\frac{t+y}{2}+n-1\big)! }\kern4cm\\
\times
 \sum _{k=0} ^{\infty} 
\frac {(1 - 2n - 
t+2k)} {(1 -2 n - t)}\frac {  (1 - 2 n - t)_k}
 {k!}\kern5cm\\
\cdot
\frac {(  1 - n - \frac{t}{2} - 
\frac{y}{2}+\ep_2)_k\,(-\frac{t}{2} + \frac{y}{2}+\ep_1)_k} 
{(  1 - n - \frac{t}{2} +
\frac{y}{2}-\ep_2)_k\,( 2 - 2 n - \frac{t}{2} - 
\frac{y}{2}-\ep_1)_k}
\(\frac{1}{1 - \ka}\)^k.
\endmultline\tag\Zt$$
We split the sum over $k$ into the three ranges $0\le k\le \frac {t-y} {2}$,
$\frac {t-y} {2}+1\le k\le \frac {t+y} {2}+2n-2$, and $\frac {t+y} 
{2}+2n-1\le k$, and then we perform the
limit $\ep_1\to0$. Thereby each summand in the range $\frac {t-y}
{2}+1\le k\le \frac {t+y} {2}+2n-2$
vanishes identically. Thus, upon writing the sum over the range 
$\frac {t+y} {2}+2n-1\le
k$ in hypergeometric notation, the expression (\Zt) becomes, after
also performing the limit $\ep_2\to0$,
$$\multline 
\frac{
\ka\,    {\left(  \ka-1 \right) }^{(t-y)/{2}}\,
 ({ \textstyle n}) _{(t-y)/2}\,(2n+t-1)_{(y-t)/2} }
{\big(\frac{t-y}{2}\big)! \,
    \big(\frac{t+y}{2}+n-1\big)! }
 \sum _{k=0} ^{(t-y)/2} \frac { (1 - 2n - 
t+2k)}
 {(1 -2 n - t)}\\
\cdot
\frac {(1 - 2 n - t)_k\,(  1 - n - \frac{t}{2} - 
\frac{y}{2})_k\,(-\frac{t}{2} + \frac{y}{2})_k} 
{k!\,(  1 - n - \frac{t}{2} +
\frac{y}{2})_k\,( 2 - 2 n - \frac{t}{2} - 
\frac{y}{2})_k}
\(\frac{1}{1 - \ka}\)^k\\
- \frac {\ka} {{\left(\ka -1 \right) }^{y+2n-1}}
\frac{
      \big({ \textstyle  \frac{t-y}{2}+1}\big) _{y + 2 n -2} \,
      \big({ \textstyle \frac{t+y}{2} +2n-1}\big) 
_{1 - \frac{t}{2} + \frac{y}{2}} }
{\big( \frac{t+y}{2} +2n-1\big)! 
\,(y + n -1)! }  \kern3cm\\
\times
      {} _{4} F _{3} \!\left [ \matrix { -1 + 2 n + y, \frac{1}{2} + n + 
\frac{y}{2},  n,-\frac{t}{2} + \frac{y}{2}}\\ { - \frac{1}{2}
 + n + \frac{y}{2}, n+y,2 n + \frac{t}{2} + \frac{y}{2}}\endmatrix ; 
{\displaystyle \frac{1}{1-\ka} }\right ] 
\endmultline\tag\Zu$$
As it turns out, 
the first sum in this expression is exactly the sum in (\Zq) (which is
seen by replacing $h$ by $\frac {t-y} {2}-k$ in the sum in (\Zq)). 
Therefore, if
we equate (\Zt) and (\Zu), we see that the sum in (\Zq) is equal to
$$\multline 
 \frac {\ka} {{\left(\ka -1 \right) }^{y+2n-1}}
\frac{(y+2n-1)!}
{\big(\frac {t-y} {2}\big)!\,\big( \frac{t+y}{2} +2n-1\big)! 
\,(y + n -1)! }  \kern3cm\\
\times
      {} _{4} F _{3} \!\left [ \matrix { -1 + 2 n + y, \frac{1}{2} + n + 
\frac{y}{2},  n,-\frac{t}{2} + \frac{y}{2}}\\ { - \frac{1}{2}
 + n + \frac{y}{2}, n+y,2 n + \frac{t}{2} + \frac{y}{2}}\endmatrix ; 
{\displaystyle \frac{1}{1-\ka} }\right ] \\
+
\underset{\ep_2\to0}\to{\lim_{\ep_1\to0}}\frac{
\ka\,    {\left(  \ka-1 \right) }^{(t-y)/{2}}\,
 ({ \textstyle n}) _{(t-y)/2}\,(2n+t-1)_{(y-t)/2} }
{\big(\frac{t-y}{2}\big)! \,
    \big(\frac{t+y}{2}+n-1\big)! }\kern4cm\\
\times
    {} _{4} F _{3} \!\left [ \matrix { 1 - 2 n - t, \frac{3}{2} - n - 
\frac{t}{2},  1 - n - \frac{t}{2} - 
\frac{y}{2}+\ep_2,-\frac{t}{2} + \frac{y}{2}+\ep_1}
\\ { \frac{1}{2} - n - \frac{t}{2},  1 - n - \frac{t}{2} +
\frac{y}{2}-\ep_2, 2 - 2 n - \frac{t}{2} - 
\frac{y}{2}-\ep_1}\endmatrix ; {\displaystyle 
\frac{1}{1 - \ka}}\right ].
\endmultline$$
Now we apply the transformation formula (\Zs) to both $_4F_3$-series.
Subsequently we perform the limit in the second term, 
and we simplify the resulting
expressions. After this is substituted back in (\Zq),
this yields exactly the expression (\Zra).\quad \quad \qed
\enddemo

\subhead 7. Asymptotic formulas for the partition function 
for watermelons of arbitrary, but fixed deviation\endsubhead
In this section we embark on the asymptotic analysis of the partition
function $Z^{(n)}_{t}(y;\ka)$ for families of $n$ vicious walkers, the $i$-th 
starting at $(0,2i-2)$ and
ending at $(t,y+2i-2)$, $i=1,2,\dots,n$, 
none of them running below the $x$-axis. Our starting point is the
double sum formula for $Z^{(n)}_{t}(y;\ka)$ in Theorem~\TK. It enables
us to compute (in principle) full asymptotic expansions for all
$\ka>0$ (see the Remark after the
proof). The technique which we employ is singularity analysis. Our
analysis reveals a phase transition at $\ka=2$. We present the first
term in the corresponding asymptotic expansions in Theorem~\TE\
below. Previous asymptotic results for $y=0$
were found by Owczarek, Essam and Brak in \cite{\OwEBAA, Eq.~(4.58)},
who gave
predictions for the order of magnitude of the partition function 
$Z^{(n)}_{2r}(0;\ka)$ as $r$ tends to $\infty$. Our results confirm
their predictions, but make them at the same time more
precise, as we are also able to provide the multiplicative constants
and, as we already mentioned,
in fact, full asymptotic expansions.

\proclaim{Theorem \TE}Let $y$ be a fixed non-negative integer.
As $t$ tends to $\infty$, the partition function
$Z^{(n)}_{t}(y;\ka)$ is asymptotically
$$\multline
2^{nt}t^{-n(2n+1)/2}{2^{2n^2-\frac {n} {2}+1}\ka}\, 
{\pi^{-n/2}}\,(2n-1)!\prod _{i=0} ^{n-2}\frac {i!\,(y+n+i-1)!} {(y+2i)!}\\
\times
\bigg(\sum _{h=0} ^{n}\frac {\binom nh\binom {y+h-2}h} {\binom {2n}h}
\frac {1} {(2-\ka)^{2n-h}}\bigg)
\(1+O(t^{-1})\),
\quad \text {if }\ka<2,
\endmultline\tag\AB$$
it is
$$
2^{nt}t^{-n(2n-1)/2}2^{2n^2-\frac {3} {2}n+1}\pi^{-n/2}
\prod _{i=0} ^{n-1}\frac {i!\,(y+n+i-1)!} {(y+2i)!}
\dsize\(1+O(t^{-1})\),
\quad \text {if }\ka=2,
\tag\AC$$
and it is
$$\multline
\(\frac {2^{n-1}\ka} {\sqrt{\ka-1}}\)^tt^{-(n-1)(2n-1)/2}
\frac {2^{(n-1)(4n-5)/2}\ka(\ka-2)^{2n-1}}
{\pi^{(n-1)/2}(\ka-1)^{\frac {y} {2}+2n-1}}\\
\kern2cm \times
 \prod _{i=0} ^{n-2}\frac {i!\,(y+n+i-1)!} {(y+2i)!}
\(1+O(t^{-1})\),
\quad \text {if }\ka>2.
\endmultline\tag\AD$$
\endproclaim

\remark{Remark}
This result shows different asymptotic behaviour for $\ka<2$, for
$\ka=2$, and for $\ka>2$. For $\ka<2$ the order of magnitude of the
partition function is $2^{nt}t^{-n(2n+1)/2}$ (i.e., the
growth rate is $2^n$, and the critical exponent is $n(2n+1)/2$), 
for $\ka=2$ it is $2^{nt}t^{-n(2n-1)/2}$, and for $\ka>2$ it is
$\(\frac {2^{n-1}\ka} {\sqrt{\ka-1}}\)^tt^{-(n-1)(2n-1)/2}$ (i.e., now
the growth rate is not constant anymore, but grows with $\ka$),
everything else is constants.
\endremark

\demo{Proof}We want to determine the asymptotic behaviour of
(\Zo). Clearly, for the product in front of the summation we just
have to apply Stirling's formula. The result is that
$$\multline
\frac {(\frac {t+y}2-1)!}{(\frac {t+y}2+2n-2)!}
\prod _{i=0} ^{n-2}\frac {(t+2i)!\,i!\,(y+n+i-1)!}
{(\frac {t+y}2+n+i-1)!\,
   (\frac {t-y}2+i)!\,(y+2i)!}\\
\sim
2^{(n-1)t}t^{-((n+1)(2n-1))/2}2^{2n^2-\frac {5} {2}n+\frac {3}
{2}}\pi^{(1-n)/2} \\
\times
 \prod _{i=0} ^{n-2}\frac {i!\,(y+n+i-1)!} {(y+2i)!}
\(1+O(t^{-1})\)
\endmultline\tag\ZK$$
as $t$ tends to $\infty$.

From now on we concentrate on the double sum in (\Zo). Let us write
$t=2r+y$ in the sequel. Since $y$ is
fixed and the summation index $k$ comes from a bounded range,
we may use the estimation
$$\frac {\big(\frac {t+y}2\big)_k} {\big(\frac {t-y}2+n-k\big)_k}
=\frac {(r+y)_k} {(r+n-k)_k}=
1+O(r^{-1})\tag\ZKa$$
to approximate the double sum by
$$\multline
\Bigg(\sum _{\ell=0} ^{r}
\sum _{k=0} ^{n-1}(-1)^k
   \binom {n-1}{k}\binom {y+2k-2}{2k}
\binom {2r+y-\ell-1}{r-\ell-k}\hskip2cm\\
\cdot
   \binom {y+\ell+2n-2}{2n-2k-1}\,
 { (2n-2k-1)!\,(2k)!}
\,\ka^{\ell+1}\Bigg)\big(1+O(r^{-1})\big).
\endmultline$$

As is well-known, the most convenient method for determining the
asymptotic behaviour of a sequence of numbers is {\it singularity
analysis}
as developed by Flajolet and Odlyzko \cite{\FlOdAA} (see \cite{\FlSeAA}
for an introduction to that method). This method requires the
generating function for the sequence of numbers to be known 
as a starting point.
In our case, this means that we need a compact expression for
$$\multline
\sum _{r=0} ^{\infty}z^{2r+y}
\sum _{\ell=0} ^{r}
\sum _{k=0} ^{n-1}(-1)^k
   \binom {n-1}{k}\binom {y+2k-2}{2k}
\binom {2r+y-\ell-1}{r-\ell-k}\hskip2cm\\
\cdot
   \binom {y+\ell+2n-2}{2n-2k-1}\,
{ (2n-2k-1)!\,(2k)!}
\,\ka^{\ell+1}.
\endmultline\tag\ZH$$

We shall make use of the fact that
$$\sum _{r=0} ^{\infty}\binom {2r+m}rx^r=\frac {C(x)^{m+1}}
{1-x\,C(x)^2},\tag\ZHa$$
where $C(x)$ is the generating function (\ZI) for Catalan numbers.
This identity follows readily from (\ZS) by writing the left-hand side
of (\ZHa) as
$$\align
\sum _{r=0} ^{\infty}\binom {2r+m}r&x^r\\
&\kern-1.5cm=
\sum _{r=0} ^{\infty}\(\binom {2r+m}r-\binom {2r+m}{r-1}+
\binom {2r+m}{r-1}-\binom {2r+m}{r-2}+-\cdots\)
 x^r\\
&\kern-1.5cm
=C(x)^{m+1}+xC(x)^{m+3}+x^2C(x)^{m+5}+\cdots.
\endalign$$
If we interchange summations in (\ZH) so that the sum over $r$ becomes
the inner-most sum, and then apply (\ZHa), 
then we obtain the expression
$$\multline
\sum _{\ell=0} ^{\infty}
\sum _{k=0} ^{n-1}(-1)^k
   \binom {n-1}{k}\binom {y+2k-2}{2k}
   \binom {y+\ell+2n-2}{2n-2k-1}\\
\cdot
{ (2n-2k-1)!\,(2k)!}
\frac {z^{y+2k+2\ell}C(z^2)^{y+2k+\ell}} {1-z^2\,C(z^2)^2}\ka^{\ell+1}.
\endmultline$$
Now we write
$$ \binom {y+\ell+2n-2}{2n-2k-1}=\sum _{h=0} ^{2n-2k-1}\binom {\ell+h}h
\binom {y+2n-h-3}{2n-2k-h-1}\tag\ZHc$$
(which follows from the Chu--Vandermonde summation formula (\ZP)) and
interchange summations so that the sum over $\ell$ becomes the
inner-most sum. Using the binomial theorem and
(\ZWa), we end up with the expression
$$\multline
\sum _{k=0} ^{n-1}\sum _{h=0} ^{2n-2k-1}(-1)^k
   \binom {n-1}{k}\binom {2n-2k-1}h 
\binom {y+2n-h-3}{2n-h-1}\\
\cdot
h!\,(2n-h-1)!\,
\frac {2^{h+1}\ka} {(2-\ka)^{h+1}}
\frac {z^{y+2k}C(z^2)^{y+2k-1}} {\sqrt{1-4z^2}
\,\(1+\frac {\ka} {2-\ka}\sqrt{1-4z^2}\)^{h+1}}
\endmultline
\tag\ZJ$$
for the generating function (\ZH).

Singularity analysis now allows to extract the
asymptotic behaviour of the coefficient of $z^{2r+y}$ as $r$ tends to
$\infty$ out of the behaviour of
the function in $z$ at its singular points. 
We have to distinguish between three cases: $\ka<2$,
$\ka=2$, and $\ka>2$.

Let first $\ka<2$. Then, in the disk $\vv{z}\le1/2$, the only
singularity of the expression (\ZJ) is at $z=1/2$. 
For $z$ close to $1/2$ we have the singular expansion
$$
\frac {C(z^2)^{y+2k-1}} {\sqrt{1-4z^2}
\,\(1+\frac {\ka} {2-\ka}\sqrt{1-4z^2}\)^{h+1}}
\sim (1-4z^2)^{-1/2}2^{y+2k-1}+D
+O\((1-4z^2)^{1/2}\),
\tag\Zv$$
where $D$ is an explicit constant independent of $z$, whose value
is however without relevance here.
The coefficient of $x^m$ in $(1-4x)^{-1/2}$ behaves like
$(4^m/\sqrt{\pi m})(1+O(m^{-1}))$ as $m\to\infty$.
The transfer theorems \cite{\FlSeAA, Theorems~5.4 and 5.5} 
then imply that the
coefficient of $z^{2r-2k}$ in the left-hand side of (\Zv) (which is
what we need when we want to consider the coefficient of $z^{2r+y}$ in
(\ZJ)) behaves like
$(2^{2r+y-1}/\sqrt{\pi r})(1+O(r^{-1}))$
as $r\to\infty$. If we substitute this in (\ZJ),
then we obtain that the coefficient of $z^{2r+y}$ in (\ZJ) is
$$\multline
\sum _{k=0} ^{n-1}\sum _{h=0} ^{2n-2k-1}(-1)^k
   \binom {n-1}{k}\binom {2n-2k-1}h 
\binom {y+2n-h-3}{2n-h-1}\\
\cdot
h!\,(2n-h-1)!\,
\frac {2^{2r+y+h}\ka} {\sqrt{\pi r}\,(2-\ka)^{h+1}}\(1+O(r^{-1})\)
\endmultline
$$
as $r$ tends to infinity.
We interchange sums and write the (now) inner sum over $k$ in
hypergeometric notation. This gives
$$\multline
\sum _{h=0} ^{2n-1}
\binom {2n-1}h 
\binom {y+2n-h-3}{2n-h-1}\\
\cdot
h!\,(2n-h-1)!\,
\frac {2^{2r+y+h}\ka} {\sqrt{\pi r}\,(2-\ka)^{h+1}}
{}_2F_1\left[\matrix \frac {1} {2}+\frac {h} {2}-n,1+\frac {h} {2}-n\\
\frac {1} {2}-n \endmatrix;1\right]\(1+O(r^{-1})\).
\endmultline
$$
Clearly, the $_2F_1$-series can be evaluated by means of the
Chu--Vandermonde summation (\ZP). As it turns out, it is only non-zero
for $h\ge n-1$. After substitution of the result of the evaluation
in the above expression,
after replacement of $h$ by $2n-h-1$, after taking into account that
$t=2r+y$, and finally combining the result
with (\ZK), the claimed result (\AB)
follows upon little rearrangement.

If $\ka=2$, then it suffices to apply Stirling's formula to the
closed form expression (\ZM).

Finally let $\ka>2$. Then, in the disk $\vv{z}\le 1/2$, the
``dominant" singularity (i.e., the singularity with least modulus;
cf\. \cite{\FlSeAA}) is $z=\sqrt{\ka-1}/\ka$. For $z\to \sqrt{\ka-1}/\ka$
we have 
$$\align
\kern-1.5pt\frac {C(z^2)^{y+2k-1}}
{\sqrt{1-4z^2}\(1+\frac {\ka} {2-\ka}\sqrt{1-4z^2}\)^{h+1}}
&\sim
\frac {\(\frac {\ka} {\ka-1}\)^{y+2k-1}}
{\frac {\ka-2}\ka\(-\frac {2(\ka-1)} {(\ka-2)^2}\)^{h+1}
\(1-z^2\frac {\ka^2} {\ka-1}\)^{h+1}}\\
&\sim
\frac {(-1)^{h+1}\ka^{y+2k}(\ka-2)^{2h+1}}
{2^{h+1}\,(\ka-1)^{y+2k+h}
\(1-z^2\frac {\ka^2} {\ka-1}\)^{h+1}}.
\tag\Zw\endalign$$
The transfer theorems \cite{\FlSeAA, Theorems~5.2, 5.4 and 5.5} 
imply that the
coefficient of $z^{2r-2k}$ in the left-hand side of (\Zw) behaves like
$$
(-1)^{h+1}\frac {\ka^{y+2k}(\ka-2)^{2h+1}}
{2^{h+1}\,(\ka-1)^{y+2k+h}}
\(\frac {\ka^2} {\ka-1}\)^{r-k}\frac {r^h} {h!}
\(1+O(r^{-1})\).$$
Hence, the asymptotically dominating terms result from the summands in
(\ZJ) where $h$ is maximal, i.e., where $h=2n-1$, which, because of
the binomial coefficient $\binom {2n-2k-1}h$ occurring in (\ZJ), 
in turn forces $k$ to be zero.
If this is combined with (\ZK), and if we again take into account that
$t=2r+y$, then the claimed result (\AD)
follows upon little simplification.

This completes the proof of the theorem.\quad \quad \qed
\enddemo

\remark{Remark} Since Stirling's formula (see \cite{\WhWaAA, Sec.~12.33})
does in fact provide a full asymptotic expansion for factorials, as
well as do the transfer theorems \cite{\FlSeAA, Theorems~5.2, 5.4 and 5.5} 
provide full
asymptotic expansions, the above approach
does in fact allow to compute full asymptotic expansions for the
partition function $Z^{(n)}_{2r}(y;\ka)$ in all three different
regions, if needed. The reader should note that, if we want to compute
more terms in the asymptotic expansion, we will have to extend
the asymptotic expansion (\ZKa). The most convenient way to do that (for the
subsequent generating function computations) is in the form
$$\multline
\frac {\big(\frac {t+y}2\big)_k} {\big(\frac {t-y}2+n-k\big)_k}
=\frac {(r+y)_k} {(r+n-k)_k}\\=
1+\frac {k(y-n+k)} {r+y+k}+
\frac {C_2} {\(r+y+k\)_2}+\cdots+
\frac {C_N} {\(r+y+k\)_{N}}+ O\(r^{-N-1}\),
\endmultline\tag\ZKb$$
with the $C_i$'s appropriate constants independent of $r$.
\endremark

\subhead 8. Exact and asymptotic results for the number of contacts
of $n$ vicious walkers\endsubhead
By definition, the {\it mean number of contacts} with the wall for families of 
$n$ vicious walkers, the $i$-th starting at $(0,a_i)$ and
ending at $(t,e_i)$, none of them running below the $x$-axis (the wall), is
$$M^{(n)}_{t}(\bold a\to\bold e;\ka):=
\sum _{\ell\ge1} ^{}N^{(n)}_t(\bold a\to\bold e;\ell)\,\ell\,\ka^{\ell},$$
where $N^{(n)}_t(\bold a\to\bold e;\ell)$ 
is the number of these families of vicious walkers {\it with exactly
$\ell$ contacts} with the wall. Clearly, we have
$$M^{(n)}_{t}(\bold a\to\bold e;\ka)=\ka\frac {d}
{d\ka}Z^{(n)}_{t}(\bold a\to\bold e;\ka),\tag\Zff$$
with $Z^{(n)}_t(\bold a\to\bold e;\ka)$ the partition function for
these vicious walkers as defined in the
Introduction.
In addition, the {\it normalized} mean number of contacts with the wall
is defined by
$$\frac {M^{(n)}_{t}(\bold
a\to\bold e;\ka)}
{Z^{(n)}_{t}(\bold a\to\bold e;\ka)}.\tag\Zf$$

In this section we shall again concentrate on 
watermelon configurations which start on the wall and have deviation
$y$, i.e., 
$n$ vicious walkers as above with $\bold a
=(0,2,4,\dots,2n-2)$ and $\bold e=(y,y+2,y+4,\dots,y+2n-2)$,
and analyse their (normalized) mean number of
contacts. Specializing (\Zf) to the above choices of $\bold a$ and 
$\bold e$, the quantity that we want to study is
$$\frac {\ka\frac {d} {d\ka}Z^{(n)}_{t}(y;\ka)}
{Z^{(n)}_{t}(y;\ka)},
\tag\ZO$$
where, as before, $Z^{(n)}_{t}(y;\ka)$ denotes the partition function for
these watermelon configurations.
Let us denote the normalized mean in (\ZO) by $\bar M^{(n)}_{t}(y;\ka)$.
We could use any of our formulas for the partition function
$Z^{(n)}_{t}(y;\ka)$ (such as (\Zo) or (\Zq)) to find explicit
representations of $\bar M^{(n)}_{t}(y;\ka)$ as a quotient of two
sums (double sums in the case of (\Zo)).

As was the case for the partition function, also these formulae for
the mean number of
contacts do not simplify in general, while they do
for $\ka=1$, and if $y=0$ also for $\ka=2$.
\proclaim{Theorem \TF}For $\ka=1$, the normalized mean number of contacts,
$\bar M^{(n)}_{t}(y;1)$, is equal to 
$$1+\frac {n(y+2n+1)(t-y)} {(y+n)(t+y+4n)}.$$ 
\endproclaim
\demo{Proof}
The value of the partition function at $\ka=1$, $Z^{(n)}_{t}(y;1)$,
was already determined in (\ZL). For computing the numerator in (\ZO),
we make use of the formula (\Zq). If this formula is differentiated
with respect to $\ka$, and subsequently $\ka$ is set equal to $1$,
then it is only two terms in the sum which survive. Some
simplification then leads to the claimed formula.\quad \quad \qed
\enddemo

Whereas there does not seem a ``nice" formula for the mean number of
contacts for $\ka=2$ for arbitrary $y$, there is one if $y=0$.

\proclaim{Theorem \TFa}For $y=0$ and $\ka=2$, 
the normalized mean number of contacts,
$\bar M^{(n)}_{2r}(0;2)$, is equal to 
$$\frac {2^{2r-1}n\binom {2n}n} {\binom {2r+2n-2}{r+n-1}}-2(n-1).$$
\endproclaim
\demo{Proof} 
Using the
representation of the partition function $Z^{(n)}_{2r}(0;\ka)=
Z^{(n)}_{2r}(\ka)$ as given in (\ZG), the normalized mean number of
contacts is
$$\bar M^{(n)}_{2r}(0;2)=2+
\frac {\sum\limits _{\ell=0} ^{r-1}\binom
{2r-\ell-2}{r-1}\binom {\ell+2n-1}\ell \,\ell\,2^{\ell+2}} 
{\sum\limits _{\ell=0} ^{r-1}\binom
{2r-\ell-2}{r-1}\binom {\ell+2n-1}\ell2^{\ell+2}}.$$
A comparison of (\ZG) and (\ZM) shows that
we already evaluated (implicitly) the denominator of this expression.
Hence,
the evaluation of the numerator of (\ZO) will be accomplished once we
have evaluated the sum in the numerator.
We begin by converting it
to hypergeometric notation, and obtain
$${{16\,n \,
     ({ \textstyle r-1}) _{r-1} }\over {({r-1})! }}
{} _{2} F _{1} \!\left [ \matrix { 1 + 2 n, 2 - r}\\ { 3 -
      2 r}\endmatrix ; {\displaystyle 2}\right ].
$$
Next we reverse the order of summation and obtain
$${{{2^{r+1}}\,
     ({ \textstyle 2 n}) _{r-1} }\over {({r-2})!
}}
{} _{2} F _{1} \!\left [ \matrix { 2 - r, r}\\ { 2 - 2 n -
      r}\endmatrix ; {\displaystyle {1\over 2}}\right ]
$$
for the numerator of (\ZO). 
We now apply the contiguous
relation
$${} _{2} F _{1} \!\left [ \matrix { a, c}\\ { c}\endmatrix ;
   {\displaystyle z}\right ]  =
  {{ c - 1    }\over {a - 1}}
  {} _{2} F _{1} \!\left [ \matrix { a - 1, c}\\
        { c - 1}\endmatrix ; {\displaystyle z}\right ] 
    + {{ a - c   }\over {a - 1}
     }
  {} _{2} F _{1} \!\left [ \matrix { a - 1,
        c}\\ { c}\endmatrix ; {\displaystyle z}\right ].
$$
Thus, the $_2F_1$-series becomes the sum of the two terms
$$    {{{2^{r+1}} \,
      ({ \textstyle 2 n}) _{r} }\over {({r-1})! }}
{} _{2} F _{1} \!\left [ \matrix { 1 - r, r}\\ { 1 - 2 n -
       r}\endmatrix ; {\displaystyle {1\over 2}}\right ]
-{{{2^{r+2}}\,n \,
       ({ \textstyle 2 n}) _{r-1} }\over {({r-1})! }} 
{} _{2} F _{1} \!\left [ \matrix { 1 - r, r}\\ { 2 - 2 n -
        r}\endmatrix ; {\displaystyle {1\over 2}}\right ].
$$
Both of the $_2F_1$-series can be summed by means of Bailey's
summation formula (cf\. \cite{\SlatAC, (1.7.1.8); Appendix (III.7)})
$$
{} _{2} F _{1} \!\left [ \matrix { a, 1 - a}\\ { b}\endmatrix ; {\displaystyle
   {1\over 2}}\right ]  =\frac { \Gamma \( {b\over 2}\)\, \Gamma \( {1\over 2} +
   {b\over 2}\)} { \Gamma \( {a\over 2} + {b\over 2}\)\, \Gamma \( {1\over 2} - {a\over 2} + {b\over
   2}\)}.
\tag\ZOa$$ 
Some simplification then yields the claimed
formula.\quad \quad \qed
\enddemo

We now embark on determining the asymptotic behaviour of the
normalized mean number of contacts. Predictions on the order of
magnitude in the case of two walkers have been already made earlier by
Brak, Essam and Owczarek in \cite{\BrEOAC, Sec.~4.4.2}. 
The next theorem solves the problem for an arbitrary number of
walkers, thereby confirming these predictions and, at the same time,
making them more precise. Again, our method
allows in fact to obtain full asymptotic expansions (see the 
Remark after the
proof of Theorem~\TE).

\proclaim{Theorem \TG}As $t$ tends to $\infty$, the normalized 
mean number of contacts\linebreak 
$\bar M^{(n)}_{t}(y;\ka)$ is asymptotically
$$1+C+O(t^{-1}),\quad \text {if }\ka<2,\tag\ZQa$$
where $C$ is the quotient
$$\ka\bigg(\sum _{h=0} ^{n}\frac {(2n-h)\binom nh\binom {y+h-2}h} {\binom {2n}h}
\frac {1} {(2-\ka)^{2n-h+1}}\bigg)\Bigg/
\bigg(\sum _{h=0} ^{n}\frac {\binom nh\binom {y+h-2}h} {\binom {2n}h}
\frac {1} {(2-\ka)^{2n-h}}\bigg),\tag\ZQb$$
it is
$$
2^{\frac {1} {2}-2n}n\binom {2n}n\sqrt{\pi
t}-2n-y+2+O\big(t^{-1/2}\big),\quad \text {if }\ka=2,
\tag\ZQc$$
and it is
$$
\frac {(\ka-2)} {(\ka-1)}\frac {t} {2}+
\frac {\ka((2-\ka) y+4n-2)} {2(\ka-1)(\ka-2)}+1+
O\(t^{-1}\),\quad \text {if }\ka>2.
$$
\endproclaim

\remark{Remark}
Thus, for the length of the walks being large, 
the normalized mean number of contacts is
proportional to a constant if $\ka<2$, it is proportional to the
square root of the length of the walks if $\ka=2$, and it is
proportional to the length of the walks if $\ka>2$.
\endremark

\demo{Proof}We start by observing that Theorem~\TK\ implies that
the normalized mean number of
contacts $\bar M^{(n)}_{t}(y;\ka)$ can be rewritten as
$$\bar M^{(n)}_{t}(y;\ka)=1+
\frac {U^{(n)}_{t}(y;\ka)} {V^{(n)}_{t}(y;\ka)},
\tag\ZR$$
where $V^{(n)}_{t}(y;\ka)$ is the double sum in (\Zo), and where
$$\multline 
U^{(n)}_{t}(y;\ka)=
 \sum _{\ell=0} ^{(t-y)/2}
\sum _{k=0} ^{n-1}(-1)^k
   \binom {n-1}{k}\binom {y+2k-2}{2k}
\binom {t-\ell-1}{\frac {t-y}2-\ell-k}\hskip2cm\\
\cdot
   \binom {y+\ell+2n-2}{2n-2k-1}
\frac {\big(\frac {t+y}2\big)_k\, (2n-2k-1)!\,(2k)!}
{\big(\frac {t-y}2+n-k\big)_k}\ell\,\ka^{\ell+1}.
\endmultline$$

We already determined the asymptotic behaviour of
the denominator of the fraction in 
(\ZR) in the proof of Theorem~\TE. What we still
need is the asymptotic behaviour of the numerator. We attack this
problem in the same way as the corresponding problem for the
denominator. Here we will also have to use singularity analysis in
the case that $\ka=2$, because there is no closed form result
available for the normalized mean number of contacts 
in that case (as opposed to for the partition function, which allowed us
to do the asymptotics for $\ka=2$ in Theorem~\TE\ by just making use
of Stirling's formula). 

Since for the cases $\ka=2$ and $\ka>2$ we aim at
computing an additional term (beyond the leading term) in the
asymptotic expansion, we will have to use the
estimation (\ZKb), with $N=1$, as a first step. 
After the substitution of this in the numerator of the fraction
in (\ZR),
we compute the generating function with the coefficient
of $z^t$ being the sum in the numerator. 
To be precise, again writing $t=2r+y$,
we want to find a compact expression for
$$\multline
\sum _{r=0} ^{\infty}z^{2r+y}
\sum _{\ell=0} ^{r}
\sum _{k=0} ^{n-1}(-1)^k
   \binom {n-1}{k}\binom {y+2k-2}{2k}
\binom {2r+y-\ell-1}{r-\ell-k}\hskip2cm\\
\cdot
   \binom {y+\ell+2n-2}{2n-2k-1}
{ (2n-2k-1)!\,(2k)!}\(1+\frac {k(y-n+k)} {r+y+k}\)
\,\ell\,\ka^{\ell+1}.
\endmultline\tag\Zx$$

Now we perform the same steps as in the proof of Theorem~\TE: we
interchange summations so that the sum over $r$ becomes the inner-most sum, 
then sum over $r$ using
(\ZHa) (respectively (\ZS) for the additional term), 
next use (\ZHc), then interchange summations so that the sum over
$\ell$ becomes the inner sum,
and finally evaluate the remaining sum over $\ell$ by means of
the binomial theorem. The result is that our generating function
(\Zx) turns into
$$\multline
\sum _{k=0} ^{n-1}\sum _{h=0} ^{2n-2k-1}(-1)^k
   \binom {n-1}{k}\binom {2n-2k-1}h 
\binom {y+2n-h-3}{2n-h-1}\\
\cdot
(h+1)!\,(2n-h-1)!\,
\frac {2^{h+2}\ka^2} {(2-\ka)^{h+2}}
\frac {z^{y+2k+2}C(z^2)^{y+2k}} {\sqrt{1-4z^2}
\,\(1+\frac {\ka} {2-\ka}\sqrt{1-4z^2}\)^{h+2}},
\endmultline
$$
plus a similar term resulting from the fraction 
${k(y-n+k)} /({r+y+k})$ occurring in (\Zx).

Now we apply singularity analysis to the function on the right-hand
side. The considerations are entirely analogous to the corresponding
ones in the proof of Theorem~\TE, except that here we have to consider
one more term in the singular expansion.
As it turns out in the end, the additional term 
${k(y-n+k)} /({r+y+k})$ does in fact
produce smaller terms as we need. We leave the details
to the reader. 

If, finally, the obtained results are combined with the corresponding
ones for the sum in the denominator of the fraction on right-hand
side of (\ZR) which were found in the proof of Theorem~\TE, then
after little simplification the claimed asymptotic expressions are
obtained.\quad \quad \qed
\enddemo

\remark{Remark}
The sums which appear in (\ZQa) in the case that $\ka<2$ can be
expressed in terms of hypergeometric $_2F_1$-series. (The sum in the
denominator is a $_2F_1$-series, while the sum in the numerator can be
written as a sum of two $_2F_1$-series.) From there we see that they
cannot be simplified, with the exception of a few special cases where
$_2F_1$-summation formulas are available (such as, for example, 
(\ZP) and (\ZOa)). For instance, if $\ka=1$, then the Chu--Vandermonde
summation formula (\ZP) is applicable. However, in that case it would
have been easier to do the asymptotics directly from the closed form
expression in Theorem~\TF.
\endremark

\enddocument